\documentstyle[12pt,epsfig,color,amssymb,latexsym]{article}

 \newtheorem{thm}{Theorem}[section]
 \newtheorem{defn}{Definition}[section]
 \newtheorem{lemma}[thm]{Lemma}
 
 \newtheorem{prop}[thm]{Proposition}
 \newtheorem{cor}[thm]{Corollary}
 \newtheorem{ex}[thm]{Example}

 \newtheorem{question}[thm]{Question}
 
 \newtheorem{remark}{Remark}[section]
 \newtheorem{observation}[thm]{Observation}

 \newcommand{\beq}[1]{\begin{equation}\label{#1}}
 \newcommand{\enq}[0]{\end{equation}}

 \newcommand{\ceiling}[1]{\lceil{#1}\rceil}
 \newcommand{\qed}[0]{\begin{flushright} $\Box $ \end{flushright}}
 \newcommand{\C}[2]{{{#1}\choose{{#2}}}}
 \newcommand{\ga}[0]{\alpha }
 \newcommand{\gb}[0]{\beta }
 
 \newcommand{\gd}[0]{\delta }

 \newcommand{\gs}[0]{\sigma}

 \newcommand{\eps}[0]{\varepsilon }
 \def\ep{\epsilon}
 \def\ulam{\Lambda}

 \newcommand{\ra}[0]{\rightarrow}

 \newcommand{\Ee}[0]{\mbox{${\bf E}$}}

 \newcommand{\cee}[0]{{\cal C}}

 \newcommand{\F}[0]{{\cal F}}
 \newcommand{\G}[0]{{\cal G}}
 \newcommand{\h}[0]{{\cal H}}

 \newcommand{\R}[0]{{\cal R}}
 \newcommand{\T}[0]{{\cal T}}
 
 \newcommand{\X}[0]{{\cal X}}

 \newcommand{\nin}[0]{\noindent}

 \renewcommand{\dots}[0]{,\ldots,}

 \def\ep{\epsilon} \def\lam{\lambda}
 
 \def\calo{{\cal O}}
 \def\calp{{\cal P}}

 \title{A Sharp Threshold for Random Graphs with a Monochromatic Triangle
 in Every Edge Coloring}  \author {Ehud Friedgut \and Vojtech R\"odl
 \and Andrzej Ruci\'{n}ski
 \and Prasad Tetali}

\begin {document}

\footnotetext{ Szemer\'edi's
Regularity Lemma. \\ Research supported by grants of the ISF, BSF, NSF and KBN.
For details see section 7.1}
 \maketitle
 \begin{abstract}
 Let $\R$ be the set of all finite graphs $G$ with the  Ramsey property that  every coloring of the edges of $G$
 by two colors yields a monochromatic triangle.  In this paper we establish a sharp threshold for random graphs
 with this property.
 Let $G(n,p)$ be the random graph on $n$ vertices with edge probability $p$.
 We prove that  there exists  a function $\widehat c=\widehat c(n)=\Theta(1)$ such that
 for any $\eps > 0$, as $n$ tends to infinity,

$$Pr\left[G(n,(1-\eps)\widehat c/\sqrt{n}) \in \R \right] \rightarrow 0$$
and
$$Pr \left[ G(n,(1+\eps)\widehat c/\sqrt{n}) \in \R\ \right] \rightarrow 1.$$
\nopagebreak A crucial tool that is used in the proof and is of independent
interest is a generalization of Szemer\'edi's Regularity Lemma to
a certain hypergraph  setting.
  \end{abstract}
  \newpage
  \tableofcontents
\newpage
 \section{Introduction}\label{intro}
 \subsection{Overview}\label{over}
 This paper brings together several important themes of
 combinatorics: Ramsey properties, threshold phenomena of random graphs, and Szemer\'edi-type regularity.

 Ramsey properties guarantee, for an arbitrary partition of a given structure,
 that a highly organized substructure can be found in at least one part of the partition.
 During the last decade of the last century there has been extensive study of Ramsey properties
 of random structures, see e.g. \cite{FR,LRV,RR1,RR2,RR3,Rado,Schur,RR4}.
 These papers were all concerned with establishing a threshold function for various Ramsey-type properties,
 either of random graphs, random hypergraphs or random sets of integers.
 For a binomial random graph $G(n,p(n))$,  for instance, they provide
 a critical edge probability $\widehat p=\widehat p(n)$ such that the limiting probability that
 every coloring of a random graph $G(n,p)$ contains certain monochromatic structures
 depends on the asymptotic ratio between $p$ and~$\widehat p$.

In all the above papers there is  a multiplicative gap left between the upper and lower bound on the threshold
edge probability $\widehat p$ (see Theorem  \ref{RRedge} below). It is therefore not surprising that the  natural
question of whether the gap can be closed has been around for some time. In other words, does there exist a
constant $\hat c$ such that the asymptotic probability that $G(n,c\widehat p)$ has a Ramsey property is either 0
or 1, depending only on whether $c > \widehat c$ or $c<\widehat c$? This question is usually phrased in the
specialized jargon as ``does there exist a {\em sharp threshold}?"

 Sharp thresholds have been established for many random graph properties,
like  connectivity, hamiltonicity and perfect matchings.
However, such precise results for Ramsey properties seemed out of hand until 1999,
when a general technique for settling these questions was introduced in \cite{F}.
Loosely speaking, the main theorem in \cite{F}
showed that the question of sharpness of threshold for a random
graph property is determined by whether the property
is related to local  or rather global graph phenomena.

 Two papers exploiting the technique from \cite{F} for coloring questions are \cite{AF} and \cite{FK}.
 The latter paper (as well as the earlier \cite{RR3}) states as the next natural candidate
 for attack, the problem of sharpness of the threshold for  property $\R$ consisting of all graphs $G$ such that in every blue-red coloring of the edges of $G$ there exists
 a monochromatic triangle.
 However, this problem turned out to be more difficult than those in \cite{AF} and \cite{FK}
 and required some new combinatorial approach.

 In this paper we add the missing tool that enables us to crack this nut.
 It is a regularity lemma for a certain class of hypergraphs, whose edges consist of small subgraphs
 of a fixed underlying sparse graph (see Lemma \ref{reglem}).
 Our lemma is a generalization of the celebrated Szemer\'edi Re\-gu\-larity Lemma for graphs \cite{Szem},
 and follows in the footsteps of  the regularity lemma for sparse graphs  presented in \cite{Kohayakawa97Szemeredi} (see also \cite{KR} and \cite{JLR}, Section 8.3)
 and of the hypergraph regularity lemma of Frankl and R\"odl \cite{FR1}.

The proof of our regularity lemma, Lemma \ref{reglem}, provides for a considerable portion of the bulk of the proof
of the following sharp threshold theorem, which is the ultimate result of our paper.

\begin {thm}
 \label{main} There exists a function $\widehat c=\widehat c(n)=\Theta(1)$ such that for every  $\eps>0$ \beq {e4} \lim _{n \ra \infty}{\rm Pr}[G(n,p)\in\R] = \left\{\begin{array}{ll}
      1& \mbox{if $p>(1+\eps)\widehat c/\sqrt{n}$}\\
      0 & \mbox {if $p<(1-\eps)\widehat c/\sqrt{n}$}.  \end{array}  \right.  \enq  \end{thm}
There is a slightly disappointing aspect of this result: although we prove that $\widehat c(n)$ is bounded,  the
natural conjecture is that $\widehat c(n)$ converges to a positive limit, and this does
not follow from our theorem. Unfortunately, an inherent property of the technique we use is that it can only
supply such existence theorems but no new information as to the exact threshold probability. We discuss various
possible extensions of Theorem \ref{main} in Section \ref{final}.

 \subsection{Ramsey Properties of Random Graphs}\label{ramsey}

Let us introduce the arrow notation, commonly used in Ramsey theory. For two graphs, $H$ and $G$, and an integer
$r\ge2$,  we write $G \ra (H)_r$  if for every coloring of the edges of $G$ by $r$ colors there exists a
monochromatic copy of~$H$. For example, it is well known that $K_6\ra(K_3)_2$. Let $\R$ be the set of all graphs
$G$ such that $G\ra(K_3)_2$.

A basic question studied in Ramsey theory is, given a graph $H$ and an integer $r\ge2$, when is $G$ ``rich'' enough
for $G \ra (H)_r$? Here richness can be interpreted either as the number of edges
of $G$, or as the ratio of edges to vertices.

In modern graph theory, problems of this type are often studied via random graphs.
The theory of random graphs addresses questions concerning typical graphs, or graphs ``on
average".
The standard model for a random graph is $G(n,p)$, a graph on $n$ vertices,
where every one of the $\C {n}{2}$ edges of the complete graph belongs to
$G(n,p)$ independently, with probability $p$.
When studying random graphs a natural problem is: Given $H$, find
a threshold function $\hat p(n)$ such that $G(n,p)   \ra (H)_r$
 with probability tending to 1
when $\hat p(n) \gg p$ and $G(n,p)  \ra (H)_r$ with probability tending to 0
when $ \hat p(n) \ll p$.
(The existence of such a threshold function is guaranteed by
a general result of Bollob\'as and Thomason  \cite{BT}.)

  In a series of papers \cite{FR,LRV,RR1,RR2,RR3}, a threshold function $\widehat p(n)$
is determined for all graphs $H$. Its culmination, paper  \cite{RR3}, establishes $\widehat p=n^{-1/m^{(2)}(H)}$
as a threshold for $G(n,p)\ra (H)_r$, regardless of $r$, for all $H$ which are not star forests. Here
$m^{(2)}(H)=\max_{F\subseteq H}(|E(F)|-1)/(|V(F)|-2)$. An analogous result in the case when the vertices (and not
the edges) are colored is given in \cite{LRV}.
In the edge--coloring setting, the first case to be settled was that of triangles
(not counting star forests which are rather trivial).

\begin{thm}[\cite{RR2}]
 \label{RRedge}
 For every integer $r\ge2$ there exist constants $c_r$ and $C_r$ such that
 \beq
 {e2}
 \lim _{n \ra \infty}{\rm Pr}[G(n,p)\ra (K_3)_r] =
 \left\{\begin{array}{ll}
      1& \mbox{if $p>C_r/\sqrt{n}$}\\
      0 & \mbox {if $p<c_r/\sqrt{n}$.}
 \end{array}
 \right.
 \enq
 \end{thm}
{\bf Remark:}  In the above theorem, and throughout the paper $G(n,p)$ is usually meant to denote
$G(n,p(n))$, hence the limits, and asymptotic notations.
As we can see, for a range of $p$, namely for $c_r/\sqrt{n}\le p\le C_r/\sqrt n$ this statement is inconclusive.
Similarly, in all the papers mentioned above there is  a multiplicative gap
left between the upper and lower bound on the threshold edge probability~$\widehat p(n)$.

In a recent paper \cite{FK} it was shown that in many cases this gap can be closed, using a general technique from
\cite{F} for proving sharpness of thresholds.
 The cases treated in \cite{FK}
cover vertex--coloring when $H$, the graph defining the Ramsey property, belongs to  a wide family of graphs
including, for example, all complete graphs. Also the case of edge--coloring when $H$ is a tree is dealt with. In
all these instances it is shown that there exists a function $\widehat p(n)$ such that for every $\eps>0$, $\lim
_{n \ra \infty}{\rm Pr}[G(n,p)\ra (H)] =1$, if $p>(1+\eps)\widehat p(n)$ and 0, if $p<(1-\eps)\widehat p(n)$.

  It is worthwhile pointing out here a difference between this
kind of a ``sharp threshold'' statement and the previous Ramsey threshold results. Although the results from
\cite{FK} show that the transition from the non-Ramsey region of the values of $p$ to the Ramsey region is swift
and sharper than what was previously proven, these are only an existence results.
They give no further information on
the critical value of $p$ which was calculated in previous works.

 The most basic case not treated in \cite{FK} is the case of graphs with a monochromatic triangle
in every edge coloring, for which a (weak) threshold is given by Theorem \ref{RRedge}.
In the present paper we prove a sharp threshold theorem
for this Ramsey property, with two colors. This is our
Theorem \ref{main} stated above.

 \subsection{Sharp Thresholds of Increasing Graph Properties}
 \label {sharp}
 In this subsection we introduce the notion of a {\em sharp threshold}
 for a graph property, as well as a technique for proving sharpness
 of thresholds.
 Consider a property $Q$ of graphs on $n$ vertices.
 We  identify $Q$ with the set of graphs having this property,
 and use the notation $G \in Q$ to denote the fact that
 $G$ has property $Q$.
 We will restrict ourselves to properties which are invariant under a graph
 automorphism and also distinguish an important class of {\it increasing}  graph
 properties, i.e. those that are preserved under edge addition.

 \begin{defn}\rm We say that $\widehat p=\widehat p(n)$ is a
{\it  threshold function} for an increasing graph property $Q$  if
 $$\lim_{n\to\infty} {\rm Pr}[G(n,p)\in Q]=\cases{1 & if $\widehat p=o(p)$\cr
 0 & if $p=o(\widehat p)$.}$$\end{defn}
 Bollob\'as and Thomason proved in \cite {BT} the existence of threshold
 functions for  all increasing set properties, and in particular
 for all graph properties.

Some properties do not have sharper thresholds in the sense that for all $p=p(n)$ which are of the same order as
$\widehat p$, we have $0<\lim_{n\to\infty}{\rm Pr}[G(n,p)\in Q]<1$. E.g., this is the case of the (increasing)
property of containing a copy of a given, balanced graph $H$, the threshold for which has been established by
Erd\H os and R\'enyi \cite{ER60} at $n^{-1/\rho(H)}$.
Here $\rho(H)$ is the edge to vertex ratio in $H$.
For example, $\rho(K_k)=(k-1)/2$, so the thresholds for
the appearance in $G(n,p)$ of a copy of $K_3$, $K_5$, and $K_6$, respectively, are $n^{-1}$, $n^{-1/2}$ and
$n^{-2/5}$.

But there are other properties, like connectivity, which enjoy much sharper  thresholds. Indeed, it has been
proved by Erd\H os and R\'enyi in \cite{ER59} that
$$\lim_{n\to\infty} {\rm Pr}[G(n,p) \mbox{ is connected}]=\cases{1 & if $np-\log n\to \infty$\cr
 0 & if $np-\log n\to -\infty$.}$$

 \begin{defn}\rm We say that an increasing property
 $Q$ has a {\em sharp threshold} if
 there exists a function  $\widehat p= \widehat p(n)$, such that for every $\eps >0$:
 $$\lim_{n\to\infty} {\rm Pr}[G(n,p)\in Q]=\cases{1 & if $p>(1+\eps)\widehat p$\cr
 0 & if $p<(1-\eps)\widehat p$.}$$
 Otherwise we say that $Q$ has a {\em coarse threshold}.
 \end{defn}
Thus, the property of being connected has a sharp threshold at $\widehat p=\log n/n$. Our main theorem in this
paper, Theorem \ref{main}, states that
 the Ramsey property $\R$ has a sharp threshold.

 In \cite {F} the first author gives a necessary
 and sufficient condition for an increasing property to have a sharp threshold.
 This condition will be the central tool used in this paper.
 Roughly stated, it says that a property with a
 sharp threshold is not statistically determined
 by a simple local reason, that is, by the presence or absence of finitely many edges. For example,
 the property of having a triangle as a subgraph is obviously local,
 and indeed has a coarse threshold (both the critical probability and the
 length of the threshold interval equal $\Theta(1/n)$), whereas it seems
 obvious that connectivity is a non-local property, even though it is statistically equivalent to the absence of isolated vertices.

For the Ramsey property $\R$, the condition in \cite {F} asserts that in order to establish the sharpness of its
threshold, one has to show that $\R$ is not influenced by the appearance of any fixed subgraph, which is likely to
be contained in $G(n,p)$, with the range of $p=p(n)$ limited by Theorem \ref{RRedge} to $p=\Theta(1/\sqrt n)$. Of
course, $\R$ is extremely influenced by the appearance of $K_6$, which, however, is very unlikely to be present in
such a $G(n,p)$.

 In this paper we will use a version of the sharpness criterion from \cite{F},
 which follows readily from the original statement
 but is more suitable for applications.
 Given a  graph $M$ and a disjoint set of $n$
 vertices, let $M^*$ be  an ordered copy of $M$ placed uniformly at random
 on one of the $n!/(n-|V(M)|)!$ possible locations.
 \begin{thm}
\label{coarse} Let $Q$ be an increasing graph property, with a coarse
threshold. Then there exist real constants $0<c<C$, $\gb>0$, a rational $\rho$,
and a sequence $p=p(n)$ satisfying $cn^{-1/\rho} < p(n) < C n^{-1/\rho}$, such that $\gb<
Pr[G(n,p)\in Q]< 1-\gb$ for all $n$. Furthermore, there exist
$\ga, \xi > 0$ and a balanced graph $M$ with density $\rho$ for which
the following holds:
\\For every graph property $\G$ such that $G(n,p) \in \G$ a.a.s.
there are infinitely many values of $n$ for which there exists a graph
$G$ on $n$ vertices for which the following holds:
\begin{enumerate}
\item[(i)] $G \in \G$ ,
\item[(ii)] $G \not \in Q$ ,
\item[(iii)] $Pr[(G\cup M^*) \in Q] > 2\ga$,
\item[(iv)] $Pr[(G \cup G(n,\xi p)) \in Q] < \ga$ .
\end{enumerate}
\end{thm}
What the theorem says, intuitively, is that in the case of a
coarse threshold one can find two graphs, $G$ and $M$ as follows: $G$ is a fixed graph on $n$ vertices
that is not a random graph but rather a pseudorandom graph, typical of $G(n,p)$
(actually a random choice of $G \in G(n,p) \setminus Q$
will work with probability close to 1); $M$ is a ``magical'' balanced graph such that
it is often the case that adding a random copy of $M$ to $G$
induces the property in question, whereas increasing the number of edges in $G$
randomly by a constant proportion $\xi$ does not induce the
property. The addition of a copy  of $M$ corresponds roughly to inducing a local
property, in contrast to increasing the number of edges which corresponds roughly to
increasing the global density of a random graph.
Therefore the conclusion
of the hypotheses of the theorem is that the property $Q$ is ``statistically local''.

The typical way in which this theorem is used to prove that a property $Q$
has a sharp threshold involves two steps:
\begin {itemize}
\item The first step is usually easy. For a coarse property $Q$, the theorem guarantees
the existence of $M$. A possible explanation to
this would be that $M$ itself has the property
$Q$. In that case, since $Q$ is a monotone increasing property it would no
longer seem ``magical'' that adding a copy of $M$ induces property
$Q$. In other words, if $M \in Q$ then Assumption (iii) in the theorem is a triviality
and does not enable one to deduce anything.
Therefore, as a start, one has to rule out this possibility by
showing that a balanced graph with the prescribed density can not
have the property.
\item The second step is typically more involved. One chooses an appropriate
property $\G$ that is typical of $G(n,p)$ and then shows that a graph $G \in \G$ with
$Pr[(G \cup M^*) \in Q] > 2\ga$ is  quite
``saturated'', i.e. is  close to having property $Q$. Therefore
adding a random copy of $G(n,\xi p)$ should induce the property $Q$
with probability  much larger than $\ga$, contradicting
condition (4) of the theorem.
\end{itemize}

 In \cite{F}, \cite{AF}, and \cite{FK} one can see variations of this scheme used
 to prove that various graph properties have a sharp threshold. Although the
 basic technique is similar, each property presents its own
 difficulties and requires a special approach. The case of property $\R$,
 handled in this paper, has by far been the most difficult and involved;
 a key technique in our approach turns out to be {\it regularization}.

 \subsection{Regularity}\label{Regular}

One of the fundamental tools in asymptotic graph theory is the well-known regularity lemma of
Szemer\'edi~\cite{Szem} (see also~\cite{RS78}).  Indeed, since its discovery in the~70s, this lemma has been
instrumental in the study of the structure of large graphs.  The reader is referred to the excellent
survey~\cite{KomlosSimonovits96} for a thorough introduction to the wide range of applications of this result.

In essence, the regularity lemma tells us that any large graph may be decomposed into a bounded number of
quasi-random, induced bipartite graphs. Thus, this lemma is a powerful tool for detecting and making transparent
the \textit{random-like behavior} of large \textit{deterministic} graphs. What makes the lemma such a powerful
tool is  that it reveals a quasi-random structure that enables one  to carry out a deep quantitative analysis.

The precise formulation of the regularity lemma is somewhat technical (see Section \ref{dense} for details).  In
this short section, we only discuss some points in broad terms.

The quasi-random bipartite graphs that Szemer\'edi's lemma uses in its decomposition are graphs in which the edges
are uniformly distributed.  The uniformity is measured by the ratio of edges to potential edges (pairs), and so
this concept becomes trivial for graphs of vanishing density.  To manage sparse graphs, one may adjust the notion
of quasi-randomness by a natural rescaling, and it is a routine matter to check that the original proof extends to
this notion, provided we restrict ourselves to graphs of vanishing density that do not contain `dense patches'
(see Section \ref{srg}).  However, the quasi-random structure that the lemma reveals in this case
is harder to exploit than in the dense case, and one needs to work harder when applying the lemma to such `sparse graphs'.
 Nevertheless, there have been some successful applications of the lemma in this context
(see~\cite{Kohayakawa97Szemeredi,kohayakawa01}).

The idea of regularity has also been extended to uniform hypergraphs.  The version that is most relevant to us is
the one in Frankl and R\"odl~\cite{FR1}, which makes it possible to decompose triple systems into quasi-random
structures made up of triples together with an `underlying' quasi-random tripartite graph. In that setting, the
density is measured by the ratio of triples to the triangles in the underlying graph (see Section \ref{fra}).
Moreover, the concept of quasi-randomness here is strong enough to allow one to prove that these quasi-random
pieces contain the same number of finite substructures as they would had they
been truly random pieces (see~\cite{FR1} and
Nagle and R\"odl~\cite{nagle:_regul}).

In this paper, we shall introduce yet another concept of
regularity, which  expands and melts together the notions of
sparse graph regularity and hypergraph regularity. In the usual
context of graph regularity, and in some more delicate versions of
it there is an invisible underlying graph behind the graph we look
at, and the regularity expresses the distribution of specified
edges with respect to the edges of the underlying graph.
Similarly, a 3-uniform hypergraph can be viewed as a distinguished
collection of triangles among all triangles of an underlying
3-partite graph. Here, we shall be interested in investigating the
structure of sparse graphs with respect to some other fixed family
of small subgraphs. Viewing these subgraphs as edges of a
hypergraph, the lemma we prove (Lemma \ref{reglem}) may be
interpreted as a sparse hypergraphs version of the regularity lemma . Our
approach is partly based on methods from
\cite{Kohayakawa97Szemeredi,kohayakawa01} and \cite{FR1}, but it
faces a further difficulty:
 the assumption of `no dense patches' in the
standard case (see~\cite{Kohayakawa97Szemeredi,kohayakawa01}) was an easy consequence of properties of random
graphs and therefore did not play any significant role; the proof of the analogous fact in the setting of
this paper, however, requires a fairly complex argument (see Section \ref{nodense}).

 \subsection{Outline of the Paper}
 \begin{itemize}
 \item In Section \ref{proof} we present the skeleton  of the proof of the main
 theorem. It is actually a formal proof which is made extremely short and compact
 by relying on lemmas which will be proven in the rest of the paper. An overview of the forthcoming proofs and an (oversimplified) illustration will also be given there.
 \item In Section \ref{TandC}, assuming  the hypothesis  of Theorem \ref{coarse}, we construct a family $\mathcal{S}$ of special
 subgraphs of $G$ (called special constellations) and show how every triangle-free coloring of $G$
 defines a set of edges of $G$ that intersects all special constellations (that is, defines a hitting set for the family $\mathcal{S}$.
 \item In Section \ref{regular} we prove a regularity
 theorem in the spirit of Szemer\'edi's Regularity Lemma that
 provides a partition of both the vertices and the edges of $G$ such that
  the special constellations are uniformly distributed with respect to this partition.
 \item In Section \ref{cores}, based on the regular partition found in the previous section,  we define a core
 which is a central notion of this proof, and show some crucial  properties of cores.
 \item In Section \ref{randomgraphs} we show various properties of
 random graphs that are needed throughout the paper.
It is there where the family $\G$ is defined, and an important lemma, Lemma \ref{pro1}, is proved.
 \item We conclude with open questions and possible
 extensions.
 \item At the end of the paper we include a glossary of symbols and definitions as well as a flowchart of constants exhibiting their mutual dependencies.
 We strongly encourage the reader to make use of both when struggling through our proof.
 \end{itemize}

 \subsubsection*{ Notation:}
 In  Sections \ref{regular}-\ref{randomgraphs} we use the following notation.
 For $0< \eps < 1$, and positive reals $x,y$,
 \[ x \stackrel{\eps}{\sim} y \ \  \mbox{  denotes that  } \ \
 (1- \eps) y \le x \le {(1+ \eps)} y.\]
 We will often abbreviate it further as follows:
 if $\eps'$ is any function of $\eps$ that tends to zero with
 $\eps$, and $x \stackrel {\eps'}{\sim} y$, then we will simply write
 $x \stackrel {\eps}\approx y$.

 Let $G=(V,E)$ be a graph, $v,u \in V$
  and $W \subseteq V$. We write $e_G(W)=|E(G[W])|$ for the number of edges in the subgraph of $G$ induced by $W$,
 ${\rm deg}_G(v,W)={\rm deg}(v,W)$ for the number of neighbors of $v$ in $W$,
${\rm deg}_G(v)={\rm deg}(v)={\rm deg}(v,V),$ for the degree of $v$ in $G$, and
 ${\rm codeg}(v,u)$ for the number of common neighbors of $u$ and $v$ in $G$,
  called the {\it co-degree} of $v$ and $u$.
The set of neighbors of $v$ in $G$ is denoted by $N_G(v)=N(v)$, while $N_G(W)$
 stands for the set of vertices outside  $W$, each having at least one neighbor in $W$
(so called {\it open} neighborhood of $W$).

For a family of sets $(A_i)_{i\in I}$, we call a set $T$ {\it a hitting set} if $T\cap A_i\neq\emptyset$ for all $i\in I$.
We will often use set partitions
 $V=V_1\cup...\cup V_t$, where
$|V_1|\le...\le |V_t|\le|V_1|+1$. Such partitions will be called here {\it equipartitions}.
Finally, all logarithms are natural and will be denoted by $\log $.

 \section{Outline of the Proof}
 \label{proof}
 \subsection{Main steps}
 Recall that $\R$ is the graph property that for every blue-red coloring of the
 edges of a graph there exists a
 monochromatic triangle. Graphs that have this property will be
 called Ramsey graphs. For a non-Ramsey graph, we call a coloring that does not have
 a monochromatic triangle {\em a triangle-free coloring}.

 We wish to prove that $\R$ has a sharp threshold.
 By Theorem \ref{RRedge}, there exist constants $c_2$ and $C_2$ such that any threshold $\widehat p$  for $\R$ satisfies
 $$c_2/\sqrt{n} < \widehat p < C_2/\sqrt{n}.$$
 This means that when applying
 Theorem \ref{coarse} to Property $\R$ we may restrict ourselves to sequences $p=p(n)$ falling into this range, and consequently
 to balanced graphs $M$ with  $\rho(M)=2$ (i.e. average degree 4.)
 Thus it suffices to prove the following result, which is a mere adaptation of  Theorem \ref{coarse} to our case.
(In fact, it is slightly stronger, since the inequality in (iv) is replaced in (\ref{conclusion}) below by
convergence to 1.)
 \begin{thm}
 \label{allweneed}
 For all $\ga,\xi>0$, all sequences $$c_2/\sqrt{n} < p=p(n) < C_2/\sqrt{n}$$
 and  all balanced graphs $M$ with  $\rho(M)=2$,
 there exists a graph property $\G$ with $\lim_{n\to\infty}{\rm Pr}[G(n,p) \in \G]=1$, and an integer $n_1$, such that for all
 $G \in \G \setminus \R$ with $|V(G)|=n > n_1$, if
 \beq
 {assumption}{\rm Pr}[G \cup M^* \in \R] > 2\ga
 \enq
  then
 \beq
 {conclusion}{\rm Pr}[(G \cup G(n,\xi p )) \in \R] = 1-o(1).
 \enq
 \end{thm}
Note that if the assumption $G \not \in \R$ (Assumption (ii) in Theorem \ref{coarse}) were false
then (\ref{conclusion}) would be trivial.
 A finer point is that if $M \in \R$
 then (\ref{assumption}) yields no information on $G$, and hence
 is useless. Fortunately the following lemma rules out this possibility.
 This is the typical ``easy step" mentioned in the introduction.
 \begin{lemma}
 \label{color} If $M$ is balanced and $\rho(M)=2$ then $M \not
 \in \R.$
 \end {lemma}
 \nin{\bf Proof:}
 It is a well known fact from the theory of random graphs (see
 \cite{JLR}, page 66) that for any {\em balanced} graph $H$ with  $\rho(H)=\rho$
 and $p=\Theta(n^{-1/\rho})$  there exists a constant $\gb>0$ such that the
 probability that $H$ appears in $G(n,p)$ is at least $\gb$.
 Hence, for any $b>0$, the probability that $M$ appears as a subgraph of
 $G(n,b/\sqrt{n})$ is bounded away from zero. If $M \in \R$ then,
 by the monotonicity of $\R$, this would mean that ${\rm Pr}[
 G(n,b/\sqrt{n}) \in \R]$ is also bounded away from zero.
 But by Theorem \ref{RRedge}, for all $b<c_2$,  ${\rm Pr}[G(n,b/\sqrt{n}) \in \R] \ra 0$. Therefore $M \not
 \in \R$.

 Alternatively, and without the use of random graphs,
 Lemma \ref{color} follows from a result in ~\cite{KuRu} which shows
 that any Ramsey graph (for $K_3$) must have a subgraph $H$ for which $\rho(H) \geq 5/2$. \hfill $\Box$
\bigskip

The rest of the paper is devoted to proving that (\ref{assumption}) implies (\ref{conclusion}).
An approach to statements like (4), which has become standard by now,
is via the so called two round exposure.
This technique originated in the seventies, in work of Posa, Ajtai, Komlos, Szemeredi,
Fernandez de la Vega, Fenner and Frieze, devoted to the existence of Hamilton paths
and cycles in random graphs (see \cite{BB85}, Chapter VIII, for references).
In the context of Ramsey properties, it was explored already in [29] and [30]
(see also [15], Sections 1.1 and 8.4).
For a graph $G=(V,E)$ and a set of edges $F \subseteq E$, let
 ${\rm Base}(F)$\label{base} be the set of edges in the complete graph on $V$
 that form a triangle with two edges of $F$, formally:
 $$ {\rm Base}(F) = \{uv :  \ uw, wv \in F \mbox{ for some }w\in V\} .$$
We often identify a subset of edges of a graph with the spanning subgraph consisting of them. So, in the above, both $F$ and ${\rm Base}(F)$
can be viewed as graphs on the same vertex set~$V$.

 Suppose we want to show that $G_1\cup G_2\in\R$ with high probability, where $G_i=G(n,b_i/\sqrt n)$, $i=1,2$, and the two random graphs are independent. (This is, in fact, our case with $b_2=\xi b_1$, except for two minor points: that
$G_1$ is not random but pseudorandom, and that $b_1$ may depend on $n$.)

First generate  the edges of $G_1$, and let them be colored by an
 adversary. Suppose that at least half of them are blue and call the set of
 blue edges $Blue$.    Clearly, if
 ${\rm Base}(Blue)$ contains a triangle $\Delta$, no proper
 coloring of $G_1 \cup \Delta$ can extend the adversarial coloring.
Therefore, it will suffice to show
 that ${\rm Base}(Blue) \cap G_2$ contains a triangle. (See Figure 1.)


\begin{figure}[hbt]
\begin{center}
      \centerline{\epsfig{figure=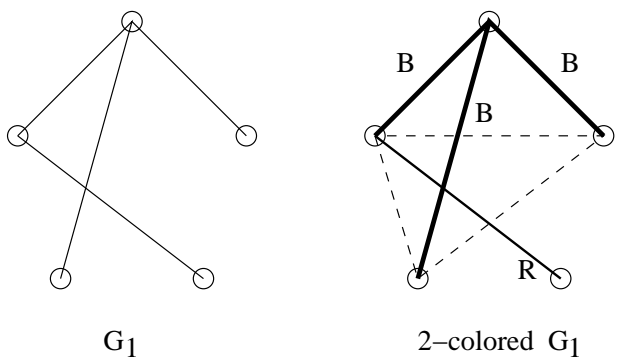}}
      \caption{No triangle-free coloring of $G_1 \cup \Delta$ can
      extend the given one}
\end{center}
\end{figure}

To this end we utilize the following lemma, proved in Section 6.
Given two real numbers $0<\lambda<1$ and $0<a<1/6$, we say that a graph $G$ has property $\T(\lambda,a)$, if for
any subgraph $F$ of $G$ with at least $\lambda |E(G)|$ edges, the graph ${\rm Base}(F)$ contains at least $a |V(G)|^3$
triangles.

\begin{lemma}\label{pro1}
For all $\lambda>0$ and $c>0$, there exists $a>0$, such that if $p\ge c/\sqrt n$ then,
 with probability $1-o(1)$, the random graph $G(n,p)$ has property
 $\T(\lambda,a)$.
\end{lemma}

Applying the above lemma with $F=Blue$, $\lambda=1/2$ and $c=b_1$ would yield $\Theta(n^3)$ triangles in
${\rm Base}(Blue)$. In the second round, by a standard application of a
correlation estimate from \cite{JLR90} (see\cite{JLR}, Section 2.2), often called Janson's inequality, at least one of these triangles will be included in
$G_2$ with very high probability. If we were allowed to take $b_2$ sufficiently large, then we could make the reciprocal of the error
probability larger than the exponential number of all bi-colorings of $G_1$, proving $G_1\cup G_2\in\R$ with
probability $1-o(1)$.

Unfortunately, in our case $b_2=\xi b_1$ depends on a given a priori $\xi$ and can be much smaller than $b_1$.
This  major difficulty demonstrates a more general problem of establishing a sharp threshold
{\em without knowing up front what the exact value of the critical constant $\hat c$ should be.}

Hence we must seek a refinement of the above approach, making use of the assumption (\ref{assumption}). And indeed,
through a special regularization of $G$ we will be able to construct a family $\rm CORE$ of subgraphs of $G$ such that
for every triangle-free coloring of $G$ at least one of these subgraphs is monochromatic. Moreover, the size of each
subgraph $K\in \rm CORE$ will be large enough to yield, via Lemma \ref{pro1} and Janson's inequality, at least one
triangle in ${\rm Base}(K)\cap G(n,\xi p)$ with probability very close to one, but at the same time the size of this
family, $|\rm CORE|$,  multiplied by the error probability will tend to 0. This is the content of the following two
lemmas which together imply Theorem~\ref{allweneed}. They are preceded by a {\em setup}, which we will often be
referring to in the paper.

\subsubsection*{Setup:}\label{setup}
 For the rest of the paper, let us fix constants $\ga, \xi>0$, a sequence $c/\sqrt{n} < p=p(n) < C/\sqrt{n}$, where $c=c_2$ and $C=C_2$, and a graph $M\not\in\R$
 as in Theorem \ref{allweneed} (it is no longer relevant that $M$ is balanced).
 We will define a graph property $\G$ in Definition \ref{defG} so that, in particular,  each $G\in\G$ has Property $\T(\lambda,a)$ with  $\lambda=\lambda(\alpha,c,C,M)$ to be specified later (or see the Glossary now) and  $a=a(\lambda,c)$  determined by Lemma \ref{pro1}.

We do not attempt to compute explicitly the integer $n_1$, promised in Theorem \ref{allweneed} and appearing in
the next lemma. In principle, it is the maximum of all values of $n_0$  encountered throughout the proof, most
notably the $n_0$'s in Theorem \ref{SSRL} and Lemma \ref{reglem}, as well as of several implicit lower  bounds on
$n$ hidden in our calculations.

 \begin{lemma}
 \label{wrapup}
 Let $G \in \G\setminus\R$ be a graph with $|V(G)|=n\ge n_1$ for which the assumption (\ref{assumption}) of Theorem \ref{allweneed} holds. Then for every $\tau>0$ there exists a family $\rm CORE$ of subgraphs of $G$
   such that
 \begin{enumerate}
 \item[(a)] For every triangle-free coloring $\chi$ of $E(G)$, there exists $K \in \rm CORE$
 which is monochromatic under $\chi$.
 \item[(b)] $|\rm CORE| \leq \exp(\tau n^{3/2})$,
 \item[(c)] For every $K \in \rm CORE$ we have $|K|>\lambda |E(G)|$.
 \end {enumerate}
 \end{lemma}
After thorough preparations in Sections \ref{TandC} and \ref{regular}, the family $\rm CORE$ will be constructed in
Section \ref{cores}. Also there, all three conclusions of the above lemma will be proven in the following manner:
 \begin{itemize}
 \item[(a)] follows from Lemma \ref{TransInEC} and Lemma \ref{CoreInColoring},
 \item[(b)] is Lemma \ref{SmallCore}, and
 \item[(c)] follows from Lemma \ref{TransInEC}, Lemma \ref{LargeCore} and Property (P3) of $\G$.
 \end{itemize}

A specific value of $\tau$, with which we apply Lemma \ref{wrapup}, is provided by another
application of Janson's inequality. We say
that  a subgraph $F \subseteq G$ {\it survives}
 if $F \cup G(n,\xi p )$ has a triangle-free coloring
 in which $F$ is monochromatic.
\begin{lemma}\label{survival}
For every $K \in \rm CORE$, the probability
 that $K$ survives is
 at most $\exp \left\{-2\tau_0 n^{3/2}\right\}$, where
$$\tau_0=\frac{a^2 \xi^6c^6}{2(a\xi^3c^3 + 2 \xi^5 C^5)}.$$
\end{lemma}
\nin{\bf Proof:}   By Lemma \ref{wrapup}(c) and the fact that $G\in\T(\lambda,a)$ there are at least $an^3$
triangles in ${\rm Base}(K)$. Let us number some $an^3$ of them by $1,2\dots an^3$ and let $I_i$ be the indicator random
variable for the event that the $i$th triangle is contained in $G(n,\xi p)$. Clearly, if $K$ survives then
$\sum_iI_i=0$. But, by \cite{JLR}, Theorem 2.18(ii),
$${\rm Pr}\left(\sum_iI_i=0\right)\le\exp\left\{-\frac{(\sum_i\Ee I_i)^2}{\sum\sum \Ee(I_iI_j)+\sum_i\Ee I_i}\right\}\ ,$$
where the double sum is over all ordered  pairs of distinct, edge-intersecting triangles in ${\rm Base}(K)$. Note that
$\sum_i\Ee I_i= an^3(\xi p)^3\ge a\xi^3c^3n^{3/2}$ and that the double sum contains at most $n^4$ summands, each
equalling at most $(\xi p)^5\le(\xi C)^5n^{-5/2}$. This completes the proof. \hfill $\Box$
\bigskip

 Let us now show how  Lemmas \ref{wrapup} and \ref{survival} imply Theorem
 \ref{allweneed}, and consequently Theorem \ref{main}.
 We refer to the members of $\rm CORE$ by the name {\it cores}.

\bigskip

\noindent{\bf Proof of Theorem \ref{allweneed}:}
  Lemma \ref{wrapup}(a)
 implies that if $G \cup G(n, \xi p)$ has a proper
 coloring $\chi$ (which, of course, induces a triangle-free coloring of $G$)
 then there exists a core which survives.
 But, using Lemma \ref{wrapup}(b) with $\tau=\tau_0$, Lemma \ref{survival} and a simple union bound, we deduce that with probability $1-o(1)$
 no core survives. Indeed,
 $${\rm Pr}[\mbox{Any core survives}] \le \mathop{\exp\{\tau_0 n^{3/2}\}}\limits_{\mathop{\uparrow}\limits_{\#\ of\ cores}} \cdot\mathop{\exp\{-2\tau_0 n^{3/2}}\limits_{\mathop{\uparrow}\limits_{\rm survival\ probability}} \} = o(1).$$
 Hence ${\rm Pr}[G \cup G(n, \xi p) \in \R] =1 -o(1)$ as required.
 \hfill $\Box$
 \subsection{Overview of the Proof Strategy }\label{sneak}
 Clearly, at such an early point in the paper, the  main idea of the
 proof may yet be obscure.
 Specifically, we have given no hint as to the connection
 between the existence of a magical graph $M$ having the property described in
 assumption (\ref{assumption}) of Theorem~\ref{allweneed}
 and the existence of a family $\rm CORE$ satisfying the three conclusions of
 Lemma~\ref{wrapup}. In order to shed some light on this connection
 (albeit it will still be a dim light), we give here
 a short explanation of the logic and motivation behind the construction of $\rm CORE$.
 The next three sections of the paper are devoted to the constructions that underlie the following scheme:
\begin{itemize}
\item The existence of $M$  such that
$$
Pr[G \cup M^* \in \R] > 2\ga
$$
implies that the set of triangle-free colorings of $G$ is very restricted;
there are many sets of vertices in $G$ such that planting
a copy of $M$ on them kills all triangle-free colorings, i.e. no triangle-free coloring of $G$ can be extended to $G \cup M$
when $M$ is placed on one of the aforementioned ``bad" sets.
\item We will associate every such ``bad" set with a set of edges in $G$,
a union of stars which we will call a special constellation.
We will fix a coloring of these special constellations in such a way that every proper
coloring of $G$ must agree with every such colored constellation on at
least one edge (see Lemma \ref{harm}).
\item The above can be translated to the language of hypergraphs:
we will construct a hypergraph whose hyperedges are the edge sets of
the special constellations, and it will turn out that every triangle-free coloring gives rise to a
hitting set of the hyperedges of this hypergraph (see Lemma \ref{TransInEC}).
\item We will show that every such hitting set (and hence every
triangle-free coloring) may be associated with a large monochromatic set called a core.
$\rm CORE$ is the family of cores (see Lemma \ref{wrapup} above).
\item The key to associating every hitting set with a core
is in showing that our hypergraph has an inherent regular structure
that may be revealed by a Szemer\'edi type partition (see Lemma \ref{reglem}).
\end{itemize}

 \subsubsection{An Illustration}\label{ill}

 To get a better feeling of how regularity  helps in creating the family $\rm CORE$,
 let us consider a simpler analogue that takes place in the well understood setting of graphs,
 in which  special constellations will be replaced simply by edges.
 Hopefully this will give some clue as to what we are doing
 in the forthcoming sections. We refer the reader to Section \ref{regular} for the notion of an $\eps$-regular
 graph and the statement of the Szemer\'edi Regularity Lemma.

 Let $H=H(n)$ be a sequence of graphs on $n$ vertices.
A {\em cover set}  in a graph is a set of vertices that intersects every edge.
In other words, it is a hitting set for the family of all edges of the graph.
In general, the number of cover sets
in an $n$-vertex graph may be exponential in $n$, and our goal is to ``capture" all cover sets
of $H$ by a smaller family of large subsets of vertices which we will call  {\it cores}.
(Elsewhere in the paper the term core is used in our larger setting, cf. Lemma \ref{wrapup}.)
We want the following properties to hold for cores:
 \begin {itemize}
\item Every cover set contains a core,
 \item The number of cores is $2^{o(n)}$,
 \item Every core is of size linear in $n$.
  \end{itemize}
In general, this is not an easy task,
as the last two properties  seem to contradict each other.

Before describing how to construct the cores, here is a hint as to {\em why} one would
want the last two conditions to hold simultaneously:
in our general setting we wish to capture the large family of all triangle-free colorings
of our graph $G$ by a smaller family of partial colorings, hence a {\em small} family of cores.
On the other hand, we wish every partial coloring to be sufficiently large to ensure that the probability
of being able to extend it to a larger graph $G \cup G(n, \xi p)$ is very small, hence
cores should be {\em large.}

 Returning to the graph $H$, suppose that $H$ is a  bipartite $\eps$-regular graph on vertex sets
 $V_1,V_2$, with $|V_1|=|V_2|=n$ and density $d(V_1,V_2)=d>0$.
 Recall that this means that for every $W_1 \subseteq V_1$ and $W_2 \subseteq V_2$
 such that $|W_1| > \eps n$ and $|W_2| > \eps n$, the density of the
subgraph between $W_1$ and $W_2$,
 is ``$\eps$-close" to $d$, and consequently,
 there is at least one edge between such $W_1$ and~$W_2$.
Hence any cover set must necessarily include at least $(1- \eps)n$ vertices from either $V_1$ or $V_2$. Then our
family of
 cores may be formed by all
 sets of size  $\lceil(1-\eps)n\rceil$ that
 are subsets of either $V_1$ or $V_2$.
 It is easy to see that every cover set contains a core, and the number of
 cores is $O(\C{n}{\eps n})$ which indeed is $2^{o(n)}$, for $\eps =
 o(1)$.

%
%

It is not  hard to generalize this construction
to the case of a multi-partite graph on vertex
 sets $V_1, \ldots , V_k$ such that for most pairs $V_i,V_j$ the
 spanned bipartite graph is dense and $\eps$-regular.   Now comes the use of the Szemer\'edi
 Regularity Lemma: since every sufficiently large and dense
 graph is very close to being of this type, the original
 problem can essentially be reduced to this case.

After seeing this example, hopefully, the reader may have a feeling as to why we will
eventually devote much energy to exposing the hidden regular structure in the hypergraph
which expresses the restrictions on the triangle-free colorings of $G$.

 \section{Tepees and Constellations}
 \label{TandC}
 Assumption (\ref{assumption}) of Theorem
 \ref{allweneed} should imply that $G$ is close to being Ramsey in the sense
 that its triangle-free colorings are quite restricted. In this section we set up a family of subgraphs
 of $G$ called special constellations that help to capture these restrictions. We will show that
 every triangle-free coloring of $G$ may be associated with a hitting set
 of this family.
 \subsection{Tepees}
 Our analysis of the restrictions imposed by (\ref{assumption}) on the colorings
 of $G$ will lead us rather naturally to the  structures
 we shall call  {\em tepees.}

 Assume that $M$, an arbitrary balanced graph with $\rho(M)=2$, has $\nu$ vertices, and thus $2\nu$ edges, and fix a generic copy of $M$ with vertices labeled by $x_1 \dots x_{\nu}$. We begin by defining a copy of $M$ ``planted" on an ordered subset of vertices of $G$.   For every sequence $X=(v_1 \dots v_{\nu})$ of distinct vertices of $G$, let $M_X$ be the copy of
 $M$ with vertex $x_a$ mapped onto $v_a$ for each $a=1\dots\nu$. We will often identify a sequence $X$ with the set $\{v_1\dots v_\nu\}$ of vertices making up $X$.

The family of all sequences $X$ satisfying $G \cup M_X \in \R$ will be denoted by~$\X$. Note that assumption
(\ref{assumption})
 of Theorem \ref{allweneed} implies that
 $$|\X| \ge 2 \ga n(n-1)\cdots (n-\nu+1) = (2-o(1))\ga n^{\nu}.$$
 Let $\X_1 \subseteq \X$ be the family of all $X\in\X$
 such that
 \begin{list}{}{}
 \item {(i)} $X$ is an independent set in $G$,
 i.e. the vertices in $X$ span no edges of $G$, and
 \item {(ii)} every vertex of $G$ has at most two neighbors in $X$.
\end{list}
Since $G \in \G$ (see Definition \ref{defG}, parts (P1) and (P2)), it follows that almost all  $X$'s have
properties (i) and (ii) above, and so we still have
$$|\X_1|  \ge (2-o(1))\ga n^{\nu}.$$
Because of (i), the property $G \cup M_X \in \R$ indicates that there should be some triangles in $G \cup M_X$
with one edge in $M_X$ and two edges in $G$ (seethe proof of Lemma \ref{harm} below). We now give a name to such
structures.

\begin{defn}\label{tepee}\rm For a pair of vertices $\{u,v\}$ of $G$
 a {\em tepee over $\{u,v\}$} is a pair of edges of $G$
 of the form $\{uw ,wv\}$. Then vertex $w$ is called {it the tip} of the tepee $\{u,v\}$.
\end{defn}
 In other words, a tepee over $\{u,v\}$ is any path of length two in $G$ with endpoints in $u$ and $v$. (Clearly, the vertices of a tepee over $\{u,v\}$  form a triangle in $G+uv$.)  We will denote such a tepee by $uwv$ for short.
 For a sequence of vertices $X\in {\cal X}_1$, let $T(X)$ denote the set of all tepees
 in $G$ over those pairs of vertices from $X$ which are edges of $M_X$. By properties (i) and (ii) above,
all  tepees in $T(X)$ are pairwise edge-disjoint and have distinct tips.

Furthermore, let $\widehat{X}$ be the graph formed
 by the edges of tepees in~$T(X)$. The vertex set of  $\widehat{X}$ consists of $X$ and he set of tips of all tepees in $T(X)$. (This notation is supposed to be
 suggestive of the tepees formed over~$X$.)
Set $t(X)=|T(X)|$. Given $t(X)=\phi$, $0\le \phi\le n-\nu$, the graph $\widehat{X}$ has $\nu+\phi$ vertices and
$2\phi$ edges, and is isomorphic to a graph which can be obtained from $M_X$ by replacing  some of its edges by multiple edges,
erasing others, and finally replacing all $\phi$ edges of the obtained multigraph by internally disjoint paths of length two. So,
there are many isomorphism types of $\widehat{X}$ possible. The next lemma states that a positive fraction of $X$'s
have the same isomorphism type of the graph $\widehat{X}$, and, moreover, the number of vertices in $\widehat{X}$
is bounded. (We
have sacrificed $\alpha_1$ for the sake of global harmony.)

\begin{lemma}\label{common} Let $q=\ceiling{10C^2\nu/\ga}$ and $\ga_2=\frac{\ga}{(2\nu+q)^q}$.
There exists an integer $\phi\le q$, a graph $\widehat M$ on  $\nu+\phi$ vertices, and a family $\X_2\subseteq
\X_1$ of size $|\X_2|=\alpha_2n^{\nu}$ such that for all $X\in{\cal X}_2$, the graph $\widehat{X}$ is isomorphic
to $\widehat M$.
\end{lemma}

\nin{\bf Proof:} We will first show that there are at most $\frac 4 5 \ga n^{\nu}$ sequences $X\in\X_1$  such that
$t(X)>q$.
 For every pair of vertices $\{u,v\} \subset V$,  let $T(u,v)$ be the set of
 tepees over $\{u,v\}$.

For a vertex $w \in V$ of degree ${\rm deg}(w)$
 there are exactly ${\rm deg}(w)\choose 2$ tepees of the form $uwv$.
This and the fact that $G \in \G$
 and, in particular, that the degrees of all vertices in $G$ are
 bounded from above by $2C\sqrt{n}$ (see Definition \ref{defG} (P3)), yields
 $$\sum_w {{\rm deg}(w) \choose 2} \leq n {{2C\sqrt{n}} \choose 2} < 2C^2n^2.$$
 To apply a simple counting argument note that
 for all $\{u,v\} \subset V$
 $$ |\{X: uv \in M_X\}| \leq   4\nu n^{\nu-2}.$$
 Consequently, by reversing the order of summation,
 $$\sum_{X\in\X_1} t(X) = \sum_{X\in\X_1} \sum _{uv \in  M_X} |T (u,v)|
 \leq 4\nu n^{\nu-2} \sum_{u,v} |T(u,v)|
 \leq 8C^2\nu n^{\nu}.$$
 This immediately implies that $t(X) > q$ holds for at most $\frac 4 5  \ga n^{\nu}$ sets $X\in\X_1$, and, therefore, at least
 $$|\X_1| - \frac45 \ga n^{\nu} >  \ga n^{\nu}$$
sets $X\in{\cal X}_1$ have $t(X)\le q$.

Given that $t(X)\le q$, the number of isomorphism types possible for the
 graph $\widehat{X}$ is no more than the number of
 ordered partitions of the integer $q$  into $2\nu+1$ nonnegative parts, corresponding to
 the decisions of how many  tepees span over each edge of $M_X$. (The $(2\nu+1)$-st part is the difference between $q$ and the actual number of tepees present.) There are ${{2\nu+q}\choose q}<(2\nu+q)^{q}$ such partitions.
 Take  as $\widehat M$ the
 most common isomorphism type among $\{\widehat{X}: X \in
 \X_1, t(X)\le q\}$ and set $\phi=|V(\widehat M)|-\nu$ and
 $$\X_2=\{X \in \X_1: \widehat{X}\mbox{ is isomorphic to } \widehat
 M\}.$$ Then
 $$|\X_2|\ge \frac{\ga}{(2\nu+q)^q} n^{\nu} = \ga_2 n^{\nu}.$$
\hfill $\Box$
\bigskip

 Let us summarize the properties of the family $\X_2$. For every $X \in \X_2$
 \begin {itemize}
 \item $G \cup M_X \in \R$,
 \item $X$ is an independent set,
 \item No three vertices in $X$ share a common neighbor in $G$,
 and hence $T(X)$ is composed of edge-disjoint tepees with disinct tips,
 \item $\widehat{X}$ is isomorphic to $\widehat M$, a graph with $\nu+\phi$ vertices and $2\phi$ edges, for some $\phi\le q$.
 \end{itemize}

 By Lemma \ref{color}, $M$ has a triangle-free coloring, that is, a blue-red coloring of its
 edges with no monochromatic triangle. Fix one such coloring and
 call it~$\gs'$. For each $X \in \X_2$, color $M_X$, the copy of $M$ planted on $X$, by  $\gs'$,  and denote the so colored copy by $M_X'$.
 This partitions the tepees in $T(X)$ into two sets,
 $RT(X)$ and $BT(X)$, as follows:  $RT(X)$  is the set
 of tepees in $G$ over the red edges of $M_X'$ and  $BT(X)$  is
 the set of tepees over the blue edges of $M'_X$.

 We broaden
 the definition of property $\R$ to partially colored graphs by saying that a  partially colored graph $F$ belongs to $\R$ if
 there is no triangle-free coloring of all edges of $F$ consistent
 with the partial coloring.
Clearly, if $F\in\R$ then $F$ with any partial coloring of $F$ belongs to $\R$ too. Hence, by the first property of $\X_2$ listed above,
for all $X \in \X_2$
 \beq
 {trivial}
 G \cup  M_X' \in \R.
 \enq
We now make a simple but crucial observation which captures  in terms of the tepees the restrictions on proper
colorings of $G$ imposed by (\ref{assumption}).
 \begin {lemma}
 \label{harm} For every triangle-free coloring $\chi$ of $E(G)$ and every $X \in
 \X_2$, there is either a tepee in $RT(X)$ colored red or
 a tepee in $BT(X)$ colored blue.
 \end{lemma}
 \nin{\bf Proof:} Let $X \in \X_2$. By (\ref{trivial}), for every triangle-free coloring $\chi$ of $E(G)$
 there is a monochromatic triangle in $G \cup  M_X'$. By the definition of $\gs'$, we know that $M_X'$ itself does not contain such a triangle.
 On the other hand, by the definition of $\X_2$,
 the vertices in $X$ do not span any edges of
 $G$.
 Therefore, we conclude that the monochromatic triangle contains exactly one edge of
 $M_X$. The other two edges form a tepee guaranteed by the lemma.
 \hfill $\Box$
\bigskip

We would now like to pre-color a subset of edges  of $G$ according to whether they belong to the tepees from
$RT(X)$ or $BT(X)$ for some $X\in \X_2$.  Based on the above lemma, we could then conclude that  every proper
coloring of $E(G)$ agrees with the pre-coloring on a set of edges which intersects every subgraph $\widehat X$ on
at least two edges, namely, the two edges making up a tepee.

 The problem with this approach is that we cannot assign a color to an edge $e$ in $G$
 once and for all, because $e$ may belong to a tepee in $RT(X_1)$ and
 at the same time to a tepee in~$BT(X_2)$ for some $X_1\neq X_2$.
 To remedy this obstacle  we now move to a more restricted family $\X_3 \subseteq \X_2$
 for which no such clashes occur.

Recall that $V(M)=\{x_1\dots x_{\nu}\}$. Let us label the remaining vertices of the prototype graph $\widehat M$
by  $u_1\dots u_{\phi}$. For each $X=(v_1\dots v_\nu)\in\X_2$ choose an isomorphism $f$ between $\widehat M$ and
$\widehat X$ which maps $x_a$ onto $v_a$ for all $a=1\dots\nu$,  and then set $w_b=f(u_b)$ for all $b=1\dots\phi$.

 Assuming that $n$ is divisible by $\nu+\phi$, let
 $\pi = \{V_1 \dots  V_{\nu}, W_{1} \dots W_{\phi}\}$ be a partition of $V=V(G)$ into $\nu+\phi$ parts
 of equal size $m=n/(\nu+\phi)$. (In reality, we put aside an arbitrary set of less than $\nu+\phi$
vertices, which will have only  a negligible effect on the  estimates
in our proofs.)

  We will call a subgraph $\widehat X$
  {\it consistent} with $\pi$ if, under the above  notational convention, $x_a \in V_a$ for $1
 \leq a \leq \nu$, and   $w_{b} \in W_{b}$ for $1 \leq b \leq \phi$. (See, e.g. Figure  2.)
Our next lemma establishes the existence of a partition with
respect to which a positive fraction of subgraphs  $\widehat X$ will be consistent. The
attached degree constraint will be utilized only in Sections \ref{regular} and \ref{cores}.

 \begin{lemma}
 \label{equi} There exists a partition $\pi$ as above
 and a family $\X_3 \subseteq \X_2$ with $|\X_3| = \alpha_3 n^{\nu}$,
 where $\alpha_3 = \alpha_2/2(\nu+q)^{(\nu+q)}$,
 such that for every $X \in \X_3$, the subgraph $\widehat X$
 is consistent with $\pi$. Moreover, for every vertex $v\in V_1\cup\cdots\cup V_\nu$
and for each $b=1\dots\phi$, we have
$\frac34mp \leq deg_G(v,W_b)\leq \frac32mp$, where $m=n/(\nu+\phi)$. \end {lemma}
 \nin{\bf Proof:}
Choose an ordered partition of the vertices of $V$
 into $\nu+\phi$ parts, uniformly at random from all such partitions.
 For a given subgraph $\widehat X$ the probability that it is
 consistent with the chosen partition is precisely
 $$\frac{m^{\nu+\phi}}{n(n-1)\ldots(n-\nu-\phi+1)} \geq (\nu+\phi)^{-(\nu+\phi)}.$$
 Therefore the expected number of subgraphs $\widehat X$ with $X \in \X_2$
 that are \emph{not} consistent with the random partition is at most
$$|\X_2|(1-(\nu+\phi)^{-(\nu+\phi)})\le (\ga_2-2\ga_3 )n^{\nu} .$$
 By Markov's inequality, with probability at least $\ga_3/(\ga_2-\ga_3)$ at least $\ga_3 n^{\nu}$ such subgraphs are consistent with the random partition.
On the other hand, Chernoff's inequality for hypergeometric distributions (see, e.g., \cite{JLR}, Theorem 2.10) and Property (P3) yield that the degree constraint is satisfied with probability $1-o(1)$. Hence, the existence of a required partition follows.
 \hfill $\Box$
\bigskip

 Let $\pi_0 = \{V_1 \dots  V_{\nu}, W_{1} \dots W_{\phi}\}$
  be a partition guaranteed by Lemma \ref{equi}
 and let $\X_3 \subseteq \X_2$ be the corresponding family.
 We have now finished refining the initial family $\X$.
 The new subfamily $\X_3$ is the final one we will work with.
 Also, the partition $\pi_0$ will now be fixed for the rest of the paper providing
 a starting point for our further constructions.

\begin{figure}[hbt]
\begin{center}
    \input{fig4_new.tex}
      \caption{subgraphs $\widehat X$ consistent with $\pi_0$}
\end{center}
\end{figure}

 \subsection{Constellations}\label{con}
Note that for all distinct $X_1,X_2\in{\cal X}_3$ we have $RT(X_1)\cap BT(X_2)=\emptyset$. In view of this and of
Lemma \ref{harm}, we are now in a position to assign a pre-coloring to the edges of $\bigcup_{X\in{\cal
X}_3}\widehat X$ in such a way that every triangle-free coloring of $G$ agrees with the pre-coloring on at least two
edges (which form a tepee) of $\widehat X$, for every $X\in{\cal X}_3$. However, for our purposes it is sufficient to
concentrate on such an agreement even on one edge of each $\widehat X$. Because of this excess, we now simplify the
structure of subgraphs we will be dealing with through the rest of the paper.

 For every $X \in \X_3$,  instead of $\widehat{X}$ we will consider a certain subgraph $S(X)$ of $\widehat X$,
 called a {\em special constellation}. (See Figure 3.)
 We will shortly define this new notion formally, but for now we note that
 $S(X)$ is a star forest (a disjoint union of stars) obtained by erasing one edge of every tepee
 in $T(X)$. The reason we now shift from tepees to constellations
 is that we will need various estimates on their number, and this boils down to counting star forests,
 a task much simpler than counting copies of $\widehat M$ which may contain many cycles.

 Let $X=(v_1\dots v_\nu)\in\X_3$ and $v_av_{a'}$ be an edge in $M_X$  with $a < a'$.
 For a tepee $v_awv_{a'}$  over $\{v_a,v_a'\}$ in $G$ we will call
 $v_aw$ the {\em left leg} of the tepee and $v_{a'}w$ the {\em right leg} of the tepee.
 Once we can distinguish between the two edges of every
 tepee in $T(X)$, it is easy to define $S(X)$.
 \begin{defn}\label{def1}\rm
For each $X\in\X_3$,  let $S(X)$ be the spanning subgraph  of $\widehat X$   whose edges are the left legs of all
 tepees in $T(X)$. The subgraphs $S(X)$, $X\in\X_3$, will be called {\it  special constellations}.
 \end{defn}
 Since for all $X \in \X_3$
 the subgraph $\widehat{X}$ is isomorphic to $\widehat M$,
 it follows that for all $X \in \X_3$, constellations $S(X)$ have the same
 isomorphism type.
\begin{defn}\label{defS}\rm
Let  $S$ denote the common isomorphism type of all special constellations $S(X)$.
 That is,  $S$ is a subgraph of $\widehat M$  with $V(S) =
V(\widehat{M})$, which is a union of $\nu$ vertex disjoint
 stars, each centered at one of the vertices of~$M$. Their degrees will be denoted
 by $\phi_1 \dots \phi_{\nu}$, where $\phi = \sum \phi_a$, and $\phi_\nu=0$ (because $x_\nu$ is never adjacent to the left
 leg of a tepee containing it).
\end{defn}
It will now be convenient to relabel the vertices of $S$ and denote the
 neighbors of $x_a$ in $S$ by $u_{ab}$, where $b=1\dots\phi_a$. Accordingly,
 we  modify the convention from the previous subsection; namely,
we label the vertices of any special constellation $S(X)$  by $v_a$ and $w_{ab}$
 in a such a way that there is an isomorphism between $S$ and $S(X)$
which maps $x_a$ onto $v_a$ and $u_{ab}$ onto $w_{ab}$ for all
 $a=1\dots\nu$ and $b=1\dots\phi_a$.
Finally, also to conform with this new notation, let us relabel the sets of the partition $\pi_0$:
$$\pi_0 = \{ V_1 \dots V_\nu, W_{11} \dots W_{1\phi_1} \dots W_{\nu-1,1}
 \dots W_{\nu-1,\phi_{\nu-1}} \},$$
so that for every  $X\in\X_3$, the special constellation $S(X)$ satisfies $v_a\in V_a$ and  $w_{ab}\in W_{ab}$ for
all $a=1\dots \nu-1$ and $b=1\dots\phi_a$.

\begin{defn} \rm
We say that a copy $S_0$ of $S$ is consistent with $\pi_0$ if there is an isomorphism $f:S\to S_0$ for which
$f(x_a)\in V_a$ and $f(u_{ab})\in W_{ab}$ for all $a=1\dots\nu$ and $b=1\dots\phi_a$.
\end{defn}
 Clearly, all special
constellations $S(X)$ are consistent with $\pi_0$, but there may be many other copies of $S$ in $G$ which are
consistent with $\pi_0$ too. We will simply call them  {\it constellations}.

\begin{defn}\label{def2}\rm
  A {\em constellation} is any copy $S_0$ of $S$ in $G$ which is
 consistent with the partition $\pi_0$. A {\em special constellation}
  is a copy $S_0$ of $S$ in $G$ for which there exists $X \in \X_3$ such that $S(X)=S_0$.
 The set of all constellations in $G$ will be denoted by $\cee$,
 and the set of all special constellations will be denoted by~$\mathcal{S}$.
 \end{defn}
 Note that the above definition of special constellations is consistent with the previously stated Definition \ref{def1}.

\begin{figure}[hbt]
\begin{center}
    \input{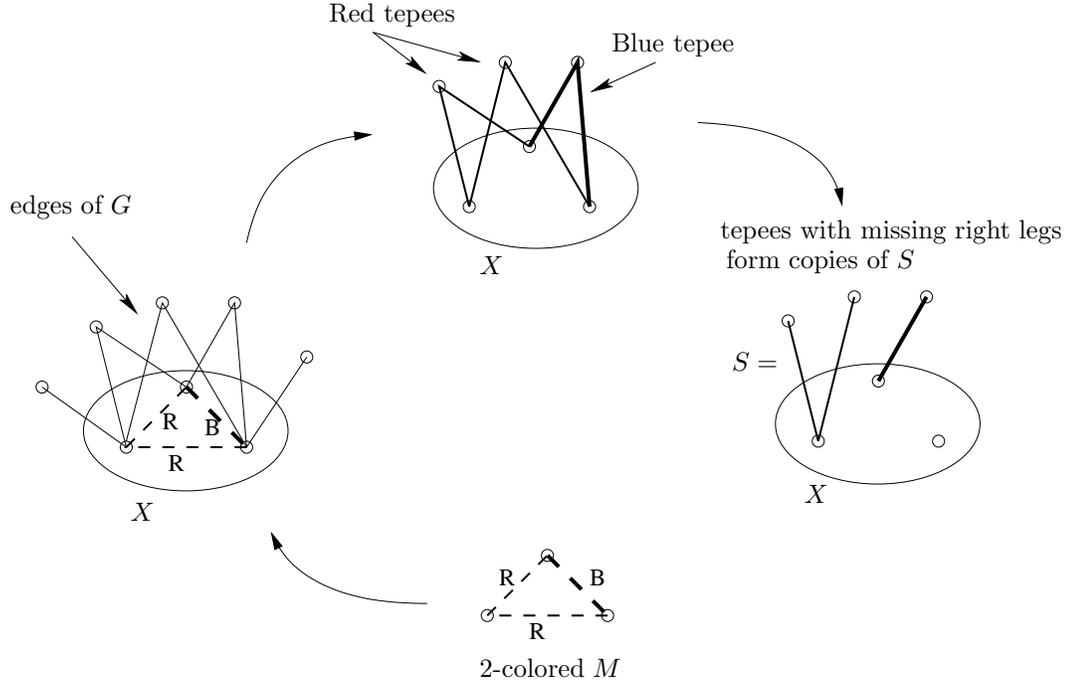}
      \caption{Tepees over $X$, with and without right legs -- creation of a special constellation $S$}
\end{center}
\end{figure}

 Although we claim that removing the right leg of each tepee
 will make our lives easier, we will need to recall them several times in the paper.
 This will be done via the following ``Missing Leg Property".
 To state it formally, we introduce  another piece of  handy  notation related to the indexing
 of the vertices of every $S(X)$, or equivalently, of $S$.
 Define $a'(ab)$ to be the unique
 index $a'$ such that $x_au_{ab}x_{a'}$ is a tepee in $\widehat M$.
 This is well defined because the tepees composing $\widehat M$ are pairwise edge disjoint.

 \begin{observation}[``The Missing Leg Property"]\label{miss}
 Let $X \in \X_3$.
 \begin{enumerate}
 \item[(a)] For every edge $v_aw_{ab}$ of $S(X)$
 the pair $w_{ab} v_{a'(ab)}$ is an edge of $G$ (``the missing leg").
 \item[(b)] $v_a$ and $v_{a'(ab)}$ are the only vertices in $X$
 which are neighbors of $w_{ab}$ in $G$ (only one missing leg per edge in $S(X)$).
 \item[(c)] The edges of $S(X)$ together with the missing legs form all the tepees
 in $T(X)$. In other words there are no other tepees over the edges of $M_X$.
 \end{enumerate}
 \end{observation}

%
%

So, the structural property that distinguishes special constellations  from
 the ordinary ones is the ``Missing Leg Property" that they inherit
 as offsprings of the $\widehat X$'s. The most important inheritance, however, is the color-related property to be formulated in Lemma \ref{TransInEC} below.

But first let us define a partial coloring $\sigma$ of $E(G)$
 which will play a crucial role in our proof. It will assign colors
only to the edges of $\bigcup_{X \in \X_3} S(X)$.
 To begin with, recall that $\sigma'$ is a fixed, triangle-free coloring of $E(M)$, and define an auxiliary coloring of the pairs $(a,b)$, $a=1\dots\nu$, $b=1\dots \phi_a$, by
 $ \sigma(a,b) = \sigma' (x_ax_{a'(ab)})$.

Now, for every edge $vw \in \bigcup_{X \in \X_3} S(X)$
 with $v \in V_a$ and $w \in W_{ab}$, let $\gs(vw) = \sigma(a,b).$
 Note that for any $X \in \X_3$ and any edge $e \in S(X)$,
 $\sigma(e)$ is red or blue according to whether $e$ belongs
 to a tepee in $RT(X)$ or $BT(X)$. Note also that the color $\sigma(e)$ is
 determined by the partition classes
 of $\pi_0$ to which the endpoints of $e$ belong.

 The partial coloring $\sigma$ is a lens through which
 we can study  complete colorings of $G$.
  For a coloring $\chi$ of the edges of $G$, let
 \beq{agree}
{\rm Agree}(\chi,\gs) = \left\{ e \in \bigcup_{X \in \X_3}S(X) : \chi(e) = \sigma(e)
\right\}.\enq
 Now comes a crucial observation that justifies the whole construction
 erected in this section.
 \begin{lemma}
 \label{TransInEC}
 Let $\chi$ be a triangle-free coloring of the edges of $G$.
 Then ${\rm Agree}(\chi,\gs)$ is a hitting set of $\mathcal{S}$, the family
 of special constellations.
 \end{lemma}
 \nin{\bf Proof:} We must show that for
 every $X \in \X_3$ there exists an edge $e \in S(X)$
 such that $\sigma(e) = \chi(e).$
 Lemma \ref{harm} states that
 for every $X \in \X_3$ and for every triangle-free coloring $\chi$, there
 exists either a tepee
 in $RT(X)$ whose edges are colored red by $\chi$, or a tepee in $BT(X)$
 whose edges are colored blue. The left leg of such a tepee is
 an edge of $S(X)$ for which $\sigma$ and $\chi$ agree, by the definition
 of $\sigma$. \qed
We may summarize this section as follows: we have shown that the triangle-free colorings of $G$ are somewhat structured, in the sense that any such coloring must ``contain" a hitting set of the family of special constallations (see Lemma \ref{TransInEC}). As it turns out, this structure imposed on triangle-free colorings is quite restrictive; in Section \ref{cores}, we will capture this restriction in a convenient way. (We are not able to say at this point what we mean by this precisely; let us simply appeal to an analogy and refer the reader to the ``cores" from Section \ref{ill}, which capture a strong structural property of cover sets of $\eps$-regular bipartite graphs.)
In what follows, it will be convenient to consider the hypergraph whose vertices are the edges of $G$ and the hyperedges are the edge sets of the special constallations. Lemma \ref{TransInEC} tells us that triangle-free colorings of $G$ give rise to hitting sets of this hypergraph (or cover sets in the terminology of Section \ref{ill}). In Section \ref{regular}, we will discuss how to ``regularize" this hypergraph; from this regularization we will be able to deduce that the structural property imposed on triangle-free colorings of $G$ is indeed quite strong (again, recall the discussion in Section \ref{ill}).

 \section{Regularity}
 \label{regular} In this section we state and prove a regularity
 theorem that is the central technical tool in this paper. It is a
 generalization of the celebrated Szemer\'edi Regularity Lemma to a
  setting in which the density is measured not merely by the number of edges but with respect to the
  number of copies of certain subgraphs.
 \subsection{Classical Regularity}
 In this subsection we introduce all the necessary  results
 concerning the background on standard graph regularity. See Section \ref{Regular}
 for more information on regularity.

 \subsubsection{Dense Regular Graphs}\label{dense}
In the classical regularity lemma of Szemer\'edi \cite{Szem}  a crucial parameter is the density of a bipartite
graph defined as a ratio of the number of edges of the graph to the number of all potential edges.

\noindent {\it In the following, $B$ stands for a bipartite graph with bipartition $(U,V)$.}

\begin{defn}\label{densit}\rm Let $U' \subseteq U$
 and  $V'\subseteq V$. {\it The density} of the pair $(U',V')$ is defined as
$$d(U',V')=\frac{e_B(U',V')}{|U'||V'|},$$
where $e_B(U',V')$ is the number of edges of $B$ with one endpoint in $U'$ and the other in $V'$ (such edges are
said to {\it belong to the pair $(U',V')$}). We will sometimes write $d(B)$ for $d(U,V)$.
\end{defn}

The $\eps$-regularity of $B$ reflects low discrepancy between the densities of large, induced subgraphs of $B$.

 \begin{defn}\label{usual_regularity}\rm
 Let $\eps >0$. A bipartite graph $B$ is called {\it $\eps$-regular}, if for each choice of $U' \subseteq U$
 and  $V'\subseteq V$, with $|U'| \ge \eps |U|$ and $|V'| \ge \eps |V|$, we have
 \[ |d(U',V') - d(U,V)| < \eps.\]
Sometimes the pair $(U,V)$ itself is called $\eps$-regular.
 \end{defn}

 The following result from \cite{DLR} and \cite{adlry}  gives a sufficient criterion
 for regularity in terms of  vertex degrees and co-degrees only.
 \begin{lemma}\label{DLR}
 Let $\eps>0$ and $d(U,V) =d>2\eps$. Assume further that
 $|U| \ge 2/\eps$ and denote by $W$ the set
 of all pairs $\{v,v'\}$ of vertices of $U$ for which

 \begin{list}{}{}
 \item[(i)]
 ${\rm deg}(v), {\rm deg}(v') > (d - \eps)|V|$,  and
 \item[(ii)]
 ${\rm codeg}(v,v') < (d + \eps)^2|V|. $
 \end{list}
 If $|W| > \frac{1}{2} (1-5\eps)|U|^2$, then $B$ is $\eps'$-regular, where $\eps' = (16 \eps)^{1/5}$.
 \end{lemma}

The next, simple fact is a direct consequence of Definition \ref{usual_regularity}.

 \begin{observation}\label{subreg}
  Let $\eps'>\eps>0$. If $B$ is
 $\eps$-regular, and  $V' \subset V$, with $|V'| = \eps'|V|$, then the pair $(U,V')$ is $\eps''$-regular, where $\eps''=\eps/\eps'$.
 \end{observation}

One of the basic applications of $\eps$-regularity is for counting  copies of a given graph in a structure
composed of several $\eps$-regular and dense pairs of vertex sets. We do not give references, as the result seems
to belong to the ``local folklore".
\begin{prop}\label{graphcount}
For every $\eps>0$, $0<d<1$ and an integer $k$ such that $\eps<d^k$ the following is
true. Let $V_i$, $i=1\dots k$, be disjoint vertex sets with $\min|V_i|>n_0=(2/d)^k$ and let $B_{ij}$, $1\le i<j\le k$, be
$\eps$-regular graphs with bipartitions $(V_i,V_j)$ and densities $d(B_{ij})=d_{ij}$, where $\min d_{ij}>d$. Then
the number $\#(K_k,\bigcup_{i,j}B_{ij})$ of complete graphs $K_k$ contained in the union $\bigcup_{i,j}B_{ij}$
satisfies, with $\eps'= \eps d^{-k}$,
$$\#(K_k,\bigcup_{i,j}B_{ij})\stackrel{\eps'}{\sim}\prod_{i=1}^k|V_i|\prod_{i,j}d_{ij},$$
\end{prop}

If a less precise estimate of the copies of $K_k$ is needed, then the $\epsilon$-regularity can be replaced with a
less restrictive notion. We say that a graph $G$ is $(\varrho,d)$-dense if for all $U\subseteq V(G)$ with $|U|\ge
\varrho|V(G)|$, we have $e_G(U)\ge d{|U|\choose2}$. When verifying this property it suffices to check only all
subsets $U$ of size $|U|=\lfloor\varrho|V(G)|\rfloor$.  The next result comes from \cite{RR3}.

\begin{prop}\label{L2.3}
For all $d>0$ and $k$ there exist $\varrho,n_0,c_0>0$ such that for every $(\varrho,d)$-dense graph $G$ with
$n=|V(G)|\ge n_0$ we have $\#(K_k,G)\ge c_0n^k$.
\end{prop}

Although we do not need it in our proofs, we now state the celebrated Szemer\'edi Regularity Lemma. This in not
exactly the original version, but one which can be compared with the sparse version from the next subsection. Recall the definition of equipartition given at the end of Section \ref{intro}.

 \begin {thm}[Szemer\'edi Regularity Lemma]
  For all $\eps>0$ and integers $t_0$ and $r$, there  exist $T_0$ and $n_0$
 such that if  $G_1 \dots G_r$ are graphs with a common vertex set $V$, \
$|V|=n > n_0$, then every  equipartition of $V$ into
$t_0$ parts has a refined equipartition  $V = W_1 \cup \ldots \cup
W_{t}$, where  $t_0\le t \le T_0$, and
 all but at most $\eps n^2$
 edges of $G_1\cup\ldots\cup G_r$ belong to the pairs $\{W_i,W_j\}$ that are
 $\eps$-regular with respect to all  graphs $G_i$,
 $i= 1 \dots r$.
 \end{thm}

 \noindent
 \subsubsection{Sparse Regular  Graphs}\label{srg}
 The graph $G$, one of the main constituents of the proof of Theorem \ref{allweneed},
  is relatively sparse, with density of order $1/\sqrt n$.
 Below we present a version of the regularity lemma designed for sequences of graphs with $n$ vertices and
  density  $p=o(1)$.

As before, let $B$ be a bipartite graph with bipartition $(U,V)$, and let $0<p\le1$.
\begin{defn}\label{sparse_dens}\rm Let $U' \subseteq U$
 and  $V'\subseteq V$. The  {\it $p$-density} of the pair $(U',V')$ is defined as
$$d_p(U',V')=\frac{e_B(U',V')}{p|U'||V'|}.$$
We will sometimes write $d_p(B)$ for $d_p(U,V)$.
\end{defn}
As in the dense case, the $(\eps,p)$-regularity of $B$ (as defined below) reflects low discrepancy between
the $p$-densities of large, induced subgraphs of $B$.

\begin{defn}\label{sparse_reg}\rm
  Let $\eps>0$ and $0<p\le1$. A bipartite graph $B$ is
 {\it  $(\eps,p)$-regular} if for each choice of $U' \subseteq U$
 and  $V'\subseteq V$, with $|U'| \ge \eps |U|$ and $|V'| \ge \eps |V|$, we have
 \[ |d_p(U',V') - d_p(U,V)| < \eps.\]
Sometimes the pair $(U,V)$ itself is called  $(\eps,p)$-regular.

 \end{defn}
There is an analog of Observation \ref{subreg} in the sparse case.
\begin{observation}\label{p-subreg}
  Let $\eps'>\eps>0$ and $0<p\le1$. If $B$  is
 $(\eps,p)$-regular, and $V' \subset V$, with $|V'| = \eps'|V|$, then the pair $(U,V')$ is $(\eps'',p)$-regular, where $\eps''=\eps/\eps'$.
 \end{observation}

A basic feature of $\epsilon$-regular graphs is that most of the vertex degrees are
under control. Here is a useful fact exploited in Section \ref{randomgraphs}.

\begin{prop}\label{degrees}
If $(U,V)$ is $(\epsilon,p)$-regular and $W\subseteq V$ satisfies $|W|\ge\epsilon|V|$, then at least
$(1-\epsilon)|U|$ vertices of $U$ have each at least $(d_p(U,V)-\epsilon)p|W|$ neighbors in $W$.
\end{prop}

The above notions of $p$-density and  $(\eps,p)$-regularity extend naturally to pairs of disjoint subsets of
vertices of non-bipartite graphs, for which we now introduce one more notion, specific just to the sparse case.

\begin{defn}\label{nodenspa}\rm
  For $D>0$ and  $0<p\le1$ we say
 that a graph $F$ on $n$ vertices has no {\em $(D,p)$-dense patches}, if
 for every two disjoint sets of vertices $W_1,W_2\subset V(F)$ with $|W_1|, |W_2| >
 n/\log n$, we have
 $$ d_p(W_1, W_2) < D.$$
\end{defn}
 This property,  called $(p;D,1/\log n)$-boundedness in \cite{JLR}, page 215, guarantees that the so called partition index is bounded from above, and thus allows to repeat the proof of the Szemer\'edi Regularity Lemma mutatis mutandis.

 But there is another profit from this assumption. In general, there is no
simple analog of Proposition \ref{graphcount} in the sparse case.
 However, there are several places in this paper where we need estimates of
 the number of stars consistent with a certain family of disjoint subsets of
 vertices. Fortunately, in this simple case, a similar counting result is
 true, provided some kind of boundedness is assumed. We prove it here in a
 specific form suitable for our purposes.

\begin{prop}\label{stars}
For every $\eps>0$, $D>0$, and integer $k\ge2$ such that $\eps\le\min\{k^{-4},2^{-6}D^{-6k}\}$ there exists
$n_0$ for which the following is true. Let $V_0$, $W_j$, $j=1\dots k$, be disjoint vertex sets of sizes $n_0\le
|W_j|/2\le|V_0|\le |W_j|$, $n=|V_0|+|W_1|+\cdots+|W_k|$ and $0<p=p(n)\le1$. Let $B_{j}$, $1\le j\le k$, be
$(\eps,p)$-regular graphs with bipartitions $(V_0,W_j)$ and $p$-densities $d_p(B_{j})=d_{j}$, where $\min
d_j\ge\eps^{1/2k}$.  Assume further that the graphs $B_j$ have no $(D,p)$-dense patches and that the maximum
degree in the graph $\bigcup_{j}B_{j}$ is at most $Dnp$.  Let $\#(S_k,\bigcup_{j}B_{j})$  be the number of copies
of the star $S_k$ contained in $\bigcup_{j}B_{j}$ and such that the vertex of degree $k$ is in $V_0$ and the
vertices of
 degree one each belong to a distinct set $W_j$.  Then we have,
 $$\#(S_k,\bigcup_{j}B_{j})\stackrel{\eps^{1/3}}{\sim} |V_0|p^k \prod_{j=1}^k |W_j|d_j .$$
\end{prop}

 \nin{\bf Proof:}  For any subset $V'$ of $V_0$, set
$S(V')=\#(S_k,\bigcup_{j}B'_{j})$, where $B'_j$ is the subgraph of  $B_j$ induced by the vertex set $V'\cup W_j$.
Let ${\rm deg}_j(v)={\rm deg}_{B_j}(v)$ be the number of neighbors of vertex $v$ in graph $B_j$. Clearly,
$$S(V')=\sum_{v \in V'} \prod_{j=1}^k {\rm deg}_j(v),$$
so we would be done if we knew that all, not just {\it almost} all, vertices of $V_0$ satisfy  ${\rm
deg}_j(v)\stackrel{\eps}{\approx}d_jp|W_j|$ for all $j$. Unfortunately, this is not so, and we have to somehow handle
the exceptional degrees. Of course, for a lower bound on $S(V_0)$ we can just leave them out. Let
 $$V_{small}^j = \{ v \in V_0  :{\rm deg}_j(v) < (d_j- \eps )p|W_j|\}.$$
 By the $(\eps,p)$-regularity of $B_j$ it follows that $|V_{small}^j| <
 \eps |V_0|.$
 Let
$$ V_{small}= \bigcup_{j=1}^k V_{small}^j.$$
 Then
 \begin{eqnarray*}
  S(V_0) &\geq& S(V_0 \setminus V_{small}) = \sum_{v \in V_0 \setminus V_{small}} \prod_{j=1}^k {\rm deg}_j(v) \geq
 \sum_{v \in V_0 \setminus V_{small}} \prod_{j=1}^k (d_j- \eps) p|W_j| \\
 &\stackrel{*}{\geq}&  |V_0|(1-k\eps)p^k \prod_{j=1}^k d_j(1-\eps^{3/4})|W_j| \geq
 |V_0|(1-k\eps ) (1-\eps^{3/4})^k \prod_{j=1}^k d_j|W_j|\\
 &\ge & |V_0|(1- (k+1)\eps^{3/4}) \prod_{j=1}^k d_j|W_j| \stackrel
 {**}{\ge} (1-\eps^{1/3})|V_0| \prod_{j=1}^k d_j|W_j|.
 \end{eqnarray*}
 Note that for inequality (*) above we use $d_j  \ge \eps^{1/4} $, which follows from
 $d_j \ge \eps^{1/2k}$ and $k \ge 2$. The assumption $1/k \ge \eps^{1/4}$ validates inequality (**).

 We now proceed to prove a corresponding upper bound
 on $S(V_0)$. We will partition the vertices of $V_0$ into three
 classes $V_{med} , V_{big} $ and $V_{huge}$, according to their
 degrees, and will see that, as expected, the majority of the contribution to $S(V_0)$
 comes from $V_{med}$, whereas the
 contribution from the other two classes is negligible. Let
  $$V_{med} = \{ v \in V_0 : \forall j \  {\rm deg}_j(v) \leq  (d_j+ \eps) p|W_j| \},$$
  $$V_{huge}=\{ v \in V_0 : \exists j \  {\rm deg}_j(v) > Dp|W_j| \}$$
  and
  $$ V_{big} = V \setminus (V_{med}  \cup V_{huge}) .$$
  Note that $S(V) = S(V_{med}) + S(V_{big}) + S(V_{huge}).$
  As in the case of the lower bound, it is immediate that
  \beq {<}
  S(V_{med}) \le |V_0|\prod_j (d_j +\eps) p|W_j|
 \stackrel{***}{\le}  |V_0|p^k (1+\eps^{1/3}/2) \prod_j d_j|W_j|,\enq
 where, as above, (***) follows from the assumptions that
 $d_j \ge \eps^{1/4}$ and $k\le \eps^{-1/4}$.

  Next, observe that by the $(\eps,p)$-regularity of $B_j$'s we have $|V_{big}| < \eps k |V_0|$.
  Hence
  \begin{eqnarray}
 \label{big}
  S(V_{big}) &<& |V_{big}|(Dp)^k\prod_j|W_j| < \eps k |V_0| (Dp)^k \prod_j|W_j| \nonumber\\
&\leq & \sqrt{\eps}
  |V_0|(\eps^{\frac{1}{2k}}Dp)^k\prod_j|W_j| \leq \sqrt{\eps}|V_0|(Dp)^k\prod_j d_j|W_j| .\end{eqnarray}

  Finally, by the assumptions on the maximum degree and the lack of $(D,p)$-dense patches, for all
   $v \in V_0$ and for all $j=1\dots k$ we have  $d_j(v)< Dnp$ and
    $|V_{huge}| < k  {n}/{\log n}$. Hence
   \beq {huge}
   S(V_{huge}) <k \frac {n}{\log n} (Dnp)^k = o(|V_0|p^k\prod_j d_j|W_j|).
   \enq
 Putting (\ref{<}), (\ref{big}) and (\ref{huge}) together gives
 $$ S(V) \leq {(1+ \eps^{1/3})} |V_0|p^k \prod_j d_j|W_j|.$$
 \hfill $\Box$
\bigskip

Of course, the above result can easily be adopted to asymptotically count  copies of star forests as well.
\begin{ex}\label{cee}\rm Consider a graph $G\in\G$, a star forest $S$,  and a partition $\pi_0$ of $V=V(G)$ as in Section \ref{TandC}. For each $a=1\dots\nu$ and $b=1\dots\phi_a$, let $V_a'\subseteq V_a$ and $W_{ab}'\subseteq W_{ab}$ be sets of size at least $n/\log n$.
By Property (P5) of $\G$ there are no $(2,p)$-dense patches in $G$. Moreover, by the same property, the bipartite
subgraphs  $F'_{ab}=G[V'_a,W'_{ab}]$ are all $(2n^{-1/5},p)$-regular. Thus, by $\nu$ simultaneous applications of
Proposition \ref{stars} (one for each star of $S$), we obtain an asymptotic formula for the number $\kappa'$ of the constellations in $G$ contained in the union of $F'_{ab}$'s:
$$\kappa'\sim\prod_{a=1}^{\nu}|V'_a|\prod_{b=1}^{\phi_a}|W'_{ab}|.$$
In particular, the number  $\kappa=|\cee|$ of all
constellations in $G$  consistent with $\pi_0$ satisfies:
$$\kappa=|\cee|\sim m^{\nu+\phi}p^\phi\quad\mbox{ where }\quad m=\frac n{\nu+\phi}.$$
(See, e.g., Figure 4)
\end{ex}


\begin{figure}[hbt]
\begin{center}
    \input{fig6_new.tex}
      \caption{Illustration of the example}
\end{center}
\end{figure}

Finally, we state a sparse version of the Szemer\'edi Regularity Lemma.
Theorem \ref{SSRL} below was independently observed by Y. Kohayakawa and  V. R\"odl
and, in fact, its proof is a simple modification of Szemer\'edi's proof. For applications of this variant of the regularity lemma see \cite{Kohayakawa97Szemeredi},
 \cite{KR} and
\cite{JLR}, Section 8.3.

 \begin {thm}[Sparse Regularity Lemma]\label{SSRL}
  For all $\eps>0$, $D>0$, and integers $t_0$ and $r$, there  exist $T_0$ and
 $n_0$  such that if $G_1 \dots  G_r$ are graphs with a common vertex set
$V$, \  $|V|=n > n_0$, and $p=p(n)$ is such that $G_i$ does not have $(D,p)$-dense patches for all $i=1,2, \ldots
, r$, then every equipartition of $V$ into $t_0$ parts has a refined equipartition $V = W_1 \cup \ldots \cup W_{t}$, where  $t_0\le t \le T_0$, and
all but at most $\eps n^2p$
 edges of $G_1\cup\ldots\cup G_r$ belong to the pairs $\{W_i,W_j\}$ that are
  $(\eps,p)$-regular with respect to all  graphs $G_i$,
 $i= 1 \dots r$.
 \end{thm}
 This theorem is essentially the same as the one presented in
 \cite{KR}.
  The theorem in \cite{KR} guarantees that all but an $\eps$-proportion of the {\em pairs}  $(W_i, W_j)$ in the partition are
 $(\eps,p)$-regular with respect to all $G_i$.
  Given the fact that there is an upper bound on the
 $p$-density of every pair (i.e., no $(D,p)$-dense patches),
 we have that at  most an $\eps D$-proportion of the {\em edges} belong to such
 pairs. Hence the theorems are equivalent (after re-scaling the
 constants).

  \subsection{The Subgraph Regularity Lemma}
 Let us summarize the state-of-art at the end of Section \ref{TandC}.
 Given  a graph $G\in\G\setminus \R$,
as a far reaching consequence of assumption (\ref{assumption}) of Theorem \ref{allweneed},
 we have constructed a star forest $S$ and a partition $\pi_0$
 of the vertices of $G$ into $|V(S)|$ sets of equal size $m=\lfloor n/|V(S)|\rfloor$ and, possibly, one set of less
than $|V(S)|$ vertices. (For clarity, we will be further assuming that $n$ is divisible
by $|V(S)|$.)
Based on this partition, we have focused on  the set $\cee$ of all constellations, that is all copies of $S$ in
$G$ that are consistent with $\pi_0$. Among the constellations we have distinguished
a set $\mathcal {S}$  of  special constellations.
These are the ones that correspond to selected copies $M_X$ of
   a certain  graph $M$, such that $G\cup M_X \in\R$.
    Later in the proof of our main result, we are going to apply a regularity lemma
    to the family  $\mathcal {S}$. The purpose of this section is to prove
    such a  regularity lemma in a more abstract setting, not using the
    concrete definition of the family  $\mathcal {S}$.

Consider a hypergraph whose vertices are the edges of $G$ and whose edges are
 the members of $\mathcal {S}$. By abuse of notation we may also denote the hypergraph
itself by $\mathcal {S}$. As we roughly explained in Section \ref{sneak}, we wish
to apply a regularity lemma to $\mathcal {S}$, and partition its vertices, i.e. the edges of $G$.
Our objective will be to find such a partition so that  various parts of the partition will carry
the edges  of $\mathcal {S}$ (special constellations) in a fairly regular manner with respect to all
constellations, just as regular pairs in a usual Sz\'emeredi partition
of a graph have edge density that is, in a sense, smoothly distributed.

This partition of $\mathcal {S}$ will rely heavily on the underlying graph $G$, and we will start out by
partitioning the vertices of $G$, thus inducing a partition of the edges of $G$
(which are the vertices of $\mathcal {S}$!). Only in the next step will we focus on the sets
thus formed and introduce a further partition of the {\em edges} of $G$ that is {\em not} induced
simply by partitioning the vertices of $G$.
We  will then iterate this procedure in a manner similar to the one introduced by Frankl  and R\"odl in \cite{FR1},
see the next subsection.
We are aware that the two roles of edges of $G$ as vertices of $\mathcal {S}$ may be confusing, and hence,
having pointed out the hypergraph structure, we will now stick to the language of subgraphs,
thus avoiding ambiguity.

 Let $F$ be the set of all edges of $G$ with one endpoint in $V_a$ and the
 other in $W_{ab}$ for all $a=1\dots\nu$ and $b=1\dots\phi_a$.   We will
 identify $F$ with the spanning subgraph of $G$ consisting of all these edges.

The regularity lemma we are aiming at will later be applied only to subgraphs
 $F$ of graphs $G\in\G\setminus \R$, created in the above way by selecting
 first a star forest $S$ and an initial partition $\pi_0$. However, it will be
 convenient to formulate it more generally for all graphs like $F$, that is,
 for graphs structured by $S$ and $\pi_0$ but with no reference to $\G$ and to assumption
 (\ref{assumption}) of  Theorem \ref{allweneed}. Moreover, with no additional effort, our regularity lemma
 can be proved with respect to subgraphs other  than $S$.
Therefore, we will soon set up a general framework in which the regularity lemma will be stated and proved,
bearing in mind that the only application we need is that for the star forest $S$. But first, in order to give
stronger foundations for our generalization, we describe one more instance of regularity.

\subsubsection{The Hypergraph Regularity of Frankl and R\"odl}\label{fra}

Recently, Frankl and R\"odl \cite{FR1} proved a powerful regularity lemma for 3-uniform hypergraphs. A simplified
version of that lemma  is described  briefly here as another piece of motivation for the general framework we are
up to.

Given three disjoint sets $V_1,V_2,V_3$, {\it a triad} is a triple $F = (F^{12},F^{23},F^{13})$ of bipartite
graphs with vertex sets $V_1 \cup V_2, ~ V_2 \cup V_3$ and $V_1 \cup V_3$. We will identify $F$ with the union of
its three component graphs. We refer to a 3-partite 3-uniform hypergraph ${\cal H}$ with a 3-partition $(V_1,V_2,
V_3)$ as a {\it 3-hypergraph}. For a triad $F$ with the same vertex partition as ${\cal H}$, we say that $F$ {\it
underlies} ${\cal H}$ if ${\cal H} \subseteq {\rm Tr}(F)$, where ${\rm Tr}(F)$ is the family of the vertex sets
of the triangles in the graph $F$.

Let ${\cal H}$ be a 3-hypergraph with an underlying triad $F$. The {\it density} of ${\cal H}$ with respect to $F$
is defined by
$$ d_{\cal H} (F) = \frac{|{\cal H} \cap {\rm Tr} (F)|}{|{\rm Tr}(F)|}.$$
In other words, the density measures the proportion of triangles of $F$ which are triples of ${\cal H}$.

Let $\delta > 0$. A 3-hypergraph ${\cal H}$ is said to be {\it $\delta$-regular} with respect to an underlying
triad $F$ if  for every subtriad ${{Q}}  = ({{Q}}^{12} , {{Q}}^{23} , { Q}^{13}) $ of $F$, ${ Q}^{12}  \subset
F^{12}, { Q}^{23} \subset F^{23}, { Q}^{13}  \subset F^{13},$ with $|{\rm Tr}({ Q} )| > \delta |{\rm Tr}(F)|$ we
have
$$|d_{\cal{H}} ({ Q})-d_{\cal{H}} ({ F})|<\delta.$$
Thus, the $\delta$-regularity of $\h$ reflects low discrepancy between the densities of ${\cal H}$ with respect to
those subtriads of $F$
 which are relatively rich in triangles. With a fixed  ${\cal H}$ in mind,
we call the triad itself $\delta$-regular. A triad which is {\em not} $\delta$-regular is called
{\it $\delta$-irregular}.

 Let $V$ be a set, and let $[V]^2$ denote the set of all 2-element subsets of $V$.  An {\it $(l,t,\epsilon)$-partition}
  $\Pi$ of $[V]^2$ consists of an equipartition $V=
  V_1\cup\cdots\cup V_t$,  together with a
  system of edge-disjoint bipartite graphs $F_k^{ij}$ with bipartitions $(V_i,V_j)$, $1\leq i<j\leq t$,
  $1\leq k\leq l_{ij}\leq l$, such that
\begin{itemize}
\item
  $\left|\bigcup_{k=1}^{l_{ij}}F_k^{ij}\right|=|V_i||V_j|$ for all~$i$, $j$,
  $1\leq i<j\leq t$, and
\item all but at most $\epsilon{{t}\choose{2}}m^2$ pairs
  $\{v_i,v_j\}$, $v_i\in V_i$, $v_j\in V_j$, $1\leq i<j\leq t$, belong to graphs $F_k^{ij}$
  which are $\epsilon$-regular.
\end{itemize}
Note that there are at most ${{t}\choose{3}}l^3$ triads $F$ made up of the graphs of $\Pi$. Let $\calp^{\rm
irr}_\Pi$ be the set of all $\delta$-irregular triads among them.

For a 3-uniform hypergraph $\h=(V,E)$, with $|V|=n$, we say that an $(l,t,\epsilon)$-partition $\Pi$ of $[V]^2$
  is {\it $\delta$-regular} if
$$\sum_{F\in\calp_\Pi^{\rm irr}} |{\rm Tr}(F)|
    <\delta{n}^3.$$

\begin{thm}\label{fr}
  For every $\delta>0$,
  integers~$t_0$ and~$l_0$,
  and all decreasing functions $\epsilon(\ell)$, there exist $T_0$, $L_0$, and $N_0$ such that any
  $3$-uniform hypergraph $\h=(V,E)$ with $|V|>N_0$ admits a
  $\delta$-regular,
  $(l,t,\epsilon(l))$-partition for some~$t$ and~$l$ satisfying
  $t_0\leq t\le T_0$ and $l_0\leq l\le L_0$.
\end{thm}

\begin{remark}\rm
In \cite{FR1} a more general and powerful hypergraph lemma guarantees, for every integer function $r=r(\ell)$, a stronger, $(\delta,r)$-regular partition. Above we discussed the case $r\equiv 1$ only.
\end{remark}

\subsubsection{General Subgraph Framework}

In what follows $H$ always stands for a fixed graph with respect to copies of which other graphs are regularized.
The reader can bear in mind that the cases of $H=K_2$ and $H=K_3$ correspond, respectively, to the Szemer\'edi and
Frankl-R\"odl regularity schemes, and that $H=S$, the star forest from Definition \ref{defS}, will be our principal (and only) application.

Given a graph $H$ with $h$ vertices and an integer $m$, let $H^m$ be the order $m$ blow-up of $H$, i.e.,
the $h$-partite graph  obtained by replacing
each vertex $x$ of $H$ by a set $V_x$ of $m$ vertices, and by replacing every edge $xy$ of $H$ by the complete
bipartite graph $K_{m,m}$ spanned between $V_x$ and $V_y$. Let $\F(H,m)$ denote the family of all subgraphs $F$ of  $H^m$.
The graphs $F\in\F(H,m)$ will be called {\it $H$-graphs} and a copy $H'$ of
 $H$ in $F$ will be declared {\it consistent} if for each $x\in V(H)$, its
 image $x'$ belongs to $V_x$ (see Figure 5 ). By a copy of $H$ in an $H$-graph we will always
 mean one which is consistent. Let $\cee(F)=\cee_{H}(F)$
 be the set of all copies of $H$ in $F$.

Besides graph $F$, the other important input to our regularity lemma is  an arbitrarily specified subfamily of
$\cee(F)$, the elements of which will be called {\it special copies}. The subfamily itself will be denoted by
${\mathcal {S}}$. By choosing this notation  we put emphasis on our primary application with $H=S$ and $\cal S$ --
the family of special constellations. To obtain the standard Szemer\'edi setting take $H=K_2$, $F=K_{m,m}$ and
$\cal S$ -- the set of edges of a bipartite graph. For the hypergraph setting of Frankl and R\"odl take $H=K_3$,
$F=K_{m,m,m}$ and $\cal S$ -- the set of triples of a 3-hypergraph.


\begin{figure}[htp]\label{blowup}
\begin{center}
    \input{fig7b.pstex_t}
      \caption{consistent copy of $H=C_4$}
\end{center}
\end{figure}

 \subsubsection{Measuring Density}
 In the classical case of graphs, regularity is related
 to the density measured by the ratio of the number of edges to all potential
edges. For 3-uniform hypergraphs, regularity will be with respect to the ratio of triples to the triangles in the underlying triad. Similarly,
we will measure the density of various (spanning) subgraphs $R$ of $F$ by the ratio of the number
of special copies of $H$ to all copies of $H$ contained in $R$. Since
our primary application is in the sparse case, the density will have to be
normalized by a scaling factor $p^*$, which will typically depend both on $H$ and on the average density
of the sparse graph in question.

Having fixed $H$, let us also fix  $F\in\F(H,m)$ with  $V(F)=\bigcup_{x\in V(H)} V_x$, ${\mathcal {S}}\subseteq \cee_H(F)$ and $0<p^*\le1$. Any
subgraph  $R$ of $F$ is itself a member of $\F(H,m)$ and thus $\cee(R)=\cee_H(R)$ is already defined,
meaning the family of copies of $H$ in $R$. As for the special copies,
though, they can only be inherited from $F$. We set
$$c_R = |\cee(R)|\mbox{ and } \quad s_R = |{\mathcal {S}}\cap\cee (R)|.$$
For all subgraphs $R\subseteq F$, with $c_R>0$, we define their {\it $({\mathcal {S}},p^*)$-density}  by
 \begin{equation}\label{density}d_R=
 \frac{s_R}{p^*c_R}.
\end{equation}
 The normalizing factor $p^*$  helps in bounding the
 density  from below. It should be related closely to the ratio $s_F/c_F$ between the numbers of special copies and all copies of $H$ in the entire graph $F$.

The price we pay for using a normalized density is the potential presence of some kind of dense patches. It is
unavoidable to have small subgraphs $R$ with unbounded $d_R$ (e.g., if $c_R=s_R=1$ and $p^*=o(1)$). Thus, all we
can do is to control the $({\mathcal{S}}, p^*)$-density of large subgraphs of $R$, where ``large" means ``rich in
copies of $H$".  Set
$$\kappa=|\cee(F)|$$
for the total number of copies of $H$ in $F$ (``total volume'' of $F$).
\begin{defn}\label{nodensePa}\rm
For fixed $H$ and $D^*>0$, an $F \in {\cal F}(H,m)$ with $n=hm$ vertices, together with a nonempty  family ${\cal
S}\subseteq{\cee}(F)$, is said to {\it have no $(D^*,p^*)$-dense $H$-patches},
 if for every subgraph $R\subseteq F$ with $c_R\ge \kappa/\log^2n$
we  have $d_R\le D^*$.
\end{defn}
Note that the assumption of no dense patches applied to $R=F$ and the definition of density in (\ref{density}) imply together that
\begin{equation}\label{p^*-lower}p^* \ge \frac{s_F}{D^*c_F} \ge \frac{h^h}{D^*}n^{-h},\end{equation}
unless $s_F=0$.

 The property of having no dense patches, together with an upper bound on the number of copies
 that one edge may belong to,  implies that subgraphs of $F$ with few copies of $H$ must also have few special ones.

\begin{defn}\label{H-uni}\rm
An $H$-graph $F$ with $n$ vertices is called {\it $H$-uniform} if no edge of $F$ belongs to more than
$\kappa/\log^3n$ copies of $H$.
\end{defn}

 \begin{prop}\label{(b)} Assume that an $H$-graph $F$ (together with  $\cal S$) has no $(D^*,p^*)$-dense $H$-patches and is $H$-uniform. Let $R\subseteq F$ with $c_R\le \kappa/\log^2n$. Then $s_R\le (1+o(1))D^* p^*\kappa/\log^2n$.
 \end{prop}

\nin{\bf Proof:} Given $R$ with $c_R \leq \kappa/\log^2n$,
 we add to it extra edges of $F$, enlarging $R$ to obtain a supergraph $R'$
 with $\kappa/\log^2n\le c_{R'} = (1+o(1))\kappa/\log^2n$. This can
 easily  be accomplished  by adding one edge at a time, and noting that
 each edge can increase the number of copies of $H$ by at most $\kappa/\log^3n$.
 Thus, because there are no $(D^*,p^*)$-dense patches in $F$, it follows that
 $d_{R'} \le D^*$, and hence $s_R \leq s_{R'} \leq (1+o(1))D^* p^* c_{R'}.$\hfill $\Box$

\subsubsection{Partitions and Polyads}
The subgraph of an $H$-graph $F$, induced by a pair $(V_x,V_y)$ corresponding to the edge $xy\in H$, will be
denoted by $F_{xy}$. The sets $V_x$ for all $x\in V(H)$ and the graphs $F_{xy}$ for all $xy\in E(H)$ form {\it the
initial partition} $\Pi_0$ of a fixed $H$-graph $F$.
 The regularity lemma we prove in this section involves
 refinements of the sets and graphs of ${\Pi}_0$ and towards this we
 define a $(t,l)$-{\it partition}, by generalizing the definition
 of ${\Pi}_0$.

 Informally speaking, a $(t,l)$-partition is an equipartition obtained by splitting
 each set $V_x$ into $t$ disjoint subsets $V_x^1,\ldots,V_x^t$. This partitions every
 subgraph $F_{xy}$ into $t^2$ {\it tubes} $$F_{xy}^{i,j} = F_{xy}[V_{x}^i, V_{y}^j], \ \ 1 \le i,j \le t,$$
which are induced bipartite
 subgraphs spanned between the partition sets. We then partition every
 tube into $l(xy,i,j) \leq l$ edge-disjoint subgraphs.

 \begin{defn}\label{tl}\rm A $(t,l)$-partition ${\Pi}$ consists of the following ingredients:
 \begin{equation}\label{vertex}
 \hbox{ equipartitions }\quad V_x=V_{x}^1\cup\cdots\cup V_{x}^t,\qquad\qquad x\in V(H)
\end{equation}
 and
 \begin{equation}\label{edge}
 \hbox{ partitions }\quad F_{xy}^{i,j} = \bigcup_{k=1}^{l(xy,i,j)} F_{xy}^{i,j,k},\quad\qquad xy\in E(H),\quad i,j=1\dots t.
\end{equation}
{\it A refinement} of $\Pi$ is any $(t',l')$-partition, with $t'\ge t$, $l'\ge l$, which consists of refinements
of (\ref{vertex}) and (\ref{edge}).
 \end{defn}

\nin In particular, any $(t,l)$-partition is a refinement of the initial  $(1,1)$-partition $\Pi_0$, and a
$(t,l)$-partition divides the vertex set into $th$ classes and the edge set into $\sum_{xy\in E(H)}l(xy,i,j)\le|E(H)|t^2l$
classes.

We now define a basic notion whose role is similar to that of bipartite subgraphs in the classical Szemer\'edi
Lemma, and to that of triads in the hypergraph regularity lemma of Frankl and R\"odl.
  \begin{defn}\label{polyad}\rm A {\em polyad} $P$ consistent with a $(t,l)$-partition $\Pi$ is a subgraph of $F$ obtained by selecting one subclass $V_x^{i(x)}$ for each $x\in V(H)$, and then one subgraph $F_{xy}^{i(x),j(y),k(xy)}$ for each $xy\in E(H)$.
More formally, given $i(x)$ for $x\in V(H)$ and  $k(xy)$ for $xy\in E(H)$,
 $$P=\bigcup_{xy} F_{xy}^{i(x),j(y),k(xy)}$$
is the corresponding polyad.
\end{defn}

 Observe that there is only one polyad consistent
 with the original partition ${\Pi}_0$, which is $F$ itself, and, more generally, there are at most
 $s=t^hl^{|E(H)|}$ polyads consistent with a $(t,l)$-partition (with equality if every tube is partitioned into  exactly
 $l$ parts). Let $\calp=\calp_\Pi$ be the set of all
 polyads consistent with a $(t,l)$-partition $\Pi$.
  Note also that every copy of $H$  belongs to exactly one
 polyad. Consequently,
 recalling that $\kappa= | \cee(F) |$ is the total number of copies of $H$  in $F$, we have
  \begin{equation}\label{kappa}\kappa=\sum_{P \in \calp}c_P.\end{equation}
Finally, notice that if $\Psi$ is a refinement of $\Pi$ then for every polyad $R\in\calp_\Psi$ there exists a
unique polyad $P\in\calp_\Pi$ such that $R\subseteq P$.

 \subsubsection{Regularity}

 Polyads are the substructures of this construction which
 play the same role as $\epsilon$-regular pairs
in a Szemer\'edi  partition of a graph.

\begin{defn}\label{dfn4.3}
 \rm  Let $\delta > 0$. Given an $H$-graph $F$ and a family $\cal S$ of special copies of $H$ in $F$, a polyad $P$ is said to be
 {\it $\delta$-regular}, if for every subgraph $Q \subseteq P$ with $c_Q
 \geq \gd c_P$, we have
 $$|d_P - d_Q| \le \delta.$$
 \end{defn}
In other words, $\delta$-regularity of $P$ reflects low discrepancy between the $({\mathcal
{S}},p^*)$-densities of those subpolyads of $F$ which are relatively rich in ordinary, consistent copies of $H$. A polyad
which is not $\delta$-regular will be called {\it $\delta$-irregular}, or just {\it irregular}.
 The set of all $\gd$-irregular polyads of $\calp_{\Pi}$ will be denoted
 by $\calp^{\rm irr}_{\Pi}$ or just  $\calp^{\rm irr}.$

We are now two notions away from our regularity lemma. The first of them ensures $(\eps,p)$-regularity of most
subgraphs constituting a $(t,l)$-partition.

 \begin{defn}\label{dfn4.1}
 \rm Let $\eps >0$ and $0<p\le1$. A $(t,l)$-partition ${\Pi}$
 is called {\it $(\eps,p)$-uniform} if all
   but $\eps|F|$  edges of $F$ are in
  $(\eps,p)$-regular subgraphs $F_{xy}^{i,j,k}$.
 \end{defn}

 The next definition specifies what kind of regularity we impose on partitions in terms of $\delta$-regular polyads.
 Roughly, the irregular polyads together cannot contain too many copies of $H$. Recall that $\kappa=|\cee(F)|$.

 \begin{defn}\label{dfn4.4}
  \rm Given an $H$-graph $F$ and a family $\cal S$ of special copies of $H$ in $F$, a $(t,l)$-partition $\Pi$ is $\delta$-regular if
 $$\sum_{P \in \calp^{\rm irr}_{\Pi}} c_P \leq \delta\kappa.$$
 \end{defn}

 We now come to the central lemma of this section, which, when applied with $H=S$, the star forest from Section \ref{TandC}, puts into action
 the mechanism that makes the whole proof tick. It describes, in
 the spirit of the Szemer\'edi Regularity Lemma, a decomposition of an $H$-graph $F$
  in which the special copies of $H$ are evenly distributed, i.e. it guarantees the
 existence of a $\gd$-regular partition. In addition, at the same time the partition that is obtained is $(\eps,p)$-uniform. This extra feature is very useful in our application.

 \begin{lemma}[Subgraph Regularity Lemma]\label{reglem}
 For all graphs $H$, constants $\delta>0$, $D>0$, $D^*\ge1$, and
for all functions $\eps(\ell)>0$,  there exist integers $T_0,L_0,n_0$ such
 that
\begin{enumerate}
\item for every $n=hm > n_0$,  $0< p=p(n) < 1$, a graph $F \in {\cal F}(H,m)$ with $|V(F)|=n$
such that

\subitem{(a)}  $F$ is $H$-uniform
\subitem{(b)}  $F$ has no $(D,p)$-dense patches,  and
\item for every $0< p^*=p^*(n) <1$ and  ${\cal S}\subseteq \cee(F)$
such that $F$ (together with ${\cal S}$)  has no $(D^*,p^*)$-dense $H$-patches,
\end{enumerate}
there exists a refinement $\Pi_1$ of the initial
 partition ${\Pi}_0$ which is a
 $\delta$-regular, $(\eps(l),p)$-uniform, $(t,l)$-partition for some $l \leq L_0$ and $t \leq
 T_0.$
 \end{lemma}

\begin{remark}\label{R2}\rm
In both classic cases, $H=K_2$ and $H=K_3$, all copies of $H$ in $H^m$ are consistent, and consequently, the corresponding regularity lemmas are true for arbitrary graphs $F$, not just for bipartite or tripartite graphs $F\in{\cal F}(H,m)$. However, for general $H$ the ``partiteness" seems necessary.
Consider, e.g., $H$ being a path on vertices $1,2,3$, and let $F=K_{3m}$ be a complete graph on vertex set $V=V_1\cup V_2\cup V_3$. Further, let $c_P$ count {\it all} copies of $H$, while $s_P$ -- only those consistent with the partition $(V_1,V_2,V_3)$.
Then no regularity lemma may guarantee a $\delta$-regular partition of $F$.
Indeed,  every polyad $P$ with large $s_P$ must contain a bipartite subgraph $R\subset F[V_1,V_2]$ with still large $c_R$ but, clearly, with $s_R=0$, making $P$ $\delta$-irregular.
\end{remark}

\begin{remark}\label{R3}\rm Our lemma admits edge partitions besides the more traditional vertex partitions. Consequently, non-induced subgraphs are considered.
It is then  closer in spirit to the hypergraph regularity lemma  from \cite{FR1}, than the classical Szemer\'edi
lemma, where pairs of vertex partition classes determine  induced  subgraphs.
\end{remark}

\begin{remark}\label{R4}\rm
If both $p$ and
$p^*$ are constants, say 1, then the assumptions about the absence of dense patches are vacuously
satisfied with $D=D^*=1$. Indeed, e.g., by (\ref{density}) we have then $d_R=s_R/c_R\le1$ and no $(1,1)$-dense $H$-patches exist (see Definition \ref{nodensePa}).
 In the sparse case, however, the assumption of no $(D^*,p^*)$-dense $H$-patches, which plays a key role in the proof,
can be very hard to verify. For the application in this paper, i.e. for $H$ being a star forest $S$ with $\phi$
edges, and $p^*=p^{\phi}$, where $p=p(n)$ is as in the Setup on page \pageref{setup},
this is done in Section \ref{nodense} with considerable technical effort.
\end{remark}

\begin{remark}\label{R6}\rm
Using the relevant definitions, it is easy to see that the condition of
not having any $(D,p)$-dense patches
bounds $p$ from below: $p \ge 1/(n^2D)$;
similarly, as observed above (see (\ref{p^*-lower})), not having $(D^*,p^*)$-dense $H$-patches yields $p^* = \Omega(n^{-h})$.
\end{remark}

 \begin{remark}\label{R5}\rm
 We will use this lemma with $\eps(\ell) \le \ell^{-4\phi}$.
 The reason for having $\eps$ decrease with $\ell$ so quickly is that we will need
 $\eps$ to be much smaller than $1/\ell$, a lower bound on the average $p$-density of
 the resulting subgraphs (seeproof of Lemma~\ref{ManyPerfect}).
 \end{remark}

  \subsection{ Proof of the Subgraph Regularity Lemma}

 The proof is based on
 a technique from the original Szemer\'edi Regularity
 Lemma: refining the initial partition $\Pi_0$ until it becomes
 $\delta$-regular. This procedure will be monitored by the
 {\em index} of a partition. We will see that as long as a current
 partition is not $\delta$-regular, it can be refined, resulting in a
 partition with a substantially larger index (this is Lemma
 \ref{index-pump}). A priori bounds on
 the index given in Lemma \ref{cor13}  will guarantee that this process
  terminates successfully after a bounded number of steps.
Throughout,  the $(\eps,p)$-uniformity of the partitions will be maintained as well.

\subsubsection{Pumping the Index}
Let \begin{equation}\label{sigma}\sigma_R=\frac{c_R}{\kappa}\end{equation}
 be the {\it relative volume} of a subgraph $R$ of $F$. It
 measures what percentage of the copies of $H$ which can be found in $F$
 are contained in $R$. Note that in view of (\ref{kappa}),
\begin{equation}\label{one}
\sum_{P\in\calp_{\Pi}^{{\rm irr}}}\sigma_P=1.\end{equation}
Let $$\mbox{Index}(\Pi)=\sum\limits_{P\in \calp_{\Pi}} \sigma_Pd_P\log d_P$$ be {\it the index of a partition}
$\Pi$.
 Thus, the index is a weighted sum of a convex function ($x\log x$) of the $({\mathcal {S}},p^*)$-densities
 of the polyads consistent with $\Pi$.
 The merit of this convex function is that whenever a
 $(t,l)$-partition is {\em not} $\gd$-regular it is
 possible to refine it into a $(t',l')$-partition with a
 substantially larger index.

To start with, let $\Psi$ be {\it any} refinement of a current partition $\Pi$ as defined in Definition \ref{tl},
and, for a given polyad $P\in\calp_{\Pi}$, let
 $\ulam(P)=\ulam(\Psi,P)$ be the set of all polyads $R$ consistent
 with $\Psi$ and such that $R\subseteq P$.  Clearly,
  $$\calp_{\Psi} = \bigcup_{P \in \calp_\Pi} \ulam(P).$$
Setting $$\mu_R=\mu_R(P)=\frac{c_R}{c_P} \quad\left(=\frac{\sigma_R}{\sigma_P}\right),$$ and using the identity
$\sigma_R = \sigma_P \mu_R$ we have
 $$\mbox{Index}(\Psi)= \sum\limits_{P\in
 \calp_{\Pi}}\sum\limits_{R\in\ulam(P)} \sigma_R d_R\log d_R=\sum_P\sigma_P\sum_{R\in\ulam
 (P)}\mu_Rd_R\log d_R.$$
Observe that $d_P=\sum\limits_{R\in\ulam(P)}\mu_Rd_R$ and
 $\sum{\mu_R}=1$. Hence, by Jensen's inequality, for each $P\in\calp_{\Pi}$,
 $$\sum_{R\in\ulam(P)} \mu_Rd_R\log d_R\ge
 \left(\sum_R\mu_Rd_R\right)\log \left(\sum_R\mu_Rd_R\right)
 =d_P\log d_P,$$
 and consequently
 $\mbox{Index}(\Psi) \geq \mbox{Index}(\Pi).$

The above application of Jensen's inequality can be sharpened significantly in the case when $P$ is a
$\delta$-irregular polyad, provided we impose an extra constraint on partition $\Psi$. Given a subgraph $Q$ of
$F$, we say that a partition $\Psi$  {\it refines} $Q$ if its edge partition (\ref{edge}) refines the partition
$(Q,F-Q)$, that is, every subgraph from (\ref{edge}) is either contained in or disjoint from $Q$.
\begin{lemma}[Defect Lemma]
\label{defect}
 Let $P\in\calp_{\Pi}^{{\rm irr}}$ be a $\delta$-irregular polyad and  $Q\subseteq P$ be such that
 \begin{list}{}{}
 \item[(i)] $c_{Q}\ge\delta c_{P}$ and
 \item[(ii)] $|d_{P}-d_{Q}|>\delta$.
 \end{list}
 If $\Psi$ is a refinement of $\Pi$ which refines $Q$ then
 $$\sum_{R\in\ulam(P)} \mu_Rd_R\log d_R\ge d_P\log d_P+\frac{2\delta^4}{d_P}.$$
 \end{lemma}
 For the sake
 of readability, we defer for a moment the technical proof of Lemma \ref{defect}. An immediate consequence of Lemma \ref{defect} is a jump in the value of the index when a partition has many $\delta$-irregular polyads and  its refinement is properly
 chosen. A subgraph $Q$ which satisfies conditions (i) and (ii) of Lemma \ref{defect} with respect to a given $\delta$-irregular polyad $P$ will be called {\it a witness of irregularity for $P$}.

 \begin{lemma}[Index Pumping Lemma]
 \label{index-pump}
 If  $\Pi$ is a $\delta$-irregular $(t,l)$-partition of $F$
 and $ \Psi$ is a refinement of $\Pi$ which also refines at least one  witness of irregularity for each  polyad
$P\in\calp^{{\rm irr}}$, then
 $$
{\rm Index}(\Psi) \ge {\rm Index}(\Pi) + \frac{\delta^5}{D^*}.$$
 \end{lemma}
\nin{\bf Proof:}  As before, but using  Lemma \ref{defect} instead of Jensen's inequality,
\begin{eqnarray*}
 \mbox{Index}(\Psi) &=& \sum_P\sigma_P\sum_{R\in\ulam
 (P)}\mu_Rd_R\log d_R\\
&\ge& \sum_P\sigma_P d_P\log d_P +
 \sum_{P \in \calp^{{\rm irr}}}
 \sigma_P\frac{2\delta^4}{d_P}\ \\
&\ge& \mbox{Index}(\Pi) + \sum_{{P\in \calp^{{\rm irr}} },\ {c_P\ge\kappa/\log^2n}}
 \sigma_P\frac{2\delta^4}{D^*},
 \end{eqnarray*}
For the last inequality above we used the assumption that there are no $(D^*,p^*)$-dense $H$-patches in $F$, and thus, if $c_P\ge\kappa/\log^2n$ then $d_P\le D^*$. By
the $\delta$-irregularity of $\Pi$ (see Definition \ref{dfn4.4}) and by (\ref{sigma}), we have
 $$ \sum_{P\in\calp^{\rm irr}} \sigma_P > \gd.$$
 Moreover, for sufficiently large $n$, recalling that $|\calp|\le s=t^hl^{|E(H)|}$,
 $$\sum_{P: c_P< \kappa/\log^2n}\sigma_P<\sum_P\log^{-2}n\le s\log^{-2}n<
 \frac\delta 2, $$
 because $t,l$ and $H$ do not depend on $n$. Hence, in view of (\ref{one}),
$$ \sum_{{P\in \calp^{{\rm irr}} },\ {c_P\ge\kappa/\log^2n}}
 \sigma_P\frac{2\delta^4}{D^*}\ge\frac {\delta^5}{D^*}.$$
 \hfill $\Box$
\bigskip

 The Index Pumping Lemma works hand in hand with the following estimates.
 \begin{lemma}\label{cor13}
 For all partitions $ \Pi$ of $F$,
$$-\frac{1}{e} \le
 \mbox{\rm Index}({\Pi})\le D^*\log D^*
 + o(1).$$
 \end{lemma}
 \nin{\bf Proof:}
 The lower bound follows easily
 from (\ref{one}) and  the fact that the function $x\log x$ achieves its minimum
 at $x=1/e$.

 For the upper bound, we split
\begin{eqnarray*}
 \mbox{Index}({\Pi})=\sum\limits_{P:c_P\ge \kappa/{\log^2n}}
 \sigma_Pd_P\log d_P +
 \sum\limits_{P:c_P< \kappa/{\log^2n}} \sigma_Pd_P\log d_P.
\end{eqnarray*}
Because $F$ has no $(D^*,p^*)$-dense $H$-patches, we can bound the first term  by
  $  D^*\log D^*$.
 As for the second term,  by Proposition \ref{(b)},
$c_P < \frac{\kappa}{\log^2n}$ implies that $s_P \le \frac{(1+o(1))D^*p^*\kappa}{\log^2n}$. Hence,
 $$\sigma_Pd_P = s_P/(\kappa p^*)\le(1+o(1))D^*/(\log^2n).$$

\nin Moreover,  for all polyads $P$ we have $\log d_P=O(\log n)$,  because by (\ref{p^*-lower})
(or see Remark \ref{R6})
 $$d_P=\frac{s_P}{p^*c_P}\le \frac{1}{p^*}=O(n^{h}).$$
Hence,
$$\sum\limits_{P:c_P< \kappa/{\log^2n}} \sigma_Pd_P\log d_P\le s
 \frac{D^*}{\log^2n}O(\log n)=o(1),$$
also using the fact that the number of polyads is at most $s=t^h l^{|E(H)|} = O(1)$.
 \hfill $\Box$

 \begin{remark}\rm Note that we would not have been able to bound
 the second term  above as easily, had we used
 the more traditional functional, $\sum_P d_P^2$, in the definition of index.
The index function was introduced by Szemeredi in [Sz] as the main tool in
the proof of his regularity lemma (Theorem 4.5 above), and later used also
in [FR] to prove the  hypergraph regularity lemma (Theorem 4.11 above).
The idea of replacing $d^2$ by any convex function in this
context was suggested by Ajtai, Komlos, Szemeredi
[E. Szemeredi, oral communication, 1978]. Such an alternative index
function first appeared in print in \cite{convex}.
\end{remark}

 The above two lemmas basically complete the proof of Lemma \ref{reglem}, provided
 we can construct a required refinement which in addition is $(\eps(l),p)$-uniform.
 Then, iterating the process of refining a current
 partition, we must reach a $\delta$-regular
 $(t,l)$-partition. The bounds in these lemmas
 imply an upper bound on the number of iterations, solely in terms of $D^*$ and $\delta$.

In order to conclude the existence of constants $T_0$ and  $L_0$ such that $t \leq T_0$ and $l \leq L_0$, we need
to show how the refinement is created. Likewise, although the asymptotic estimates we encountered so far
 contribute toward the value of $n_0$,
 the most demanding constraint on $n_0$ comes from the refinement which invoke  Theorem \ref{SSRL}.
  A description of how the refinement is constructed
 is deferred to the next subsection.
 We complete this subsection by giving a proof of the Defect Lemma.

\bigskip
 \noindent{\bf Proof of Lemma \ref{defect}:}
Let $\ulam(Q)$ be the set of all polyads $R$ consistent
 with $\Psi$ and such that $R\subseteq Q$.
 Recall that $\mu_R=c_R/c_P$ and set
 $$\mu=\sum_{R\in\ulam(Q)}\mu_R.$$
Since $\Psi$ refines $Q$, we have $\sum_{R\in\ulam(Q)}c_R=c_Q$ and $\sum_{R\in\ulam(Q)}s_R=s_Q$. Thus, by (i),
$\mu=c_Q/c_P \ge \gd>0$ and
$$d_Q=\frac1\mu\sum\limits_{R\in\ulam(Q)}
 \mu_Rd_R.$$
Moreover, $\mu=1$ then $c_Q=c_P$, $s_Q=s_P$, and consequently, $d_P=d_Q$, contradicting (ii).
Hence, we may assume that $\mu<1$, and by Jensen's inequality,
 $$\sum_{R\in\ulam(Q)}\mu_Rd_R\log d_R\ge \mu\cdot d_Q\log d_Q$$
 and
 $$\sum_{R\in\ulam(P)\backslash \ulam(Q)} \mu_Rd_R\log d_R\ge
 (d_P-\mu d_Q)\log\left(\frac{d_P-\mu d_Q}{1-\mu}\right).$$
 So,  setting $d=d_P$, $d'=d_Q$, and $\displaystyle d''=\frac{d-\mu d'}{1-\mu}$ for convenience,
 \begin{equation}
 \label{index1}
 \sum_{R\in\ulam(P)} \mu_Rd_R\log d_R\ge \mu d'\log d'+ (1-\mu)
 d''\log d''.
 \end{equation}
 Consider the function
 $$f(x)=x\log \left(\frac x\mu\right)+(d-x)\log \left(\frac{d-x}{1-\mu}
 \right).$$
 Clearly, $f(\mu d) = d \log d$, and further,
 since $\mu d' + (1-\mu)d'' = d$,
 \begin{eqnarray}
 f(\mu d')& = &
 \mu d'\log d' + (d- \mu d')\log\left(\frac{d-\mu d'}{1-\mu}\right)\nonumber\\
 &=& \mu d'\log d' + (1- \mu) d''\log d''.
 \end{eqnarray}
 Thus in view of (\ref{index1}),
 \[\sum_{R\in\ulam(P)} \mu_Rd_R\log d_R - d\log d \ge
 f(\mu d') - f(\mu d),\]
 and we will be done if
 \begin{equation}
 \label{changeinf}
 f(\mu d') - f(\mu d) \ge \frac{\delta^4}{2d}.
 \end{equation}
 Towards this, compute
 $$f'(x)=\log\left(\frac x\mu\right) -\log \left(\frac{d-x}{1-\mu}\right),$$
 and
 $$f''(x) = \frac{1}{x} + \frac{1}{d-x},$$
 and note that $f'(\mu d)=0$, and that $f''(x) \ge 4/d$, for all $0<x<d$. Thus, by
 the Lagrange formula for the remainder,
 we have, for some $0<\mu d < x_0 < \mu d' < d$,
 \begin{eqnarray*}
 f(\mu d') - f(\mu d)  =  f''(x_0) \frac{(\mu d'-\mu d)^2}{2}
 \ge  \frac{4}{d} \frac{\delta^4}{2} = \frac{2\delta^4}{d},
 \end{eqnarray*}
 since $\mu \ge \delta$ and $|d'-d|>\delta$.
 This proves (\ref{changeinf}), completing the proof of Lemma~\ref{defect}.
 \hfill $\Box$

 \subsubsection{Refining a Partition}
 We will now finish the proof of Lemma \ref{reglem}. As already observed prior to the proof of Lemma \ref{defect}, it remains to construct
 a
 $(t',l')$-partition $\Psi$ with the following properties. Recall that a subgraph $Q$ which satisfies conditions (i) and (ii) of Lemma \ref{defect} with respect to a given $\delta$-irregular polyad $P$ is called {\it a witness of irregularity for $P$}. The new partition $\Psi$ must
\begin{enumerate}
\item
be a refinement of
 a current $\delta$-irregular $(t,l)$-partition $\Pi$,
\item
refine at least one witness of irregularity for each irregular polyad of $\Pi$,
\item
be $(\eps(l'),p)$-uniform.
\end{enumerate}
  Recall that a partition $\Pi$ gives rise to at most $s=t^hl^{|E(H)|}$ different polyads.
 Suppose that
 polyads $P_1,\ldots, P_{s'}$ are $\delta$-irregular with witnesses of their irregularity $Q_1,\ldots, Q_{s'}$.
The requirement that the new partition must refine each $Q_g$, $g=1\dots s'$, is easy to fulfill.
 Let Venn$\{Q_1,\ldots, Q_{s'}\}=\{R_1,\ldots, R_{s''}\}$,
 $s''\le 2^{s'}$, be the family of all graphs of the form $\bigcap\limits^{s'}_{g=1}
 Q^{\ep_g}_i$, for $\ep_g \in \{0,1\}$, where
 $$Q^{\ep_g}_i=\cases{Q_g & $\ep_g=1$\cr
 F-Q_g & $\ep_g=0$.}$$
 The graphs $R_g, g=1 \ldots s''$, partition (by intersection) every
 subgraph $F_{xy}^{i,j,k}$ into smaller subgraphs  $F_{xy}^{i,j,k,g} =F_{xy}^{i,j,k} \cap R_g$, that is,
 $$F_{xy}^{i,j,k}=\bigcup_{g=1}^{s''} F_{xy}^{i,j,k,g}.$$
This  new partition ${F_{xy}^{i,j,k,g}}$, clearly, refines  every subgraph $Q_g$, $g=1\dots s'$.
 Note also that the previous partition of every tube into at most $l$ subgraphs
 is now
 refined into at most $l'$ finer subgraphs, with $l' = ls''\le l2^{s'}\le l2^{s}$.

 So far we have only partitioned
{\it the edge sets} of $\Pi$. With the help of Theorem \ref{SSRL} (the Sparse Regularity Lemma) we will now refine
the  vertex equipartition of $\Pi$ to  ensure that the new equipartition is $(\eps(l'),p)$-uniform (and still refines
every $Q_g$).
By the assumption on $F$, none of the resulting graphs $F_{xy}^{i,j,k,g}$  has
 $(D,p)$-dense patches ($F_{xy}^{i,j,k,g}$ is a subgraph of $F$). Thus, we may apply  Theorem \ref{SSRL} with $\eps=\eps(l')$, $D=2$, $t_0=t$ and $r=|E(H)|t^2l2^s$, to the graphs $\{F_{xy}^{i,j,k,g}\}$ with the vertex equipartition given  by ${\Pi}$, where $xy\in E(H)$, $i,j=1\dots t$, $k=1\dots l(xy,i,j)\le l$ and $g=1\dots s''$. This way we obtain  a finer equipartition of the vertices
 which together with the described above partition of the subgraphs, defines an $(\eps(l'),p)$-uniform $(t',l')$-partition $\Psi$. This new partition  consists of the following parts:
 \begin{itemize}
 \item Vertex sets defined by the equipartition given by Theorem \ref{SSRL}:
 $$\{\overline V_{x}^{i}, i\in [ t']\}\prec \{V_{x}^{i}, i\in
 [t]\},$$
 where $\prec$ denotes refinement.
 \item Edge sets
 $$\overline{F}_{xy}^{i,j,k,g} = F_{xy}^{i,j,k,g}[\overline{V}_x^i,\overline{V}_y^j].$$
 \end{itemize}
The new parameter $t'$ is bounded by $T_0$ of Theorem \ref{SSRL} with the above inputs, and so depends on $l$,
$t$, $H$ and the function $\eps(\ell)$.
Starting from $t=l=1$ and iterating for a bounded number of times, as determined in the previous subsection,
yields the promised constants $T_0$ and $L_0$ of Lemma \ref{reglem}. We do not attempt to compute any of these
constants explicitly.
This completes the proof of Lemma \ref{reglem}.

 \subsection{No Dense Patches}\label{nodense}
In the previous section we proved Lemma \ref{reglem}. Now we are going to verify that it can be applied in the
setting relevant to the proof of Lemma \ref{wrapup}.
 Let us summarize again the state-of-art at the end of Section \ref{TandC}. As a consequence of assumption (\ref{assumption}) of Theorem \ref{allweneed} we have selected a star forest $S$ consisting of $\nu$ stars with, respectively, $\phi_1\dots \phi_{\nu}$ arms.
 We set
$$|E(S)|=\phi=\phi_1+\ldots+\phi_{\nu}$$
and so $|V(S)|=\nu+\phi$. For clarity, we will be further assuming that $n$ is divisible by $\nu+\phi$. Under this assumption, Lemma~\ref{equi} provided us with
 a partition of the
 vertices of $G$ into sets of equal size $m=n/(\nu+\phi)$, further denoted
\[\pi_0= \{V_1 \dots  V_\nu, W_{11}
 \dots W_{1\phi_1} \dots W_{\nu,1} \dots W_{\nu,\phi_{\nu}}\}.\]
Based on this partition, we have focused on  the set $\cee$ of all constellations, that is all copies of $S$ in
$G$ that are
 consistent with $\pi_0$. Among the constellations we have distinguished
 the set $\mathcal {S}$  of special constellations, related to the special role of certain additions of a copy of $M$ to $G$.

Recall that $M$ is an arbitrary (balanced) graph with $\rho(M)=2$ with vertices labeled $x_1\dots x_\nu$, and that
$\widehat M$ is a graph built upon $M$ and consisting of a number of tepees spreading above the edges of $M$. It
is $\widehat M$ which gave rise to the star forest $S$ by deleting one leg of each tepee (seeThe Missing Leg
Property -- Observation \ref{miss}). The vertices of $S$ not belonging to $M$ (i.e. the tips of tepees) are
denoted by $u_{ab}$, where $a=1\dots\nu$ and $b=1\dots\phi_a$, in such a way that $u_{ab}$ are the neighbors of
$x_a$. The graph $G$ is an arbitrary member of family $\G$ defined in Definition \ref{defG}. For all $a=1\dots\nu$
and $b=1\dots\phi_a$, set $F_{ab}=G[V_a,W_{ab}]$ and $$F= \bigcup_{ab}F_{ab}.$$ A copy of $S$ in $F$ is consistent
if each $x_a$ is mapped onto a vertex of $V_a$ and each $u_{ab}$ is mapped onto a vertex of $W_{ab}$.

Note that $F\in{\cal F}(S,m)$. We want to apply the Subgraph Regularity Lemma \ref{reglem}, with $H=S$, to the
graph $F$, family ${\cal S}$, and sequences: $p=p(n)$ as in the Setup preceding Lemma~\ref{wrapup} and
$p^*=p^\phi$. This normalizing factor $p^*$ is to be expected because there are $\phi$ edges that must appear in
$G$ to turn a
 constellation into a special constellation (the missing legs of
 the tepees), and
 informally (if $G$ were random), each edge would appear with probability~$p$.

In order to verify the assumptions of Lemma \ref{reglem} for $F$, we have to show that for some $D>0$ and $D^*\ge1$
(independent of $G$, and thus independent of $S$ and $\phi$ in particular),
\begin{enumerate}
\item[1(a):]  $F$ has no $(D,p)$-dense patches
\item[1(b):] $F$ is $S$-uniform
\item[2:] $F$, together with given $\cal S$, has no $(D^*,p^\phi)$-dense $S$-patches.
\end{enumerate}
All these properties of $F$ will be derived from  Definition \ref{defG} of family $\G$ to which $G$, a
supergraph of $F$, belongs. Property 1(a) follows immediately with $D=2$ from Property (P5) of $\G$, while property 1(b) is an easy consequence of Property (P3), that is, of the  bounds  on the vertex degrees in $G$. Indeed,
by Example \ref{cee}, there are
$$\kappa\sim m^{\nu+\phi}p^\phi$$
 (consistent) copies of $S$ in $F$, while a given edge belongs to at most
\begin{equation}\label{atmost}
n^{\nu-1}(2pn)^{\phi-1}=O\left(\frac\kappa{n^{3/2}}\right)=o\left(\frac\kappa{\log^3n}\right)
\end{equation} of them, for large $n$. The next lemma establishes assumption 2 with $D^*=2^q+1$, where, recall,
$q=\ceiling{10C^2\nu/\ga}$ is an upper bound on $\phi$ (see Lemma \ref{common}).
Recall also that $c_R$ is the number of constellations (i.e., consistent with $\pi_0$ copies of $S$) contained in a
subgraph $R$.
\begin{lemma}\label{nodensepa}
For every subgraph $R\subseteq F$ with $c_R\ge \kappa/\log^2n$ we  have $d_R\le 2^q+1$.

\end{lemma}

Let us set $U_a=V(R)\cap V_a$, $R_{ab}=R\cap F_{ab}$ and  ${\rm deg}_{ab}(v)={\rm deg}_{R_{ab}}(v)$ for all $a=1\dots\nu$, $b=1\dots\phi_a$. It is straightforward to express $c_R$ in terms of the degrees ${\rm deg}_{ab}(v)$. Indeed, we have
\begin{equation}\label{key}
c_R= |U_\nu|\prod^{\nu-1}_{a=1} \left(\sum_{v\in U_a}
 \prod^{\phi_a}_{b=1}{\rm deg}_{ab}(v)\right).
\end{equation}
This simple formula plays a key  role in our proof. First, it trivially implies that
\begin{equation}\label{cRlow}
c_R \ge \prod\limits^{\nu}_{a=1}\left( |U_a|\prod\limits^{\phi_a}_{b=1}\delta_{ab}\right),
\end{equation}
 where $\delta_{ab}=\min_v{\rm deg}_{ab}(v)$.
 Second, in connection with the estimate of $\kappa$, (\ref{key}) yields lower bounds on $|U_{a}|$ and $\Delta_{ab}=\max_v{\rm deg}_{ab}(v)$ in terms of $c_R$. We will use them later in the proof.

 \begin{prop}\label{prop16}
 For all $1 \leq a \leq \nu$,
 \begin{list}{}{}
 \item[(a)] $\displaystyle |U_a|=\Omega\left(\frac{c_R}{\kappa}  n\right)$, and
 \item[(b)]  $\displaystyle \Delta_{ab}=\Omega\left(
 \frac{c_R}{\kappa}  np\right)$ for each $b=1\dots\phi_a$.
 \end{list}
 \end{prop}

 \nin{\bf Proof.} Using Property (P3) of $G$, we can further estimate (\ref{key}) as follows. For each $1\le a\le \nu$ and $1\le b\le \phi_a$ we have
 \begin{eqnarray*}
 c_R&=& |U_v|\prod^{\nu-1}_{a'=1} \left(\sum_{v\in U_{a'}}
 \prod^{\phi_{a'}}_{b'=1}{\rm deg}_{a'b'}(v)\right)\\
 &\le & \cases{m\prod\limits^{\nu-1}_{a'=1} \left(\sum\limits_{v\in
 U_{a'}}(2np)^{\phi_{a'}}\right) =O(|U_{a}|n^{\nu-1}(np)^\phi)\cr
  m^{\nu}\prod\limits^{\nu-1}_{a'=1}\left(\prod\limits^{\phi_{a'}}_{b'=1} \Delta_{a'b'}\right)
 =O\left( n^{\nu}(np)^{\phi-1}
 \Delta_{a b} \right).}
 \end{eqnarray*}
By Example \ref{cee}, we have $\kappa\sim m^{\nu+\phi}p^\phi$, which completes the proof.
 \hfill $\Box$

\bigskip
In view of (\ref{cRlow}) to prove Lemma \ref{nodensepa} we need to prove a matching upper bound on $s_R$. However, our current proof techniques allow us only to set a bound of the form
\begin{equation}\label{sRup}
s_R \le (1+o(1))p^{\phi}\prod\limits^{\nu}_{a=1}\left( |U_a|\prod\limits^{\phi_a}_{b=1}\Delta_{ab}\right),
\end{equation} (see later in this section). So, we can only obtain the desired bound on $d_R$ if $\delta_{ab}$ and $\Delta_{ab}$ are at most a factor of two  apart. This leads to the following definition. A subgraph $R$ of $F$ will be called {\it degree sandwiched} if there exist numbers $r_{ab}$, $a=1\dots\nu$, $b=1\dots\phi_a$, such that for all $v\in U_a$, we have $r_{ab}\le{\rm deg}_{ab}(v)\le2r_{ab}$.

 To put ourselves in a more comfortable position, we will cover $R$ with its degree sandwiched subgraphs in such a way that each constellation in $R$ is contained in exactly one of them.
  Remembering that $\phi\le q$, the following lemma will imply Lemma
\ref{nodensepa}.
\begin{lemma}\label{sand}
For every degree sandwiched subgraph $R\subseteq F$ with $c_R\ge \kappa/\log^{\phi+3}n$ we  have $d_R\le
2^{\phi}+o(1)$.

\end{lemma}
In order to see how Lemma \ref{sand} implies Lemma \ref{nodensepa} we need one more fact.
\begin{prop}\label{new} For every degree sandwiched subgraph $R$ with $c_R<\kappa/\log^{\phi+3}n$ there exists a degree sandwiched subgraph $R'$ such that $R\subset R'\subseteq F$ and
$$\kappa/\log^{\phi+3}n\le c_{R'_i} = (1+o(1))\kappa/\log^{\phi+3}n.$$
\end{prop}
\nin{\bf Proof.}
We will construct a sequence of subgraphs $R=R^1\subset R^2\subset \cdots \subset F$
such that each $R^j$ is degree sandwiched and is obtained from $R^{j-1}$ by adding at most $2np=O(\sqrt n)$ edges. Note that by Lemma \ref{equi}, $F$ itself is degree sandwiched between $\frac34mp$ and $\frac32mp$.
By (\ref{atmost}), the increments in $c_{R^j}$ are at most $O(\kappa/n)=o(k/\log^{\phi+3}n)$, and so there must be some $j$ such that $R'=R^j$ is the desired subgraph. (Note that it is certain that such an index $j$ will be found many steps before reaching $F$. Indeed, by (\ref{cRlow}),
we will be done when,  for example, all $U_a$ will be of order $\Omega(n)$ and all degrees of order at least $np/\log n $.)

Our construction is performed in two phases. In phase one, for each $a$ and $b$, we turn $R\cap F_{ab}$ into an induced subgraph of $F_{ab}$, and then in phase two we will add one by one all remaining vertices of $\bigcup_a V_a$. A single step of phase one consists of fixing a vertex $v$ in $U_a$ of minimum degree $r$ in $R^{j-1}\cap F_{ab}$ and adding some $\min(r,{\rm deg}_{F_{ab}}(v)-r)$ edges $vu\in F_{ab}\setminus R^{j-1}$. This way the process maintains the ``degree-sandwichness'' throughout, and is bound to arrive at $R^j=F_{ab}[U_a]$ for some $j$. A single step of phase two consists of adding a vertex $v\in \bigcup_a V_a\setminus V(R^j)$ together with all its neighbors in $F$. \hfill $\Box$
\bigskip

\nin{\bf Proof of Lemma \ref{nodensepa}:}
 Let $R\subseteq F$ with $|C(R)|=c_R\ge \kappa/\log^2n$.
Without loss of generality assume that for all  $a=1\dots\nu$, $b=1\dots\phi_a$, and all  $v\in U_a$, we have ${\rm deg}_R(v,W_{ab})>0$, and define for each
$j=0\dots \lfloor\log \Delta_{ab}\rfloor$
$$U_{ab}^j=\left\{v\in U_a: \frac{\Delta_{ab}}{2^{j+1}}\le {\rm deg}_{ab}(v)\le\frac{\Delta_{ab}}{2^m}\right\}.$$
 This
 subdivides every set $U_a$  into $O(\log^{\phi_a}n)$ classes
$$U_a^{j_1\dots j_{\phi_a}}=\bigcap_{b=1}^{\phi_a}U_{ab}^{j_b}$$
 so that the degrees in $R$ of all vertices from one class, along every
 subgraph $F_{ab}$ are sandwiched between a number and its double.
By selecting one partition class from each set $V_a$ and taking all the edges of $R$ with endpoints in these
classes, this yields a covering of  $R$
 by $O(\log^\phi n)$  subgraphs $R_1, R_2,\ldots$ which are all degree sandwiched. Set $d_i=d_{R_i}$, etc., for convenience. Clearly, we have $c_R=\sum_ic_i$, $s_R=\sum_is_i$, and
 $$d_R=\sum\limits_{i}\frac{s_i}{p^\phi c_R}=\sum_{c_i\ge \kappa/\log^{\phi+3}n}d_{i}\frac{c_i}{c_R}+\sum_{c_i< \kappa/\log^{\phi+3}n}\frac{s_i}{p^\phi c_R}.$$
Applying Lemma \ref{sand} in the first sum, we bound it by $2^\phi+o(1)$. For the second sum, we use a similar idea
as in the proof of Proposition \ref{(b)}. For each $i$, let $R_i'$ be the degree sandwiched supergraph of $R_i$ as in Proposition \ref{new}.
Applying Lemma \ref{sand} to $R'_i$, we obtain
  $$s_i \leq s_{R'_i} \leq (2^\phi+o(1)) p^* c_{R'_i}\le(1+o(1))\frac{ 2^\phi p^* \kappa}{\log^{\phi+3}n}$$
and so
$$\sum_{c_i< \kappa/\log^{\phi+3}n}\frac{s_i}{p^\phi c_R}=O\left(\frac{\kappa\log^\phi n}{c_R\log^{\phi+3}n}\right)=O(\log^{-1}n)=o(1),$$
by the assumption that $c_R\ge\kappa/\log^2n$. \hfill $\Box$

%
%

\bigskip
The rest of this section is devoted to proving Lemma \ref{sand}. Let $R$ be a degree sandwiched subgraph of $F$
with $c_R\ge \kappa/\log^{\phi+3}n$. As for its vertex set, we have,  for each $a$ and $b$, $V(R)\supset W_{ab}$.
 and $R_{ab}=R\cap F_{ab}$. The degrees in $R$ are such that for all $a$ and $b$,
 and all $v \in U_a$
 $$r_{ab}\le{\rm deg}_{R_{ab}}(v)\le 2r_{ab},$$
 for some $r_{ab}>0$, where recall ${\rm deg}_{R_{ab}}(v)={\rm deg}_R(v,W_{ab})={\rm deg}_{R_{ab}}(v)$.

 In view of (\ref{cRlow}), to prove
 Lemma~\ref{sand}, it suffices to show that
\beq {numb} s_R\le(1+o(1))
 p^\phi  2^\phi \prod^{\nu}_{a=1}\left( |U_a| \prod^{\phi_a}_{b=1} r_{ab}\right).
\enq
The following corollary of Proposition \ref{(b)} will be of some use.
\begin{cor}\label{lowbounds} If $c_R\ge \kappa/\log^{\phi+3}n$, then for all $1\le a\le \nu$,
 \begin{list}{}{}
 \item[(a)] $\displaystyle |U_a|\ge\frac n{\log^{\phi+4} n}$, and
 \item[(b)]  $\displaystyle r_{ab}\ge
 \frac{np}{\log^{\phi+4}n}$ for each $b=1\dots\phi_a$.
 \end{list} \hfill $\Box$
\end{cor}
 Let us now briefly outline the rest of the proof of Lemma~\ref{sand}: we want to bound the number of special constellations in
 $R$, or more specifically we want to prove~(\ref{numb}), where $U_a=V(R)\cap V_a$. Recall from Section \ref{TandC} that every special constellation is
 of the form $S(X)$ for a set of vertices $X \in \X_3$ which contains one vertex from each set $V_a$.
 Such a constellation corresponds via the Missing Leg Property (Observation \ref{miss})
 to a set of $\phi$ edges of $G$ (the missing legs.) These edges, together with the edges of the constellation itself form a set of $\phi$ tepees over the edges of
  $M_X$, a copy of $M$ planted at $X$ (see Figure 3).

 We will construct a $\nu$-partite auxiliary graph $A$
 on the vertex sets $U_1\dots U_\nu$ that actually contains these graphs $M_X$ for which $S(X)\subset R$, but possibly  also other copies of $M$.
 Hence, the number of consistent copies
 of $M$ in $A$ (such that $x_a$ is mapped  onto a vertex of $U_a$) will be an upper bound on the number of special
 constellations in $R$.

More precisely, for  every  pair $a,a'$ for which $x_ax_{a'}\in E(M)$ define the set  $\Gamma_{aa'}$  of indices
 $b$ such that $x_au_{ab}x_{a'}$ is a tepee in $\widehat M.$
 Formally, using notation $a'(a,b)$ from Section \ref{TandC},
 \beq
 {gam}
 \Gamma_{aa'} = \{ b : a'=a'(a,b) \}.
 \enq
We will draw an edge in $A$ joining $v\in U_a$ with $v'\in U_{a'}$,
 if for each $b\in\Gamma_{aa'}$ the vertex $v'$ is connected by an edge of $G$ with a vertex of $W_{ab}$ which is a neighbor of $v$ in $R$.
 If each bipartite subgraph $A_{aa'}=A[U_a,U_{a'}]$ was highly regular, then using Proposition \ref{graphcount}
 it would be
 easy to count the number of copies of $M$ in $A$.

To achieve this high regularity of $A_{aa'}$ we will consider its supergraph $ A'_{aa'}$ with the entire $V(G)$ in
place of $U_{a'}$. Then the degree of $v\in U_a$ in $ A'_{aa'}$ will be asymptotically determined, via Property
(P6) of $G$, by the sizes of its neighborhoods in $W_{ab}$ for all $b\in\Gamma_{aa'}$. We know that for distinct
$v$ they are all sandwiched between $r_{ab}$ and $2r_{ab}$, but may vary, obstructing the desired regularity. To
remedy this problem we will enlarge these neighborhoods so that they are all equal $2r_{ab}$. Since we still want
to keep the co-degrees fairly low (as they are by Property (P4) of $G$), we will make these enlargements randomly.
Only after this is done, we will formally define the auxiliary graph, prove that it is highly regular and
finally count copies of $M$ in it.

 Fix $a,b$ and set $r=r_{ab}$. Let $k=|U_a|$, $U_a=\{v_1\dots v_k\}$, and
$N_i=N_i(b)$ be the neighborhood of $v_i$ in $R_{ab}$. Note that $r\le |N_i|\le2r$ and,
by Definition \ref{defG} (P4),  for all $i\neq j$, $|N_i \cap N_j| \leq 3 \log n.$ Note also that by Corollary \ref{lowbounds} we have, say, $r\ge n^{0.49}$. The following lemma gives us exactly what we need.
 \begin{lemma}
 \label{enlarge} For a sufficiently large integer $n$, let $N_1,\dots,N_k$, $k\le n$, be a set system of subsets of $W$, where $\theta_1n\le|W|\le n$,  for all $i$, $ n^{0.49}\le r\le |N_i|\le 2r\le\theta_2\sqrt n$, and
for all $i\neq j$, $|N_i \cap N_j| \leq 3 \log n.$
Then there exist sets $\overline N_i=\overline N_i$, $i=1\dots k$, such that for all $i$
 \begin{list}{}{}
\item[(a)] $N_i\subseteq \overline N_i$
 \item[(b)] $|\overline N_i|=2 r$
 \item[(c)]
  for all $i\neq j$, $|\overline N_i\cap \overline N_j| \leq  12 \log n.$
 \end{list}
 \end{lemma}
 \nin{\bf Proof: }
  Pick $\widehat
 N_1\dots \widehat N_k$ randomly, uniformly and independently from
 $[W]^{2r}$, the family of $2r$-element subsets of $W$. Then with probability $1-o(1)$,
 \beq
 {int}
 |\widehat N_i\cap \widehat N_j|\le 3\log n,\quad \forall\ i<j\quad\mbox{ and}\quad
 |\widehat N_i\cap N_j|\le 3\log n\quad \forall\ i,j.
 \enq
   This is because for a fixed
 set $N\subset W$, and a random set
 $\widehat N\subset W$, $|\widehat N|=2r$, the random variable $Z=|N\cap \widehat N|$ has
 hypergeometric  distribution  with expectation $|N||\widehat N|/|W|
 \leq 4\theta_2^2/\theta_1$. By Chernoff's inequality (e.g., (2.11) and Thm.~2.10 from
 \cite{JLR}),
 $${\rm Pr}(Z>3\log n)<e^{-3\log n}=n^{-3}.$$
 As there are fewer than $2n^2$ pairs of sets $(\widehat N_i,\widehat N_j)$,
 $i<j$, and $(\widehat N_i, N_j)$, the probability that at least one
 pair will violate (\ref{int}) is less than $\frac{2n^2}{n^3}\to 0$.
 From (\ref{int}), clearly,
 $$|(N_i\cup \widehat N_i)\cap (\widehat N_j\cup N_j)|\le 12\log n.$$
 Since $3r-3\log n>2r$ for large $n$, the required sets $\overline N_1\dots \overline N_k$ can now be chosen
 arbitrarily in such a way that $|\overline{N}_i|=2r$ and that
 $N_i\subseteq \overline N_i\subseteq N_i\cup\widehat N_i$.
 Then
 \[\overline{N}_i \cap \overline{N}_j \subseteq (N_i \cup \widehat{N}_i) \cap
 (N_j \cup \widehat{N}_j), \]
 and thus (c) holds as well.
 \hfill $\Box$


\begin{figure}[htb]
\begin{center}
    \input{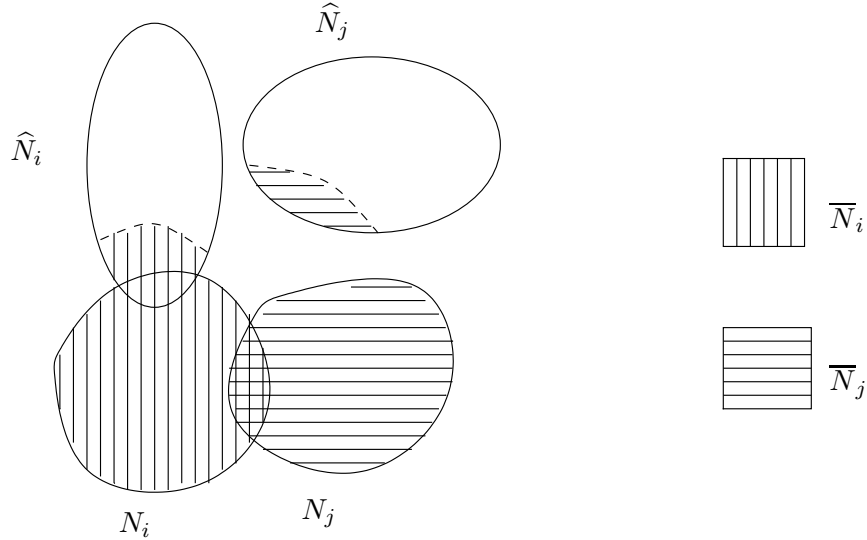}
      \caption{ Enlarging $N_i$'s}
\end{center}
\end{figure}

\bigskip
Let $\widehat N_i(b)$ be the extentions of neighborhoods $N_i(b)$ of vertices $v_i\in U_a$ in $W_{ab}$ whose existence is guaranteed by Lemma \ref{enlarge}.
 We now define the auxiliary graphs.
 For every $aa'$ such that $x_ax_{a'} \in E(M)$ let
 $  A'_{aa'}=(U_a,V(G), E'_{aa'})$,
 where for $v_i \in U_a, v' \in V(G) $:
 $$v_iv'\in   E'_{aa'}\Leftrightarrow |\overline{N}_ i(b)\cap
 N_G(v')|\geq 1\
  \mbox{for all} \ b \in \Gamma_{aa'}.$$
 \begin{remark}\rm
 Since $U_a \subseteq V(G)$, the above definition involves a slight
 abuse of notation. Formally we are creating
 a graph between a copy of $U_a$ and $V(G)$. We use $V(G)$ and not $V(G) \setminus U_a$
 for technical reasons. However,  we will momentarily consider the subgraph
 of ${A}'_{aa'}$ induced between $U_a$ and $U_a'$, so that this conflict will
 be resolved.
\end{remark}
 \begin{remark}\rm In the case that $\Gamma_{aa'}= \emptyset$ we take ${A'}_{aa'}$
 to be the complete bipartite graph, as the condition in the definition
 is fulfilled vacuously.
 \end{remark}
 Next define the bipartite subgraph $A_{aa'} \subset {A'}_{aa'}$
 as the induced graph between the  vertex sets  $U_a$ and $U_{a'}$, and finally
 let  $$A=\bigcup_{x_ax_{a'}\in E(M)} A_{aa'}.$$
 See Figure 7 below for an illustration of $A$.
 Our definition
 of $A_{aa'}$ guarantees that every special constellation in $G$ yields
 a (not necessarily induced) subgraph of $A$ isomorphic to $M$.  More formally
 let
 $$\#(M,A) = \mbox{ the number of (consistent) copies of $M$ in $A$} .$$
 \begin{prop}\label{keycount}
 $$s_R \le  \#(M,A)$$
  \end{prop}
 \nin{\bf Proof:}
 This follows immediately from the definition of $A$ and the Missing Leg Property.
 Every special constellation $S(X)$ has one vertex exactly in
 each of the sets $U_a$.
 The definition of $A$ is designed so that a copy of $M$ is placed
 on $X$ if condition (a) of the Missing Leg Property holds for $S(X)$.
 Note also that every copy of $M$ corresponds to at most one special constellation,
 else the corresponding copy of $S(X)$ would violate
 condition (b) or (c) of the Missing Leg Property.
 \hfill $\Box$


\begin{figure}[htb]
\begin{center}
      \input{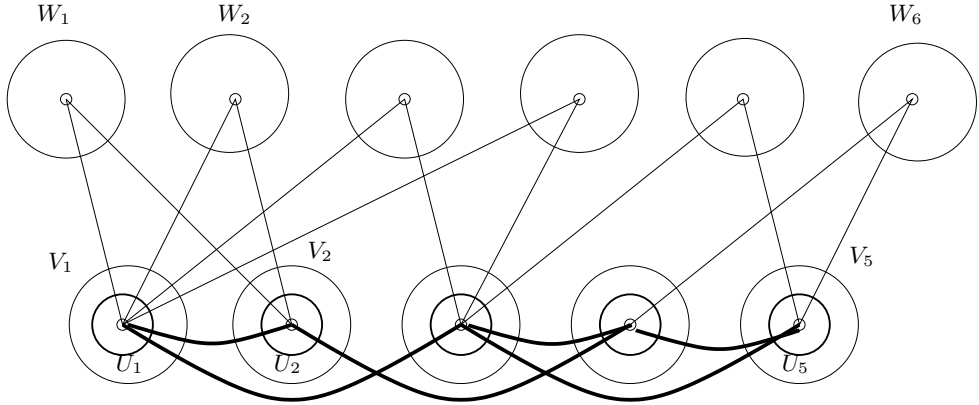}
      \caption{The auxiliary graph $A$}
\end{center}
\end{figure}

\bigskip
\begin{remark} \rm  Note that most likely there will be many ``fake'' copies of $M$ in $A$, that is, copies of $M$ which do not correspond to any special constellation  $S(X)$. There are two reasons for this overcount. First,  we may create a copy of $M$ in $A$ over a set $X\not\in\X_3$, or for which $\widehat X$ is not isomorphic to $\widehat M$, but rather contains $\widehat M$ as a proper subgraph. Second, in veryfying the definition of an edge in $A$ we may have used the vertices of $\widehat N_i(b)\setminus N_i(b)$, which, in fact, are not neighbors of $v_i$.
\end{remark}
 We will show now that
 $ A_{aa'}$ inherits from $G$ some random-like properties.
 This in turn will allow us to use Proposition
 \ref{graphcount} to easily estimate the number of copies of $M$ in $A$.
Targeting Property (P6) of $G$, for each $aa'$ with $x_ax_{a'}\in E(M)$,  define
$$\gamma_{aa'}=\prod_{b\in\Gamma_{aa'}}(1-(1-p)^{2r_{ab}}),$$
and note that
\begin{equation}\label{bounds}
(1+2C^2)^{-\phi}\prod_{b\in\Gamma_{aa'}}2pr_{ab}\le\prod_{b\in\Gamma_{aa'}}\frac{2pr_{ab}}{1+2pr_{ab}}
\le\gamma_{aa'}\le\prod_{b\in\Gamma_{aa'}}2pr_{ab},
\end{equation}
where the leftmost inequality follows from the fact that $2pr_{ab}\le p(2np)\le2C^2$. For our next lemma recall notation $x\stackrel{\epsilon}{\sim}y$ defined at the end of Section~\ref{intro}.
 \begin{lemma}
 \label{rege} For each $aa'$ such that $x_ax_{a'}\in E(M)$, the graph $ A_{aa'}$ is
 $n^{-1/30}$-regular  with density
 $$d(A_{aa'})\stackrel{n^{-1/30}}{\sim}\gamma_{aa'}.$$
 \end{lemma}
 \nin{\bf Proof:}
 Our main tools will be Lemma \ref{DLR} and Property (P6) of Definition~\ref{defG}.
   This property implies that for all $v_i,v_j\in U_a$, the degrees ${\rm deg}(v_i)$
 and co-degrees ${\rm codeg}(v_i,v_j)$ in $ A'_{aa'}$  are  quite close to their
 averages. More precisely,  let $\Gamma_{aa'}=\{b_1\dots b_g\}$. Applying Property (P6) first with
 $k=g$, and
 $S_l = \overline{N}_i(b_l)$, $l=1\dots k$, and then again with
$k=2g$, $S_l = \overline{N}_i(b_l)$ and $S_{k+l} = \overline{N}_j(b_l)$,  $l=1\dots k$, we get, respectively,

 \begin{list}{}{}
 \item[(I)]
 $\forall v \in U_a \ |{\rm deg}(v)-n\gamma_{aa'}|< n^{4/5}$, \ \ and

 \item[(II)]
 $\forall v,v' \in U_a \ |{\rm codeg}(v,v')-n\gamma^2_{aa'}|< n^{4/5}.$
 \end{list}

 Properties (I) and (II) imply, respectively, (i) and (ii) of
  Lemma~\ref{DLR} with $d=\gamma_{aa'}$,  and, say,  $\eps = n^{-11/60}$ (and with $W=[U_a]^2$).
Indeed, by (I),
$${\rm deg}(v)>n\gamma_{aa'}- n^{4/5}=n(\gamma_{aa'}- n^{-1/5})>n(d-\epsilon),$$
which is precisely condition  (i) of
  Lemma~\ref{DLR}.
Also, by Corollary \ref{lowbounds}(b), and the above lower bound (\ref{bounds}) on $\gamma_{aa'}$,
$$\gamma_{aa'}\ge(1+2C^2)^{-\phi}\prod_{b\in\Gamma_{aa'}}2pr_{ab}>\frac1{\log n}
\left(\frac{np^2}{(\log n)^{\phi+4}}\right)^g>\frac1{\log n}\frac{c^{2\phi}}{(\log n)^{\phi(\phi+4)}}>n^{-1/60},$$
 and so
 $$\gamma^2_{aa'} + n^{-1/5} < \left(\gamma_{aa'} + n^{-11/60}\right)^2.$$
 Thus, by (II),
$${\rm codeg}(v,v')<n\gamma^2_{aa'}+ n^{4/5}=n(\gamma^2_{aa'}+ n^{-1/5})<n(d+\epsilon)^2,$$
and so,
condition (ii) of
  Lemma~\ref{DLR} is true for {\it all} pairs $v,v'$.
Moreover,  Corollary \ref{lowbounds}(a) implies that for all $1\le a\le\nu$
  $$|U_a|\geq \frac n {\log^{\phi+4} n} > 2n^{11/60}=\frac2{\epsilon}.$$
Hence, by Lemma~\ref{DLR} the graph
 $A'_{aa'}$ is $(16 n^{-11/60})^{1/5}$-regular. Recall that  $A_{aa'} \subset {A'}_{aa'}$ is the induced (bipartite) subgraph of $A'_{aa'}$
 on the vertex set  $U_a \cup U_{a'}$.
 Therefore, we infer by Observation \ref{subreg} with $\eps'=|U_{a'}|/n\ge(\log n)^{-\phi-4}$,
 that $A_{aa'}$ is  $n^{-1/30}$-regular.
  \hfill $\Box$

\bigskip
 \nin{\bf Proof of Lemma~\ref{sand}}: In view of Lemma \ref{rege}, we apply Proposition \ref{graphcount} with $k=\nu$ and $B_{aa'}=A_{aa'}$ if $x_ax_{a'}\in E(M)$ and otherwise $B_{aa'}$ being the bipartite complete graph between $U_a$ and $U_{a'}$, so that
$\#(M,A)=\#(K_\nu,\bigcup_{a,a'}B_{aa'})$. Hence,
 $$\#(M,A)=(1+o(1))\prod^{\nu}_{a=1}|U_a|\prod_{aa'} \gamma_{aa'}
 \leq(1+o(1))  (2p)^\phi\prod^{\nu}_{a=1}|U_a|\prod^{\phi_a}_{b=1} r_{ab},$$
 which, by Proposition \ref{keycount}
 proves (\ref{numb}) and thus completes the proof of Lemma \ref{sand}. \hfill $\Box$

 \section{The Core Section (Proof of Lemma \ref{wrapup})}
 \label{cores}
In this section we will define the notion of a core, show that
 every hitting set of the family $\cal S$ of all special constellations contains a core, prove that there are few cores
 and that every core is large. In other words, we shall prove Lemma \ref{wrapup}. Fix an arbitrary $\tau>0$.

\subsection{What is a Core?}
Cores will be defined through a $\delta$-regular partition $\Pi_1$ of $G$ guaranteed by Lemma~\ref{reglem}. To
make this definition more user-friendly, after applying the regularity lemma we will give suggestive names to all
the parts of the final partition, and only then define  cores formally.

\subsubsection{Applying the Regularity Lemma}\label{apprl}

Let $\ga_3$ and $q$ be as defined in Section \ref{TandC}, and constants $\nu,c,C$  as in the Setup on page \pageref{setup} (see also the Glossary). For an arbitrary $\tau>0$, let $\eta$ be any constant satisfying  $0<\eta<1$ and
$$\eta \log(1/\eta) \leq \frac {\tau}{4C\nu}.$$
Furthermore, let us set
$$\alpha_4=\frac{(\nu+q)^{\nu+q}}{3C^{2q}}\ga_3,$$
$$\gd = \min \left\{ 0.4 \eta^q ,\left(\frac{\alpha_4}{8D^*}\right)^2\right\},$$
 and
  $$\eps(\ell)= \min\left\{\ell^{-4q}, \left( \frac{\delta
 \alpha_4}{D^*2^{q+3}q (\nu+q)^{\nu+q-2}}\right)^{4q} ,q^{-4},2^{-6-6q}\right\}.$$

For each star forest $S$ with at most $q$ edges there are constants $T_0=T_0(S)$, $L_0=L_0(S)$ and  $n_0=n_0(S)$ determined by Lemma \ref{reglem} with $D=2$, $D^*=2^q+1$,  the above defined $\gd$ and $\eps(\ell)$, and with $H=S$. The meaning of these constants is that  for every $n > n_0$, for all  $0<p,p^*<1$,  all $S$-uniform, $S$-graphs
$F$  with $|V(F)|=n$,  and no  $(2,p)$-dense patches, and for all ${\cal S}\subseteq\cee(F)$ with no
$(2^q+1,p^*)$-dense  $S$-patches,  there exists a refinement $\Pi_1$ of the initial
 partition ${\Pi}_0$. This refined partition $\Pi_1$ is a
 $\delta$-regular, $(\eps(l),p)$-uniform, $(t,l)$-partition for some $l \leq L_0$ and $t \leq
 T_0$. Set $T_1=\max_ST_0(S)$, $L_1=\max_SL_0(S)$ and $n_1=\max_Sn_0(S)$.

 Let $G\in\G$ be as in Lemma \ref{wrapup} (in particular, $|V(G)|=n>n_1$) and let $F$ be a subgraph of $G\in\G$ determined by the
partition $\pi_0$. Furthermore, let  $\cal S$ be the family of special constellations in $F$  generated by assumption (\ref{assumption}) of Theorem \ref{allweneed} as shown in Section \ref{TandC}.
Recall that $S$, a star forest with $\nu$ stars and $\phi\le q$ arms altogether, is the isomorphism type of the
constellations. For clarity, we will be further assuming that $n$ is divisible by $\nu+\phi$, so that the initial equipartition $\Pi_0$ has, in fact, all sets of equal size.
We will apply Lemma \ref{reglem} to the pair $(F,{\cal S})$  with $p=p(n)$ as in the Setup on page \pageref{setup} and $p^*=p^\phi$.
We have proved in Section \ref{nodense} that in such a setting $(F,{\cal S})$ does fulfill all the above assumptions of Lemma \ref{reglem}.
For the purpose of regularization, set $F_{ab}=G[V_a,W_{ab}]$ and let ${\Pi}_0$  denote the {\it initial
partition} made of $\pi_0$, the partition of vertices,  and of the edge sets between the corresponding vertex
sets.
 In other words, the elements of $\Pi_0$ are:
\begin{itemize}
 \item vertex sets: $V_a$, $a=1,\ldots, \nu$; $W_{ab}$, $b=1,\ldots, \phi_a$,
 for $a=1,\ldots , \nu-1$.
 \item edge sets: $E(F_{ab})$, $a=1,\ldots, \nu$, $b=1,\ldots, \phi_a$.
 \end{itemize}
Finally, let $$F= \bigcup_{ab}F_{ab}$$ and note that $F\in\F(S,m)$, where
$$m=\frac n{\nu+\phi}=|V_a|=|W_{ab}|.$$

 Let $\Pi_1$ be the final $(\eps(l),p)$-uniform, $\delta$-regular $(t,l)$-partition, $t\le T_1$, $l\le L_1$,
  guaranteed by  Lemma \ref{reglem}.
  It consists of vertex equipartitions
 $V_a=V_{a}^{1}\cup\cdots\cup V_{a}^{t}$, $a=1\dots \nu$, and
 $W_{ab}=W_{ab}^{1}\cup\cdots\cup W_{ab}^{t}$, $a=1\dots \nu-1$,
 $b=1\dots \phi_a$, and graph partitions
 $$F_{ab}^{i,j}=G[V_a^i, W_{ab}^j]=\bigcup^{l(ab,i,j)}_{k=1} F_{ab}^{i,j,k},\quad
 a=1\dots \nu-1, \quad b=1\dots \phi_a,$$
 $$l(ab,i,j) \leq l,\qquad i,j=1\dots t.$$
For simplicity, we assume that $m$ is divisible by $t$ and thus for all
$a,b,i$ and $j$ we have
$$|V_a^i|=|W_{ab}^j|=\frac mt.$$
  The number of edges in
 different graphs $F_{ab}^{i,j,k}$ may vary.
 Note that $F = \bigcup_{a,b,i,j,k} F_{ab}^{i,j,k}$.

\begin{remark}\label{divide}\rm
The assumption that $t$ divides $m$ is not very serious. Since $t\le T_1$, every $m$ which is divisible by $T_1!$ is also divisible by $t$, no matter what $t$ turns out to be.
Thus, our proof is absolutely formal for all $m$ which are the multiples of $T_1$, that is, for all $n$ which are the multiples of $T_1!(\nu+\phi)$. The argument for other values of $n$ is essentially the same with only some obvious, though tedious, notational complications.
\end{remark}

\subsubsection{Anatomy and Polyads}\label{ana}

We use suggestive terminology throughout this section. The
definitions of tubes and polyads given below are equivalent to
those introduced in Section~\ref{regular} in a more general
setting of $H$-graphs. Figure 8  illustrates these notions.
 \begin{itemize}

\item The sets $V_a^i$ will be called {\em palms}.
\item The sets $W_{ab}^j$ will be called {\em fingernails}.
 \item The subgraphs $F_{ab}^{ijk}$ will be called {\em fingers}.
 \item The induced subgraphs $F_{ab}^{ij}=\bigcup^{l(ab,i,j)}_{k=1} F_{ab}^{ijk}$ will be called {\em tubes}.
\item {\it A polyad} $P$ consists of $\nu$ palms $V_a^P$ (one for each $a$),  $\phi$ fingernails $W_{ab}^P$ ($\phi_a$ for each $a$), and $\phi$ fingers $F_{ab}^P$ (one from each tube $F_{ab}^{ij}$, where $V_a^P=V_a^i$ and $W_{ab}^P=W_{ab}^j$).
\item {\it An olympiad} $\calo$  consists of $\nu$ palms $V_a^\calo$ (one for each $a$),  $\phi$ fingernails $W_{ab}^\calo$ ($\phi_a$ for each $a$), and $\phi$ tubes $F_{ab}^\calo=F_{ab}^{ij}$, where $V_a^\calo=V_a^i$ and $W_{ab}^\calo=W_{ab}^j$.
\item The partition $\Pi_1$ will be called {\it myriad} as it contains all the elements of the puzzle, that is, it consists of all $t\nu$ palms, $t\phi$ fingernails, and at most $t^2l\nu\phi$ fingers.
 \end {itemize}
Note that there are precisely $t^{\nu+\phi}$ olympiads in the myriad, and at most $l^\phi$ polyads in any given
olympiad. Throughout, we will treat the olympiads (and the myriad) both as a set of polyads and as a subgraph of $F$
(which is the union of those polyads.)


\begin{figure}[hbt]
\begin{center}
     \input{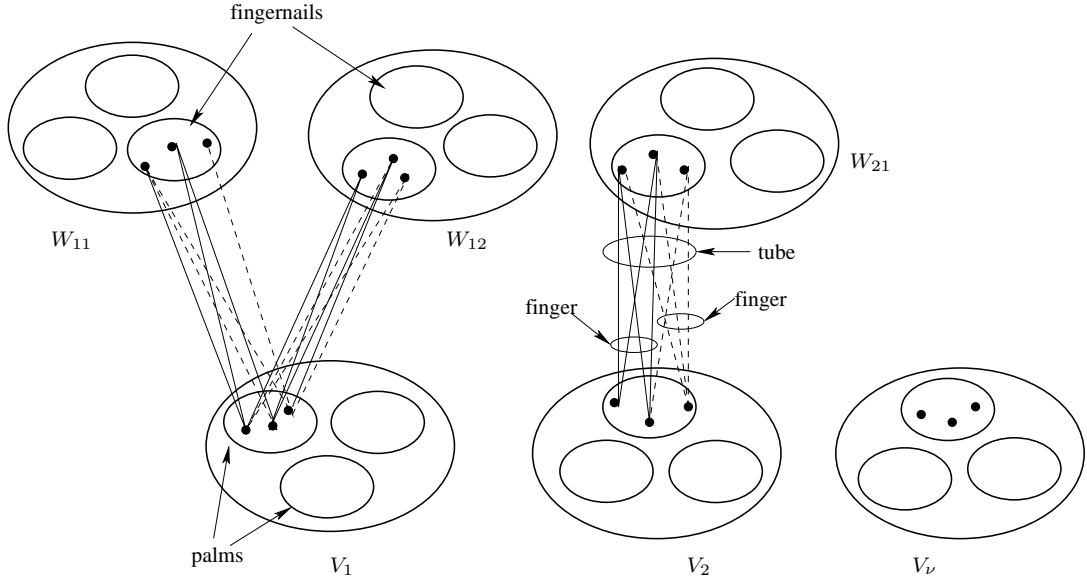}
      \caption{The Anatomy}
\end{center}
\end{figure}

\



\subsubsection{The Definition of a Core}\label{defcor}

Let  $\eps=\eps(l)$.
 We  say that a polyad $P$ is {\em nice} if all the  fingers in $P$ are $(\eps,p)$-regular, and have
 $p$-density at least $ 3\eps^{\frac{1}{2\phi}}$.
 We  say that a polyad $P$ is {\em healthy} if the following three conditions hold:
 \begin{itemize}
 \item $P$ is  nice,
 \item $P$ is $\gd$-regular,
 \item $ d_P >\delta$.
 \end{itemize}
 We define  precores first. A \emph{precore} will be a subgraph which contains a substantial
 part of every healthy polyad $P$. Roughly speaking these parts will be  formed
 by choosing at least a half of the vertices  of just one palm
 $V_a^P$, and for each vertex $v \in V_a^P$ choosing a finger
 $F^P_{ab}$, where $b=b(v)$, and taking more than a $(1-\eta)$-fraction  of the ${\rm deg}_P(v,W_{a,b}^P)={\rm deg}_{F^{P}_{ab}}(v)$
 edges incident to $v$ in that finger, where $\eta$ is defined at he beginning of Section \ref{apprl}.

 More formally,
 denote the set of all healthy polyads by
 $\calp^{healthy}.$ For each  $ P\in\calp^{healthy}$ define a family ${\rm PRECORE}^P$ of subgraphs
 of $P$ as follows: $J\in
 {\rm PRECORE}^P$ if $J\subseteq P$ and
 \begin{list}{}{}
 \item[(i)] $\exists\ 1\le a_J\le \nu:|V(J)\cap V^P_a|\cases{=0& for
 $a\ne a_J$\cr \ge \frac m{2t}=0.5|V_a^P| & for $a=a_J$}$
 \item[(ii)] $\forall v \in V(J)\cap V^P_a\  \exists\ b(v) \in \{1\dots \phi_a\}$
 such that
 $${\rm deg}_J(v,W_{a_J,b}^P)\cases{=0& for
 $b\ne b(v)$\cr > (1-\eta) {\rm deg}_P(v,W_{a_J,b}^P)& for $b=b(v)$}.$$
\end{list}
 Now we are ready to define the family {\rm PRECORE} of subgraphs of $F$:
 \beq{defprecor}
 {\rm PRECORE} = \{J: J=\bigcup_{P\in \calp^{healthy}}J^P,\quad J^P\in {\rm PRECORE}^P\}
 \enq
 The graphs in ${\rm PRECORE}$ are called {\it precores}. Thus, a precore $J$ is a union of subgraphs $J^P$, one for every healthy polyad $P$.

 We are one step from the definition of a core.

 Recall that  the (pre)coloring $\gs$ of all the edges of
$\bigcup_{X \in \X_3}S(X)\subseteq F$, that
is, of all the edges of special constellations, is determined by a fixed triangle-free coloring $\gs'$ of the graph $M$. For any precore $J$ contained in $\bigcup_{X \in \X_3}S(X)$, the coloring $\gs$
 partitions its edges
 into $J_{red}$ and $J_{blue}$.
 Define $Maj(J)$ to be the bigger of the two sets (ties are
broken arbitrarily).
  Then
 $${\rm CORE} = \{ Maj(J) : J \in {\rm PRECORE}, J\subseteq\bigcup_{X \in \X_3}S(X)\}.$$
 The graphs in $\rm CORE$ are called {\it cores}.
The reader is encouraged to draw an analogy between the above
definition and that given in Section \ref{ill}.

In the forthcoming subsections we will verify that the   family $\rm CORE$ defined above satisfies conditions (a-c) of Lemma \ref{wrapup}.

 \subsection{Every Hitting Set of $\mathcal{S}$ Contains a Core}
 We have come to the heart of our construction. The following lemma, together with
 Lemma \ref{TransInEC}, immediately implies Lemma \ref{wrapup}(a) -- see below.
Recall that a hitting set of a family of graphs is a set of edges which intersects the edge set of every graph in the family.
 \begin{lemma}\label{CoreInColoring}
 Let $T$ be a hitting set of the family $\mathcal{S}$ of special constellations. Then there exists a precore $ J \in
 {\rm PRECORE}$ such that $J \subseteq T$.
 \end{lemma}
 \nin{\bf Proof:}  Consider an arbitrary subset $T$ of $E(F)$ not containing any precore  $ J\in {\rm PRECORE}$.
We are going to show that there is a special constellation disjoint from $T$.
Because $T$ contains no precore, there exists  $ P_0\in \calp^{healthy}$ such that
 for all $J\in {\rm PRECORE}^{P_0}$ we have $J\not\subseteq T$ (otherwise, $T$ would contain a union as in (\ref{defprecor})).
 That is, for all $a=1\dots \nu$,
the set
 $$\widetilde V^{P_0}_a=\{v\in V^{P_0}_a:\forall\ b=1\dots \phi_a:
 {\rm deg}_{T\cap F^{P_0}_{ab}}(v)
 \le (1-\eta){\rm deg}_{F^{P_0}_{ab}}(v)\}$$
has size $$|\widetilde V^{P_0}_a |>  \frac {m}{2t}\,.$$
In other words, more than half of the  vertices  of each palm retain in $T$ at most $(1-\eta)$-fraction of their neighbors they had within each finger of $P_0$.
 Defining $Q_0=P_0-T$, we see that for all $a=1\dots \nu$,
 if $v\in\widetilde V^{P_0}_a$ then for all $b=1\dots \phi_a$
 $${\rm deg}_{Q_0\cap F^{P_0}_{ab}}(v)\ge
 \eta{\rm deg}_{F^{P_0}_{ab}}(v)\,.$$
Recall that for a subgraph $R$ of $F$, symbols $c_R$ and $s_R$ stand, respectively, for the numbers of constellations and
special constellations contained in $R$.
 Let $\widetilde Q_0$ and $\widetilde P_0$ be induced subgraphs of, respectively, $Q_0$ and $P_0$ with sets $V^{P_0}_a$ restricted to $\widetilde V^{P_0}_a$. Then, trivially $c_{Q_0}\ge c_{\widetilde Q_0}$ and, by formula (\ref{key}), $c_{\widetilde Q_0}\ge\eta^{\phi}c_{\widetilde P_0}$.

Since  $P_0$ is nice, all fingers $F^{P_0}_{ab}$ are
 $(\eps,p)$-regular with $p$-density at least $3\eps^{\frac{1}{2\phi}}$. Moreover, by Observation \ref{p-subreg}, all subfingers $F^{ P_0}_{ab}\cap\widetilde P_0$ are
 $(2\eps,p)$-regular with $p$-density at least $3\eps^{\frac{1}{2\phi}}-\eps>(2\eps)^{\frac{1}{2\phi}}$ and also $d_p(F^{P_0}_{ab})\stackrel{\eps}{\approx}d_p(F^{ P_0}_{ab}\cap\widetilde P_0)$.
 Hence, by two applications of Proposition \ref{stars}, to $P_0$ and to $\widetilde P_0$,
$$\frac{c_{\widetilde P_0}}{c_{P_0}}>(1-\eps^{1/4})(0.5)^\nu>(0.4)^\nu ,$$ and consequently,
$$c_{Q_0}\ge c_{\widetilde Q_0}\ge\eta^{\phi}c_{\widetilde P_0}\ge(0.4)^\nu\eta^\phi c_{P_0}>\delta c_{P_0},$$
by our definition of $\delta$.
 Therefore,
 since $P_0$ is a $\delta$-regular polyad of
 $(S,p^*)$-density greater than $\delta$,
 $$\frac{s_{Q_0}}{p^\phi c_{Q_0}}=d_{Q_0}> d_{P_0}-\delta  >0,$$
 implying that  $s_{Q_0}\ge 1$.
 This means that there exists a promised special constellation $S\in \mathcal{S}$
  such that $S\subseteq Q_0\subseteq F-T$.\hfill $\Box$

\nin{\bf Proof of Lemma 2.4(a):} Note that for every triangle-free coloring $\chi$
of $G$ the set ${\rm Agree}(\chi,\sigma)$ is a hitting set of the family
$\mathcal{S}$. Thus, by Lemma \ref{CoreInColoring}, there is a precore $J$ contained in ${\rm
Agree}(\chi,\sigma)$. Since the latter set is contained in $\bigcup_{X \in \X_3}S(X)\subseteq F$,
the set $Maj(J)$ is a monochromatic core under $\chi$.\qed

 \subsection{There Are Few Cores}
 Now we give a quick proof of Lemma \ref{wrapup}(b).
 \begin{lemma}\label{SmallCore} $|\rm CORE|<\exp\{\tau n^{3/2}\}$.
 \end{lemma}
  \nin{\bf Proof:}
  We will actually show that $|{\rm PRECORE}| <\exp\{\tau n^{3/2}\}$,
  which is stronger. It follows from the definition of  precores that
  for all $ \ J\in {\rm PRECORE} $, all
   $ v\in V(J)\cap \bigcup\limits^{\nu}_{a=1} V_a$,
   and for all $a,b,i,j$ and $k$
 $${\rm deg}_{J\cap F_{ab}^{ijk}} (v)\cases{=0, & or\cr
 \ge (1-\eta) {\rm deg}_{F_{ab}^{ijk}}(v).}$$ We will bound the number of {\it all}
 subgraphs of $F$ with that property.  For every vertex $v\in V_a$ there are at most $
 2^{\phi tl}$ choices of the ``substantial" fingers $F_{ab}^{ijk}$ along which $v$ has  positive degree.  For
 each choice of substantial fingers, if the degree of $v$ in their
 union is $r$, then the number of choices of neighbors
 of $v$ in $J$ (which is also the number of choices of the
 non-neighbors) is  at most
$$\sum_{k=0}^{\eta r}{r\choose k}\le\eta r{r\choose \eta r} < n
 {2pn\choose \eta 2pn},$$ by Property (P3) of $G$.  Thus
 $$|{\rm PRECORE}|<\left(2^{\phi tl}
 n{2pn\choose \eta 2pn}\right)^{\nu n}<\left(2^{\phi tl}n\right)^{\nu n}
\left(\frac{e}{\eta}\right)^{2C\eta\nu}<e^{\tau n^{3/2}}$$
 for $n$ sufficiently large, by the definition of $\eta$.\hfill $\Box$

 \subsection{Cores Are Large}
In this section we will show the most difficult part of Lemma \ref{wrapup}, that is, part~(c). We will prove it
with
$$\lambda =\frac1{160} \left(\frac{\alpha_4}{D^*q(\nu+q)} \right)^2,$$
where
$$\ga_4= \frac{(\nu+q)^{\nu+q}}{3C^{2q}}\ga_3\qquad\mbox{ and }\qquad D^*=2^q+1.$$
Recall that olympiads were defined in Section \ref{ana} and can be viewed as sets of polyads.
 The plan is to cover a precore $J$ by subgraphs arising from different olympiads:
$$J=\bigcup\limits_{\calo \in \Pi_1} J_\calo,\quad \mbox{ where }\quad J_\calo=\bigcup\limits_{P\in\calo^{\rm healthy}} J^P\quad\mbox{ and }\quad\calo^{\rm healthy}=\calo\cap\calp^{\rm healthy},$$
 and then to show
(see Lemma \ref{clm5.43} below) that for a certain  type of olympiads,  called perfect, we have
$$|J_\calo |\ge\frac{\alpha_4 pm^{2}}{10D^*\phi^2t^2}.$$
We will  prove that there are at least $(\alpha_4/8D^*)t^{\nu+\phi}$ perfect olympiads (see Lemma \ref{ManyPerfect} below). This and the obvious
fact that an edge may belong to at most $t^{\nu+\phi-2}$ olympiads, yields the following result.
\begin{lemma}\label{LargeCore}
 For all $J\in {\rm PRECORE}$ we have  $|J| \ge 2\lam pn^2$.
 \end{lemma}
\nin{\bf Proof:} By Lemmas \ref{ManyPerfect}  and \ref{clm5.43} we have
$$|J|\ge \frac{\alpha_4}{8D^*} t^{\nu+\phi}\times\frac{\alpha_4pm^2}{10D^*\phi^2t^2}\times  t^{-\nu-\phi+2} = 2\lambda pn^2.$$ \qed
Lemma \ref{LargeCore}, together with the fact that $|E(G)|<pn^2$ (see Property (P3) of
$\G$), imply Lemma \ref{wrapup}(c).
\subsubsection{Perfect Olympiads}

 Recall that $\kappa=\sum_{P\in\calp_{\Pi_1}} c_P $ is the  total volume of the myriad.
 By Proposition \ref{stars} (see Example \ref{cee})
$$\kappa\sim m^{\nu+\phi}p^\phi.$$
 Let $$\kappa_0=\frac \kappa{t^{\nu+\phi}}$$ be the
 average volume of an olympiad.
 Setting $$c_{\calo}=\sum\limits_{P\in\calo} c_P,$$
 we have, again by Proposition \ref{stars}, that for all olympiads $\calo$
 $$c_{\calo}\sim\kappa_0.$$
 Let
 $\calo^{\rm reg}$ denote the set of $\delta$-regular polyads contained in
 $\calo$, and $\calo^{\rm irr}=\calo\setminus\calo^{\rm reg}$ be the set of irregular ones. Write
 $$c_{\calo}=c_{\calo}^{\rm reg} + c_{\calo}^{\rm irr}$$ where
 $c^{\rm reg}_\calo= \sum_{P \in \calo^{\rm reg}}
 c_P$ and $c_{\calo}^{\rm irr}=\sum_{P \in \calo^{\rm irr}} c_P$, and
 set $$s_\calo=\sum\limits_{P\in\calo} s_P.$$
 Finally, let $\calo^{\rm nice}$ denote the set of nice polyads contained in
 $\calo$
 and $$c_{\calo}^{\rm
 nice}=\sum_{P \in \calo^{\rm nice}} c_P.$$
\begin{defn}\rm We call an olympiad $\calo$ {\it perfect} if the following three conditions hold:
 \begin{equation}\label{bronze}
 c^{\rm nice}_\calo \geq (1-\delta)\kappa_0
\end{equation}
\begin{equation}\label{silver}
  c^{\rm reg}_\calo \geq (1-\sqrt\delta)\kappa_0,
\end{equation}
and
 \begin{equation}\label{gold} s_\calo\ge \alpha_4 p^\phi\kappa_0
 \end{equation}
\end{defn}
\begin{lemma}
 \label{ManyPerfect} At least $\displaystyle\frac{\alpha_4}{8D^*} t^{\nu+\phi}$ olympiads are
 perfect.
 \end{lemma}
 \nin{\bf Proof:} By Property (P5) of $\G$,
$$|E(F)|\sim\phi m^2p.$$
Let us begin with a count of the edges of $F$ belonging to fingers which are either $(\eps,p)$-irregular or of
$p$-density smaller than $3\eps^{{1}/{2\phi}}$. Let us call such fingers {\it nasty}. First, from the
 $(\eps,p)$-uniformity of partition $\Pi_1$, there are at most $\eps |E(F)|$ edges belonging to fingers that are $(\eps,p)$-irregular.

Next note that  every tube is an induced subgraph of $F$ split into at most $l$ fingers,  so even if they were all
of $p$-density smaller than $3\eps^{{1}/{2\phi}}$, they
 would contribute no more than a $3l \eps^{{1}/{2\phi}}$-fraction of the edges. As we have $l \le
 \eps ^{-{1}/{4\phi}}$ and  $\eps\le \eps^{{1}/{4\phi}}$, altogether there are fewer than $4\eps^{{1}/{4\phi}}|E(F)|$
edges which belong to nasty fingers of $F$.

Furthermore, by Property (P3) of $\G$, any given edge may belong to at most $n^{\nu-1}(2np)^{\phi-1}$
constellations. Call a constellation {\it spoiled} if at least one of its edges belongs to a nasty finger. Then, recalling the definition of a nice polyad from Section \ref{defcor},
$c_\calo-c_\calo^{\rm nice}$ counts precisely  those constellations of olympiad $\calo$ which are  spoiled. Thus,
$$\sum_{\calo}(c_\calo-c_\calo^{\rm nice})\le(1+o(1))4\eps^{\frac{1}{4\phi}}\phi m^2pn^{\nu-1}(2np)^{\phi-1}.$$
On the other hand,
denoting by $\Upsilon^-_{1}$ the number of olympiads which {\it do not} satisfy condition (\ref{bronze}), we have
$$\Upsilon^-_{1}\delta\kappa_0\le\sum_{\calo}(\kappa_0-c_\calo^{\rm nice}),$$
which, using the definitions of $\eps$, $\kappa_0,\kappa$, and $m$, and the relation  $c_{\calo}\sim\kappa_0$, yields the bound
$\Upsilon^-_{1}\le(1+o(1))(\alpha_4/4D^*)t^{\nu+\phi}$.

 Next, let $\Upsilon^-_{2}$ be the number of olympiads which {\it do not} satisfy condition (\ref{silver}). We claim that $\Upsilon^-_{2}\le\sqrt\delta t^{\nu+\phi}\le(\alpha_4/8D^*) t^{\nu+\phi}$, where the latter inequality follows just by the definition of $\delta$.
 Indeed, if $\Upsilon^-_{2}>\sqrt\delta t^{\nu+\phi}$
  then
 $$\sum_{P\in\calp_{\Pi_1}^{\rm irr}}c_P=\sum\limits_{\calo}\sum\limits_{P\in\calo^{\rm irr}} c_P>\sqrt\delta
 t^{\nu+\phi}\sqrt\delta \kappa_0 =\delta \kappa$$ -- contradiction with the
 $\delta$-regularity of $\Pi_1$ (compare Definition \ref{dfn4.4}).

Finally, let $\Upsilon^+_{3}$ be the number of olympiads fulfilling condition (\ref{gold}).
We claim that $\Upsilon^+_{3}\ge\frac23(\alpha_4/D^*) t^{\nu+\phi}$. Indeed, by  Proposition~\ref{(b)}
and the fact that $F$ has no $(D^*,p^\phi)$-dense $S$-patches, we have, for every olympiad $ \calo$,
 \begin{eqnarray*}
  s_\calo&=&\sum_{c_P<\kappa/{\log^2n}}s_P+\sum_{c_P\ge\kappa/{\log^2n}}
 s_P\le
 D^*p^\phi\kappa/{\log^2n} + D^*p^\phi\sum_{P\in\calo} c_P\\
 &=& o(p^\phi\kappa_0)+D^*p^\phi\kappa_0<\frac32D^*p^\phi\kappa_0.\hspace{2in}(*)\end{eqnarray*}
 Recall from Section \ref{TandC} that $|{\mathcal{S}}|=\alpha_3n^{\nu}$, and suppose that $\Upsilon^+_{3}<\frac23(\alpha_4/D^*)t^{\nu+\phi}$.
  Then
 $$|{\mathcal{S}}| <\frac23\frac{\alpha_4t^{\nu+\phi}}{D^*}\, \frac32D^*p^\phi\kappa_0 + t^{\nu+\phi}
 \alpha_4\, p^\phi \kappa_0=2\alpha_4p^\phi\kappa<\alpha_3n^{\nu}$$ -- a
 contradiction.
Hence, we conclude that there are at least
$$\Upsilon^+_{3}-\Upsilon^-_{1}-\Upsilon^-_{2}\ge\left(\frac23\frac{\alpha_4}{D^*}-
\frac{\alpha_4}{4D^*}-\frac{\alpha_4}{8D^*}-o(1)\right)t^{\nu+\phi}>\frac{\alpha_4}{8D^*}t^{\nu+\phi}$$
perfect olympiads.\hfill $\Box$
\bigskip

 Recall that  $\calo^{\rm healthy} =
 \calo \cap \calp^{healthy}$ is the set of healthy polyads contained in an olympiad $\calo$. The next lemma will be used in the proof of Lemma \ref{clm5.43}.
 \begin{lemma}\label{clm5.41} For a perfect olympiad $\calo$
 $$\sum\limits_{P\in \calo^{\rm healthy}}
 c_P\ge\frac{\alpha_4}{2D^*}\, \kappa_0.$$
 \end{lemma}
 \nin{\bf Proof:} Because there are no $(D^*,p^\phi)$-dense $S$-patches in $F$,
we have $c_P\ge s_P/(D^*p^\phi)$ for each polyad $P$ with $c_P\ge\kappa/\log^2n$. Hence
$$\sum_{P \in \calo^{healthy}}c_P\ge\frac1{D^*p^\phi}
 \sum_{\stackrel{P\in\calo^{healthy}}{c_P\ge  \kappa/{\log^2n}}}
 s_P.$$
Clearly,
$$ \sum_{\stackrel{P\in\calo^{healthy}}{c_P\ge  \kappa/{\log^2n}}}
 s_P\ge s_\calo- \sum_{\stackrel{P\not\in\calo^{healthy}}{c_P\ge  \kappa/{\log^2n}}}s_P-
\sum_{c_P<  \kappa/{\log^2n}} s_P.$$
Recalling the definition of a healthy polyad, to estimate further we must then
 subtract from $s_{\calo}$ the contribution coming from polyads  with $c_P\ge\kappa/\log^2n$ that are either not nice,
 or $\delta$-irregular, or of $(S,p^*)$-density less than $\delta$, and also from polyads with $c_P<\kappa/\log^2n$. First, since $\calo$ satisfies
condition (\ref{bronze}),
 $$\sum_{\stackrel{P\not \in\calo^{\rm nice}}{c_P\ge \kappa/{\log^2n}}} s_P\le D^*p^\phi\sum_{P\not \in\calo^{\rm nice}}c_P= D^*p^\phi(c_\calo-c_\calo^{\rm nice})\le(1+o(1))\delta
 D^*p^\phi\kappa_0.$$
 Second, using the
 fact that $\calo$ obeys condition (\ref{silver}),
 $$\sum_{\stackrel{P\in\calo^{\rm irr}}{c_P\ge  \kappa/{\log^2n}}} s_P<D^*p^\phi\sum_{P \in\calo^{\rm irr}}c_P=D^*p^\phi(c_\calo-c_\calo^{\rm reg})\le(1+o(1))\sqrt\delta D^*p^\phi\kappa_0.$$   Third,
 $$\sum_{P:d_P<\delta} s_P<\delta p^\phi\sum_{P:d_P<\delta}c_P\le\delta p^\phi c_\calo=(1+o(1))\delta p^\phi\kappa_0.$$  Last, by Proposition \ref{(b)}
$$\sum_{c_P<  \kappa/{\log^2n}} s_P \le(1+o(1))D^*p^\phi\frac\kappa{\log^2n}= o(1)p^\phi\kappa_0.$$
Putting all this together
 and using the fact that $\calo$ fulfills (gold) and thus $s_\calo\ge \alpha_4p^\phi\kappa_0$, we finally get
 \begin{eqnarray*}
 \sum_{\stackrel{P\in\calo^{healthy}}{c_P\ge  \kappa/{\log^2n}}}
 s_P&\ge& \left(\alpha_4-\delta D^*-\sqrt\delta D^*-\delta -o(1) \right)
 p^\phi\kappa_0\\
 &>&\left(\alpha_4-3\sqrt\delta D^*-o(1)\right) p^\phi\kappa_0
>\left(\alpha_4-\frac38\alpha_4-o(1)\right) p^\phi\kappa_0\\
&>&\frac12\alpha_4p^\phi\kappa_0
\end{eqnarray*}
 and consequently
$$\sum_{P \in \calo^{healthy}}c_P\ge\frac{\alpha_4}{2D^*}\, \kappa_0.\quad$$\hfill$\Box$

\subsubsection{Precores Are Large}

 We now  concentrate on the contribution to the precore $J$
 coming from healthy polyads of one perfect olympiad.
 Recall our notation
 $$J_\calo=\bigcup\limits_{P\in\calo^{\rm healthy}} J^P.$$
The following lemma, together with Lemma \ref{ManyPerfect} completes the proof of Lemma \ref{LargeCore}, which, in
turn, finishes off the proof of Lemma \ref{wrapup}, yielding Theorem \ref{allweneed}, and thus our main result --
Theorem \ref{main}.

 \begin{lemma}\label{clm5.43}
 For every precore $J\in {\rm PRECORE}$ and for every perfect olympiad $\calo$, we have
 $$|J_\calo|\ge\frac{\alpha_4pm^2}{10D^*\phi^2t^2}.$$
 \end{lemma}

 \noindent
 {\bf Proof}.
The idea of the proof of Lemma~\ref{clm5.43} is to take only one finger from each polyad $P$ -- the one chosen by
the majority of the vertices in the definition of $J^P$ and then show, using Lemma \ref{clm5.41} and Proposition
\ref{stars}, that many of those fingers will belong to the same tube. As such, they will be edge-disjoint and
adding up their contributions to the respective $J^P$'s will give the desired lower bound on $|J_\calo|$.

Recall the definition of the family ${\rm PRECORE}^P$ from Section \ref{defcor}.
For every $P\in\calo^{\rm healthy}$ and for every $J^P\in{\rm PRECORE}^P$, let $a_P$ and $b_P$ be such
that
$$|V(J^P)\cap V_{a_P}^\calo|\ge \frac m{2t}\mbox{ and } |\{v\in V(J^P)\cap V_{a_P}^\calo: b(v)=b_P\}|\ge\frac m{2\phi t}.$$
Then $a_P$ and $b_P$ determine a finger of $P$, denoted further by $F^P=F_{a_Pb_P}^P$, belonging to the tube
$F_{a_Pb_P}^\calo$. Obviously,
$$\left|\bigcup_{P\in\calo^{\rm healthy}}J^P\right|\ge \left|\bigcup_{P\in\calo^{\rm healthy}}(J^P\cap F^P)\right|.$$
Note that by the definition of $J^P$, by the $(\eps,p)$-regularity of fingers and the inequality
$\eps\le1/(2\phi)$,
\begin{equation}\label{fin}
|J^P\cap F^P|\ge\frac1{2\phi}(1-\eta)(1-\eps)|F^P|>\frac{|F^P|}{3\phi},
\end{equation} where the last inequality follows because, say, both $\eta<1/6$ and $\eps<1/6$.
For every $a,b$ consider the edge-disjoint union $\overline F_{ab}^\calo$ of selected fingers $F^P$ belonging to
the tube $F_{ab}^\calo$, that is
$$\overline F_{ab}^\calo=\bigcup_{a_P=a,b_P=b}F^P\subseteq F_{ab}^\calo.$$
We claim that for some $a,b$
\begin{equation}\label{ab}
|\overline F_{ab}^\calo|\ge\frac{\alpha_4}{3D^*\phi}p\left(\frac mt\right)^2.
\end{equation}
Suppose to the contrary that the opposite inequality is true for all $a,b$, and consider the subgraph $R$ of $F$
which is a union of the complements $F^\calo_{ab}-\overline F^\calo_{ab}$ of $\overline F_{ab}^\calo$ within the
respective tubes. Then, since every healthy polyad of $\calo$ must intersect  some $\overline F_{ab}^\calo$, the
subgraph $R$ is a union of unhealthy polyads, and by Lemma \ref{clm5.41}
$$c_R\le c_\calo-\frac{\alpha_4}{2D^*}\kappa_0=\left(1+o(1)-\frac{\alpha_4}{2D^*}\right)\kappa_0.$$
On the other hand, by property (P5) of Definition \ref{defG}, tubes  are $(\eps,p)$-regular of $p$-density
$1-o(1)$, and the graphs $F^\calo_{ab}-\overline F^\calo_{ab}$ are complements of unions of at most $l$ $(\eps,p)$-regular
fingers within the tube $F^\calo_{ab}$. Thus, they are themselves $(\eps',p)$-regular, where $\eps'\le l\eps\le\eps^{1-1/4\phi}$, with $p$-densities at least $1-\alpha_4/(3D^*\phi)+o(1)$.
Thus, by Proposition \ref{stars} and by (\ref{ab}),
$$c_R\stackrel{\eps}{\approx}\left(\frac mt\right)^{\nu+\phi}p^\phi\prod_{ab} d_p(F^\calo_{ab}-\overline F^\calo_{ab})\ge\left(1+o(1)-\frac{\alpha_4}{3D^*\phi}\right)^\phi\kappa_0,$$
contradicting the above upper bound.

Hence, let $a_0,b_0$ be such that (\ref{ab}) holds. Then
$$\left|\bigcup_{P\in\calo^{\rm healthy}}(J^P\cap F^P)\right|\ge\left|\bigcup_{a_P=a_0,b_P=b_0}(J^P\cap F^P)\right|=\sum_{a_P=a_0,b_P=b_0}|J^P\cap F^P|.$$
Using (\ref{fin}) and (\ref{ab}) we  easily estimate  further to obtain the desired bound:
$$\sum_{a_P=a_0,b_P=b_0}|J^P\cap F^P|\ge\sum\frac{|F^P|}{3\phi}=\frac{|\overline F_{a_0b_0}^\calo|}{3\phi}\ge\frac{\alpha_4pm^2}{10D^*\phi^2t^2}.$$
\hfill $\Box$
%
%
%

 \section{Random Graphs}
 \label{randomgraphs}
  In this section we state and prove some
 assertions about random graphs which constitute the family $\G$ of graphs playing an important role in our proof.

  Let $G=(V,E)$ be a graph, and let $v,u \in V$,
  $U,W \subseteq V$, where $U\cap W=\emptyset$, and $F\subset E$.
Recall our notation
 $${\rm deg}_G(v,W)={\rm deg}(v,W) = | \{w \in W : vw \in E\} |,$$
$${\rm deg}_G(v)={\rm deg}(v)={\rm deg}(v,V),$$
 $${\rm codeg}(v,u) = | \{ w : vw, uw \in E \} |,$$
and $e_G(U)$ and $e_G(U,W)$ for the number of edges of $G$ with, resp., both endpoints in $U$, and one endpoint in
$U$ and the other in $W$. Furthermore,  ${\rm Base}(F)$\label{base1} is the set of edges in the complete graph on
$V=V(G)$
 that form a triangle with two edges in $F$, formally:
 $$ {\rm Base}(F) = \{uv :  \ uw, wv \in F \mbox{ for some }w\in V\} .$$

 \subsection{The Graph Family $\G$}
 Given a constant $\alpha$, a sequence $c/\sqrt{n} \leq p=p(n) \leq C/\sqrt{n}$ and an integer $\nu\ge5$, let $q$
 and $\lambda$ be as in  The Glossary and let $a=a(\lam,c)$ be a constant determined by Lemma \ref{pro1}.
We now define the property $\G$ and show that it holds  for
 $G(n,p)$ asymptotically almost surely (a.a.s. for short), that is, with probability tending to 1 as $n\to\infty$.
 \begin{defn}\rm
 \label{defG} A graph $G=(V,E)$ on $n$ vertices is said to have
 property $\G=\G(p,\nu,q,\lambda,a)$ if the following hold:
 \begin{enumerate}
 \item[(P1)] The number of sets of $\nu$ vertices that span an edge
 (i.e. are not  independent sets) is $o(n^{\nu})$.
 \item[(P2)]  The number of sets of $\nu$ vertices such that some three
 of them share a common neighbor is $o(n^{\nu})$.
 \item[(P3)] For all $v \in V$,
$$(1-n^{-1/5})np \leq {\rm deg}(v)\le(1+n^{-1/5})np \leq
 2np.$$
 \item[(P4)] For all $v, u \in V$  ${\rm codeg}(v,u) \leq 3\log n.$
 \item[(P5)] For all pairs of disjoint sets $U , W \subset V$
 with $|U| , |W| \geq n/\log n$,
 $$(1-n^{-1/5})|U||W|p\le e_G(U,W) \leq(1+n^{-1/5})   |U||W|p\le 2 |U||W|p,$$
and
$$e_G(U)<|U|^2p.$$
 \item[(P6)]
  For all $1\le k\le 2q$, where $q$ is as in the Glossary, and all choices of subsets  $S_1 \dots S_{k}$
 of $V$, such that $|S_i| = s_i \leq 2np$, for all $i=1\dots k$, and
$|S_i \cap S_j| \leq 12 \log n$ for all $1\le i<j\le k$, the number
\[Z=Z(S_1,\ldots, S_{k})\] of vertices with at least one neighbor in each $S_i$ satisfies
 $$|Z-\mu| \leq n^{4/5},$$
 where
  $\mu =n\prod_{i=1}^k \mu_i$ and $\mu_i = 1-(1-p)^{s_i}$.

 \item[(P7)] $G$ has property $\T(\lambda,a)$, that is for any subgraph $F$ of $G$ with at least $\lambda |E|$ edges,
 the set ${\rm Base}(F)$ contains at least $a |V|^3$ triangles.
 \end{enumerate}
 \end{defn}

%
%

 \begin{lemma}\label{aas} For all sequences $p=p(n)$ such that $c/\sqrt{n} \leq p \leq C/\sqrt{n}$
 $$\lim_{n\to\infty}{\rm Pr}(G(n,p) \in \G)=1.$$
 \end{lemma}
 \nin{\bf Proof: } We will prove that each of the above properties holds
 a.a.s., so, of course, their conjunction also holds a.a.s.
The proofs of (P1,P2) rely on Markov's inequality, whereas those of (P3-P6) -- on various versions of the Chernoff
bound for which we refer the reader to Section 2.1 of \cite{JLR}. Recall that the set of neighbors of $v$ in $G$ is
denoted by $N_G(v)=N(v)$, while $N_G(W)$  stands for the set of vertices
 {\it outside}  $W$, each having at least one neighbor in $W$. Setting $G=G(n,p)$, note that $|N_{G}(W)|$ is binomially distributed with expectation $(n-|W|)(1-(1-p)^{|W|})$.

\medskip

 \nin{\bf (P1) }
 The expected number of such sets is $O(n^\nu p)$. Thus, by Markov's inequality
 there are a.a.s. no more than, say, $n^\nu\sqrt p=o(n^\nu)$ of them.
 \medskip

  \nin{\bf  (P2) }
The expected number of such sets is $O(n^{\nu+1} p^3)$. Thus, by Markov's inequality there are a.a.s. no more
than, say, $n^{\nu+1} p^{5/2}=o(n^\nu)$ of them.

\medskip

  \nin{\bf (P3) } For any $v \in V$,  the degree ${\rm deg}(v)$ is a binomial
 random variable with expectation $(n-1)p$. Since $c \sqrt{n} \leq
 np \leq C \sqrt{n}$, the usual Chernoff bound ((2.9) in \cite{JLR}) gives
 that
 $${\rm Pr}\left({\rm deg}(v) \not \in [(1-n^{-1/5})np, (1+n^{-1/5})np ]\right) \leq \exp\left\{-\Theta(n^{1/10})\right\}.$$
 A simple union bound shows that a.a.s. $(1-n^{-1/5})np\le {\rm deg}(v)\le (1+n^{-1/5})np$  for all
 vertices $v$ simultaneously.
\medskip

  \nin{\bf (P4)} As above, for fixed $u,v$  the random variable in
 question, ${\rm codeg}(v,u)$,  is binomial with
 expectation $(n-2)p^2 = \Theta(1)$. From  a handy version of Chernoff's bound ((2.11) in \cite{JLR})  it
 follows that
 $${\rm Pr}\left( {\rm codeg} (u,v) \geq 3\log n \right) \leq e^{-3\log n} =n^{-3}.$$
 Hence, a.a.s. the required condition holds for all $\C{n}{2}$ pairs $u,v$
 simultaneously.
\medskip

  \nin{\bf (P5)}
 Again, given two disjoint sets of vertices $U , W \subset V$
 with $|U| , |W| \geq n/\log n$, the number $e_G(U,W)$ is a binomial random variable with expectation $|U||W|p=\Omega(n^{3/2}/\log^2n)$. And, again, by the usual Chernoff bound,
$${\rm Pr}\left( d_p(U,W)\not \in [(1-n^{-1/5}), (1+n^{-1/5}) ]\right) \leq \exp\left\{-\Theta(n^{12/11})\right\}.$$
Since there are less than $2^{2n}$ choices of $U,W$, the required condition holds  a.a.s. for all
 such pairs simultaneously. The second statement is proved similarly.

\medskip

 \nin {\bf(P6)}
Split $Z=Z'+Z''$, where
 $Z'$ is the number of vertices $v$ counted by $Z$ such that $v\not\in\bigcup_{i=1}^k S_i$ and
  $ {\rm deg}(v, S_i\cap S_j)=0$ for all\ $i<j$.
 Set $$ \mu' =\Ee Z', \  \mu'' = \Ee Z''\mbox{ and } \bar{\mu}= \Ee Z = \mu' + \mu''.$$

First, by the crude Chernoff bound ((2.12) in \cite{JLR}) we have
$$Pr(|Z'-\mu'|\ge n^{4/5}/2)\le 2e^{-2n^{8/5}/(n-s)}\le
 2e^{-2n^{3/5}},$$
where $s=\sum s_i$. Observe also that there are at most
 \[\left[2np{n\choose 2np} \right]^{2q} \le
 \left[ 4C\sqrt{n} {n\choose 4C\sqrt{n}}\right]^{2q}
 =\exp\left\{O(\sqrt{n}\log n)\right\} \]
 choices of $S_1,\ldots, S_k$. Hence
 $$|Z'-\mu'|\le n^{4/5}/2$$
 holds a.a.s for all such choices.

 Second, by (P3) we have $\mu''=\Ee Z''=O((\log n)\sqrt n)$. Let $U_{ij}= S_i \cap S_j$,
 for $i,j=1, \ldots , k$.  Observe that if
 $Z'' > n^{2/3}$ then for some $i\neq j$, we have
$$\left|\bigcup_{v \in U_{ij}} N_G(v)\right| > (n^{2/3}-s)/k^2 >n^{5/9}.$$
  Consequently, by the handy Chernoff bound ((2.11) in \cite{JLR})
 \begin{eqnarray*}
 &&P(\exists\ S_1,\ldots, S_k : Z''>n^{2/3})\\
 &&\qquad\le P\left(\exists U\subset V, \ |U|\le 12 \log
 n : |N_G(U)|>n^{5/9}\right) \le n^{12\log n} e^{-n^{5/9}}=o(1).
 \end{eqnarray*}
Therefore, a.a.s. for all choices of $S_1\dots S_k$, $1\le k\le 2q$, we have
$$|Z-\bar\mu|\le |Z'-\mu'|+Z''+\mu''\le n^{4/5}/2+n^{5/9}+O((\log n)\sqrt n)<2n^{4/5}/3.$$
Finally note that $|\mu-\mu'|\le s=O(\sqrt n)$.

   \medskip
  \nin{\bf  (P7)} This is Lemma \ref{pro1}, proved in the next section.

\medskip
\nin We have thus completed the proof that $G(n,p) \in \G$ a.a.s. \hfill $\Box$

\subsection{Proof of Lemma \ref{pro1}}
Here we prove Lemma \ref{pro1}, which turned out to be both difficult and interesting.
 Work on it led to a separate paper \cite{FKRRT} where its extention is utilized in a wider context of dynamic Ramsey-type colorings
 of random graphs. The proof provided in \cite{FKRRT} is based on a counting technique for sparse random graphs. In this paper we prefer our original approach which, besides a regularity lemma, uses an upper tail estimate
established in \cite{RR3}. We state it here in a general setting of random subsets.

Given a finite set $\Gamma$ and $0<p<1$, we denote by $\Gamma_p$ a random binomial subset of $\Gamma$. Note that
$([n]^2)_p$ is nothing else but a random graph $G(n,p)$. Let ${\cal H}\subseteq[\Gamma]^h$, where $h$ is a
positive integer, and set $\mu=\Ee|\{H\in{\cal H}: H\subset\Gamma_p\}|$.

\begin{lemma}\label{del}
For all integers $\beta>0$, with probability at least $1-\exp\left\{-\frac{\beta}{2h}\right\}$ there is
$E_0\subset \Gamma_p$ with $|E_0|=\beta$, such that $\Gamma_p\setminus E_0$ contains fewer than $2\mu$ sets
from~${\cal H}$.
\end{lemma}

For more about upper tail estimates see  Section 2.6 in \cite{JLR}, \cite{JR1} and \cite{JOR}, and for further development of the above
deletion technique see \cite{JR2}.

\begin{cor}\label{cor}
A.a.s. for each pair $U,W\subset[n]$, $U\cap W=\emptyset$, there exists $E_0\subset G(n,p)$ with $|E_0|=n\log n$,
such that the bipartite subgraph of $G(n,p)\setminus E_0$ spanned between the sets $U$ and $W$ contains no more
than
$$2{|U|\choose2}{|W|\choose2}p^4$$
copies of the 4-cycle $C_4$.
\end{cor}

{\bf Proof:} Apply Lemma \ref{del} to $\Gamma=U\times W$ and $\cal H$ -- the family of the edge sets of all
4-cycles in $\Gamma$, with $s=4$, $\beta=n\log n$ and $\mu={|U|\choose2}{|W|\choose2}p^4$. \qed

\nin{\bf Proof of Lemma \ref{pro1}}

Given $c$ and $\lambda$, let $d=\min\{\lambda^3c^2/270,\lambda^4/729\}$,  and let $\varrho=\varrho(d)$ and
$c_0=c_0(d)$  be determined via Proposition \ref{L2.3} (applied with $k=3$). Further, let
$\epsilon=\min\{\varrho,\lambda/100\}$, and let $T_0=T_0(\epsilon)$  be determined by Lemma \ref{SSRL} with $D=2$,
$t_0=3$ (say), $r=1$, and  the above $\epsilon$ as inputs. We promise to prove Lemma \ref{pro1} with
$a=c_0/T_0^3$.

Let $G$ be a graph from the space $G(n,p)$ which satisfies  Properties (P3) and (P5) and the property stated in
Corollary \ref{cor}. Note that Property (P5) implies that $G$, as well as each of its subgraphs has no
$(2,p)$-dense patches. Fix $F\subseteq G$ with $|F|>\lambda |G|>(\lambda/3)pn^2$ and apply Lemma \ref{SSRL} to $F$
with $D$, $t_0$, $r$, and  $\epsilon$ as above, obtaining a partition $V=W_1\cup\cdots\cup W_t$ with $3\le
t\le T_0$ such that all but $\ep n^2p$ edges of $F$ belong to the $(\ep,p)$-regular pairs $\{W_i,W_j\}$. Let $F'$
be the subgraph of $F$ consisting of those edges, and let $x{t\choose 2}$ be the number of pairs $(i,j)$ such
that $e_{F'}(W_i,W_j)\ge\frac13\lambda p(n/t)^2$. For clarity of presentation, let us assume that $t$ divides $n$, and  so
$|W_i|=n/t$ for each $i=1\dots t$ (cf. Remark \ref{divide}). We will now estimate $x$ from below.

By  (P5) and the definition of $\ep$,

$$(\lambda/3)pn^2\le|F'|\le2x{t\choose 2}(n/t)^2p+\left({t\choose 2}-x{t\choose 2}\right)(\lambda p/3)(n/t)^2$$
and so
$$\lambda\le x+(1-x)\lambda/4\ .$$
Solving the last inequality for $x$, we obtain $x\ge\lambda/6$. Hence,  at least
$$(\lambda/6){t\choose 2}>(\lambda/10)t^2$$
 pairs $(W_i,W_j)$ are $(\epsilon,p)$-regular, each satisfying
$$e_F(W_i,W_j)\ge(\lambda p/3)(n/t)^2\ .$$
 By simple averaging, there must be an index $i_0$ and a set $J\subseteq[t]\setminus\{i_0\}$ of cardinality $|J|=\lambda t/10$ such that $(W_{i_0},W_j)$ is as above.
Set $W=W_{i_0}$ and $U=\bigcup_{j\in J}W_j$, and note that $|W|=n/t$ and
$$|U|=|J|(n/t)=\lambda n/10\ .$$

Let $B$ be the bipartite subgraph of $F$ induced by all edges with one endpoint in $U$ and the other in $W$,
shortly $B=F[U,W]$. Then ${\rm Base}(B)[W]\subseteq {\rm Base}(F)$. We will show that ${\rm Base}(B)[W]$ is $(\varrho,d)$-dense, and
hence, by Lemma \ref{L2.3}, for $n$ large enough, contains at least $c_0(n/t)^3\ge (c_0/T_0^3)n^3$ triangles,
which will complete the proof.

To this end, assume for clarity that $\varrho n/t$ is an integer, and pick any $W'\subset W$ with
$|W'|=\varrho|W|=\varrho n/t$. Our ultimate goal is to show that
$$N=|{\rm Base}(B)[W']|\ge d{|W'|\choose2}.$$

Since for each $j\in J$ the pair $(W_j,W)$ is $(\epsilon,p)$-regular with density at least $$\pi=\lambda p/3$$ and
$|W'|\ge\epsilon|W|$, by Proposition \ref{degrees}, all but at most $\epsilon|W_j|$ vertices of $W_j$ have each at
least $(\lambda/3-\ep)p|W'|)=(1-\epsilon')\pi|W'|$ neighbors in $W'$, where $\epsilon'=3\epsilon/\lambda>\ep$.
Consequently, all but at most $\epsilon|U|$ vertices of $U$ have each at least $(1-\epsilon')\pi|W'|$ neighbors in
$W'$.

Let $E_0\subset G[U,W']$ be as in Corollary \ref{cor}, i.e.  $|E_0|=n\log n$ and there are at most
$2{|U|\choose2}{|W'|\choose2}p^4$ copies of $C_4$ in $G_0=G[U,W']\setminus E_0$. Clearly, the same is true for the
subgraph $B_0=F[U,W']\setminus E_0$ of $G_0$. As, say, only  at most $n^{2/3}$ vertices of $U$ can be each
incident to more than $n^{1/3}\log n$ edges of $E_0$, still, say, all but at most $2\epsilon|U|$ vertices of $U$
have each at least $(1-2\epsilon')\pi|W'|$ neighbors in $W'$. Put in another way, for all but at most
$2\epsilon|U|$ vertices  $u\in U$
$$d_u={\rm deg}_{B_0}(u)>(1-2\epsilon')\pi|W'| .$$

Let $[W']^2=\{e_1,\dots,e_{{|W'|\choose2}}\}$. Denote by $x_i$ the number of vertices $u\in U$ which are adjacent
in $B_0$ to both elements of the pair $e_i$, that is which are the tips of the tepees of $B_0$ over $e_i$.
Clearly, $N=|\{i:x_i>0\}|$. For convenience, assume that $x_i>0$ for $i=1,\dots,N$. Observe that for sufficiently
large $n$, by the choice of $\ep$,
$$\sum_{i=1}^Nx_i=\sum_{u\in U}{d_u\choose2}>\frac12(1-O(n^{-1/2}))(1-2\epsilon')^3|U|\pi^2|W'|^2>\frac13\pi^2|U||W'|^2\ , $$
and that $\sum_{i=1}^N{x_i\choose2}$ equals the number of copies of $C_4$ in $B_0$. Consequently,
$$\sum_{i=1}^N{x_i\choose2}\le2{|U|\choose2}{|W|\choose2}p^4<\frac12p^4|U|^2|W'|^2\ .$$

By the Cauchy-Schwarz inequality,
\begin{equation}\label{C-S}
\sum_{i=1}^N{x_i\choose2}\ge N{\sum x_i/N\choose2}=\frac{\sum x_i}2\left(\frac{\sum x_i}N-1\right)\ .
\end{equation}
Consider two cases:

I. $N\ge\sum x_i/2$. Then
$$N\ge\frac16|U|\pi^2|W'|^2\ge\frac1{3}\times\frac{\lambda n}{10}\times\frac{\lambda^2p^2}9{|W'|\choose2}\ge d{|W'|\choose2},$$ as required.

II. $N\le\sum x_i/2$. Then, by (\ref{C-S}),

$$N\ge\frac{(\sum x_i)^2}{4\sum{x_i\choose2}}\ge\frac{(1/9)\pi^4|U|^2|W'|^4}{2p^4|U|^2|W'|^2}\ge d{|W'|\choose2}\ ,$$ as needed.\qed

 \section{  Summary, Further Remarks, Glossary}
 \label{final}
 When embarking on the journey of solving this problem we
 did not realize the range of techniques which
 would be required, nor did we predict the extent of technical
 difficulties involved. Already one paper,
 \cite{FKRRT}, has sprung as a side effect, solving a problem which arose
 during our work, and addressing themes related to it.
 We hope in the future to follow the various trails that
 we have encountered here.

 There are several directions in which the results of this paper
 may be generalized or extended. We state below the corresponding open problems.
 It seems that in order to solve some of them one will need to develop a better
 understanding of sparse regular graphs and ``online" Ramsey theory as described
 in \cite{FKRRT}.

 The natural questions that come to mind are with regard to generalizing the main theorem of this paper to
 the case of coloring with more than two colors, or the case of Ramsey properties where the defining forbidden
 monochromatic graph is not a triangle but some other graph. We conjecture that an analog
 of our main theorem in this paper holds for all these cases, that is, that for all such natural Ramsey properties
 there exists a sharp threshold. These problems seem to be within the grasp of the techniques
 of this paper but present some serious technical difficulties. It is our hope that we will be able to
 overcome these difficulties in future work.

 It is also of interest to study thresholds of Ramsey properties of
 random sets of integers, such as having a monochromatic arithmetic
 progression in every bi-coloring of the set (see\cite{RR3,KLR,Rado}).

 A different question which seems at the moment quite hard is whether
 one can improve the result in this paper by establishing the exact
 threshold probability for the Ramsey property, or even proving that it exists:
 \begin{question}
 Does the function $\widehat{c}(n)$ defined in Theorem \ref{main} tend to a limit
 as $n$ tends to infinity? If so what is the exact value of the constant it tends to?
 \end{question}
 It is worthwhile noting that this natural question has not been answered in any of the cases
 where the existence of a sharp threshold has been established using the techniques from \cite{F}
 (i.e. in \cite{F}, \cite{AF} and \cite{FK}).
\subsection{Acknowledgments}
The authors wish to thank the anonymous referee for his/her
serious and professional effort, and the resulting suggestions.

The work on this paper was supported by the following grants:
\begin{itemize}
\item Ehud Friedgut:
Research supported in part by the Israel Science Foundation grant
no. 0329745, and by the Binational Science Fund grant 2002166.
\item Vojtech R\"odl: Research supported by NSF grant NSF Award
DMS - 0300529, and by the Binational Science Fund grant 2002166.
\item Andrzej Ruci\'{n}ski: Research supported by KBN grant 2
P03A 015 23; part of this research performed during several visits
at Emory University
\item Prasad Tetali: Research supported
in part by NSF grants DMS-0100298 and DMS-0401239.
\end{itemize}

 \subsection{Glossary}
  Any fan of 19'th century Russian literature is familiar with the
 frustration of reading a novel with a huge number of characters,
 each having two or three names and nicknames. The timid western
 reader often finds comfort in a short list added by the editors to
 the translation, a list of all characters, each accompanied by a
 reminder of their family background. We present here such a list
 that may help the reader of this paper, who may refer back to it
 while reading. We also define here the appropriate choices to
 ensure the needed relations between various constants. We end the
 glossary with a diagram which gives the dependency between the
 constants.

 {\center \large \bf Constants}
 \begin{itemize}
 \item  $c,C$ are absolute constants given by Theorem 1.2; we have $c=c_2=1/e$ and $C=C_2=10^4$; their role is to bound the threshold probability $p(n)$ scaled by $\sqrt n$ (see Theorem \ref{allweneed});
 \item $\xi>0$ and $0<\ga <1 $ are arbitrary constants at the outset of Theorem~\ref{allweneed};

 \item $\nu$ is the number of vertices of $M$, an arbitrary balanced graph with average degree 4 at the outset of Theorem \ref{allweneed};
 \item $q = \ceiling{ 10C^2\nu/\ga}$ is an upper bound on the number $\phi$ of edges in the star forest $S$, a prototype of constellations; used to define Properties (P6) and (P7) of~$\G$;
 \item $\ga_1$ -- no such thing;
 \item $\displaystyle\ga_2 = \frac{\ga}{(2\nu+q)^q}$ is an intermediate constant in Section \ref{TandC};
 \item $\displaystyle\ga_3 =\frac{\ga_2}{2(\nu+q)^{(\nu+q)}}$ -- there are at least $\ga_3n^\nu$ special constellations in each~$G\in\G$;
 \item $\displaystyle\ga_4 = \frac{(\nu+q)^{\nu+q}}{3C^{2q}}\ga_3=\frac{\ga}{6C^{2q}(2\nu+q)^q}$ is an intermediate constant in Section \ref{cores};
\item $D$ is a parameter bounding from above  the $p$-density of patches (large subgraphs) of $G$ (see Definition \ref{nodenspa}) which is a crucial assumption in all sparse regularity lemmas; in our application of these lemmas, owing to Property (P5) of $\G$, we take $D=2$;
\item $D^*$ is a parameter bounding from above  the $(S,p^*)$-density of $S$-patches of $F$ (see Definition \ref{nodensePa})  which is a crucial assumption in the Subgraph Regularity Lemma \ref{reglem}; in our application of this lemma we take
$D^*=2^q+1$ (see Lemma \ref{nodensepa});
 \item $\displaystyle\lambda =\frac1{160} \left(\frac{\alpha_4}{D^*q(\nu+q)} \right)^2 $ -- every core has size at least $\lam pn^2$;
 \item $a=a(\lam)$ is a constant determined by Lemma \ref{pro1}; both, $\lambda$ and $a$ are used to define Property (P7) of $\G$;
 \item  $\displaystyle\tau_0  =\frac{a^2 \xi^6c^6}{2(a\xi^3c^3 + 2 \xi^5 C^5)}$ is a value of $\tau$ with which we apply Theorem~\ref{allweneed};
 \item  $\tau>0$ is an arbitrary constant at the outset of Theorem  \ref{allweneed};

 \item $\eta>0$ is chosen so that $\displaystyle\eta \log(1/\eta) \leq \frac {\tau}{4C\nu}$; this parameter serves to define cores and the above bound guarantees that there are few of them;

 \item  $\displaystyle\gd = \min \left\{ 0.4 \eta^q ,\left(\frac{\alpha_4}{8D^*}\right)^2\right\}$,
is an input of the Subgraph Regularity Lemma \ref{reglem}; the final partition $\Pi_1$ resulting from  Lemma
\ref{reglem} with $H$ being any star forest with at most $q$ edges, is guaranteed to be $\delta$-regular;

 \item
 $\displaystyle\eps(\ell)= \min\left\{\ell^{-4q}, \left( \frac{\delta
 \alpha_4}{D^*2^{q+3}q (\nu+q)^{\nu+q-2}}\right)^{4q} ,q^{-4},2^{-6-6q}\right\},$
   is a function of variable $\ell$, at the outset of the Subgraph Regularity Lemma \ref{reglem};
 \item $T_0,L_0$ are output constants   of the Subgraph Regularity Lemma  \ref{reglem};
 \item $t,l$ are not really constants but vary with $G$; the final partition $\Pi_1$ resulting from the Subgraph Regularity Lemma  \ref{reglem} with $H$ being any star forest with at most $q$ edges, is a $(t,l)$-partition for some $t\le T_0$ and $l\le L_0$;
do not confuse $l$ with $\ell$ -- the argument of the function $\eps(\ell)$ defined above.
 \item $\eps =\eps(l)$ -- the final partition $\Pi_1$ resulting from the Subgraph Regularity Lemma  \ref{reglem} with $H$ being any star forest with at most $q$ edges, is guaranteed to be $(\eps,p)$-uniform; also, this $\eps$ satisfies the assumptions of Proposition \ref{stars} with $D=2$ and $k=q$;

\item $n_0$ is a lower bound on $n$, different at different places in the paper, in particular it is an output constant of the Subgraph Regularity Lemma  \ref{reglem};
\item $n_1$
 is the lower bound on $n$ promised in Theorem \ref{allweneed}; it is the maximum of all values of $n_0$  encountered throughout our proof, most notably the $n_0$ of Lemma \ref{reglem}, as well as of several implicit lower  bounds on $n$ hidden in our calculations;
\end{itemize}

{\center \large \bf  Graphs }
\begin{itemize}
\item  $\displaystyle\rho(H)=\frac{|E|}{|V|}$ is half of the average degree in a graph  $H=(V,E)$.
\item $\R$ is the family of all graphs such that every bi-coloring of their edges results in a monochromatic triangle; they are referred to as Ramsey graphs;
\item $\G$ is a family of graphs which depends on $c,C,\ga,\nu,p=p(n)$ and is defined in Definition \ref{defG}; it is crucial for our proof that $G(n,p)\in\G$ almost surely (see Theorem \ref{allweneed} and Lemma \ref{aas});
\item  $\T(\lambda,a)$ is a family of graphs $G=(V,E)$ such that for any subgraph $F$ of $G$ with at least $\lambda |E|$ edges, the set ${\rm Base}(F)$
contains at least $a |V|^3$ triangles (for the definition of ${\rm Base}(F)$ see pages \pageref{base} or
\pageref{base1}); we have $\G\subseteq\T(\lambda,a)$ with $\lambda$ as above;
\item $M$ is an arbitrary, balanced graph with $\nu$ vertices and $2\nu$ edges, at the outset of Theorem \ref{allweneed}; we put $V(M)=\{x_1\dots x_\nu\}$;
\item $G=(V,E)$ is an arbitrary member of $\G$ with more than $n_0$ vertices; it is assumed that $G\not\in\R$, but $G$ satisfies (\ref{assumption}); the whole proof boils down to showing that it also satisfies (\ref{conclusion});
\item $M_X$ is an ordered copy of $M$ with vertex set $X=\{v_1\dots v_\nu\}\subset V$ and such  that
$x_a$ is mapped onto $v_a$, for each $a=1\dots\nu$;
\item $M^*$ is a random copy  $M_X$, where $X$ is chosen uniformly over all $(n)_\nu$ possibilities;
\item $M'_X$ is $M_X$ with edges colored by a fixed coloring $\sigma'$;
\item $T(X)$ is the set of all tepees in $G$ formed over all pairs of vertices from $X$ which are edges of $M_X$ (see Definition \ref{tepee} for the definition of a tepee);
\item $\widehat X$ is the subgraph of $G$ formed by the edges of $T(X)$;
\item $\widehat M$ is a common isomorphic type of all $M_X$ with $X\in\X_2$; it has $\nu+\phi$ vertices and $2\phi$ edges, $1\le\phi\le q$, that is, it consists of $\phi$ tepees over~$M$;
\item $S(X)$, for $X\in\X_3$, called a special constellation, is a subgraph of $\widehat X$  consisting of the left legs of all tepees from $T(X)$ (see Section \ref{con});
 \item $\mathcal{S}$ is the set of all special constellations;

\item $S$ is the common isomorphism type of all $S(X)$ with $X\in\X_3$; $S$ is a forest of $\nu$ stars;
\item $\phi$ is the number of edges of the star forest $S$; the stars forming $S$ have $\phi_1\dots\phi_\nu$ edges, where
$\phi=\phi_1+\cdots+\phi_\nu$;
\item $\pi_0 = \{ V_1 \dots  V_\nu, W_{11} \dots W_{1\phi_1} \dots W_{\nu,1}
 \dots W_{\nu,\phi_{\nu}} \},$
 is an initial equipartition of $V$; for clarity, we further assume that all the sets are of equal size $m=n/(\nu+\phi)$;
\item $\cal C$ is the set of all constellations, that is copies of $S$ which are consistent with $\pi_0$;
\item $F_{ab}=G[V_a,W_{ab}]$, $a=1\dots\nu$, $b=1\dots\phi_a$, is the bipartite subgraph of $G$ consisting of all edges with one endpoint in $V_a$ and the other in $W_{ab}$;
\item $F=\bigcup_{ab}F_{ab}$ is the subgraph of $G$ consisting of all edges connecting  sets $V_a$ with $W_{ab}$; in Section \ref{regular},
 $F$  also stands for an arbitrary member of $\F(H,m)$;
\item $H$ always stands for an arbitrary fixed graph with $h$ vertices, mostly in Section \ref{regular}, where the
general form of the Subgraph Regularity Lemma \ref{reglem} is proved; it is later applied with $H=S$;
\item $\F(H,m)$ denotes the family of all spanning subgraphs of  $H^m$, where $H^m$ is the $h$-partite graph  obtained by replacing each vertex $x$ of $H$ by a set $V_x$ of $m$ vertices, and by replacing every edge $xy$ of $H$ by the complete bipartite graph $K_{m,m}$ spanned between $V_x$ and $V_y$;
\item $\#(H,F)$ is our notation for the number of copies of $H$ in $F$, sometimes it only counts copies which are consistent with a specified partition of~$V(F)$;
\item $\Pi_0$ is an initial partition of $F$ consisting of $\pi_0$ and the (trivial partitions of) subgraphs $F_{ab}$;
\item $\Pi_1$ is a final partition of $F$ resulting from the Subgraph Regularity Lemma \ref{reglem}; it consists of vertex equipartitions $V_a=V_{a}^{1}\cup\cdots\cup V_{a}^{t}$ and
 $W_{ab}=W_{ab}^{1}\cup\cdots\cup W_{ab}^{t}$ (under a simplifying assumption, all these sets are of equal size $m/t$), and subgraph partitions
 $F_{ab}^{i,j}=G[V_a^i, W_{ab}^j]=\bigcup_{k} F_{ab}^{i,j,k};$
 note that $F = \bigcup_{a,b,i,j,k} F_{ab}^{i,j,k}$;
\item $\Pi$ stands for an arbitrary (or current) $(t,l)$-partition  (see Definition \ref{tl});
\item $\Psi$ stands for a refinement of $\Pi$  (see Definition \ref{tl});
\item $P$ always stands for a polyad -- a special kind of a subgraph of $F$, always related to (consistent with) a partition $\Pi$ (see Definition \ref{polyad}); in Section \ref{cores} all polyads are consistent with the final partition $\Pi_1$;
\item $\calp_{\Pi}=\calp$ is the set of all polyads consistent with
 a partition $\Pi$.
\item $R$, generally, is a subgraph of $F$, often a subgraph of a given polyad $P$; in the proof of Lemma \ref{reglem},
$R\subset P$ represents  polyads which are consistent with a partition $\Psi$ which refines $\Pi$;

\end{itemize}

{\center  \large \bf Quantities depending on $n$ }
\begin{itemize}

\item $n> n_0$ is the number of vertices of $G$;
 \item $m= n/ (\nu+\phi)$ is the number of vertices in each part of the
 partition $\pi_0$ of the vertices of $G$ (we assume this is an integer);
\item $p=p(n)$ is a sequence of edge probabilities satisfying $c<p\sqrt n<C$; it is also used for an abstract scaling factor in the formula for $p$-density (see Section \ref{regular});
\item $e_G(U,W)$ is the number of edges of a graph $G$ with one endpoint in $U$ and the other in $W$; here $U$ and $W$ are disjoint subgraphs of $V$;
\item $d(B)=d(U,V)=\displaystyle\frac{e_G(U,V)}{|U||V|}$ is the density of the pair $(U,V)$; also called the density of the bipartite graph $B$, where $B=G[U,V]$;
\item $d_p(B)=d_p(U,V)=\displaystyle\frac{e_G(U,V)}{p|U||V|}$ is, for $0<p\le1$, the $p$-density of the pair $(U,V)$; also called the $p$-density of the bipartite graph $B$, where $B=G[U,V]$;
\item $p^*$ is an abstract normalizing factor for the $(H,p^*)$-density $d_R$ of subgraphs $R$ of $F$; for our application with $H=S$ we take $p^*=p^\phi$;
\item $\kappa=|\cee|$ is  the total number of constellations (consistent with $\pi_0$ copies of $S$) in $F$; we have
  $\kappa\sim m^{\nu+\phi} p^\phi$;
\item $c_R = |C(R)|$, where $C(R)$ is the set of constellations in $R$;
 \item $s_R = |S(R)|$, where $S(R)$ is the set of special constellations in $R$;
 \item $d_R = \displaystyle\frac{s_R}{p^* c_R}$ is the (normalized) $(H,p^*)$-density of
 $R$; for our application with $H=S$ we take $p^*=p^\phi$;

 \end{itemize}

\

{\center \large \bf The flowchart of constants}

\

For the convenience of the reader we have included, in the following
a chart indicating the interdependencies of various constants that
we have used in our proofs.


\begin{figure}
\begin{center}
    \input{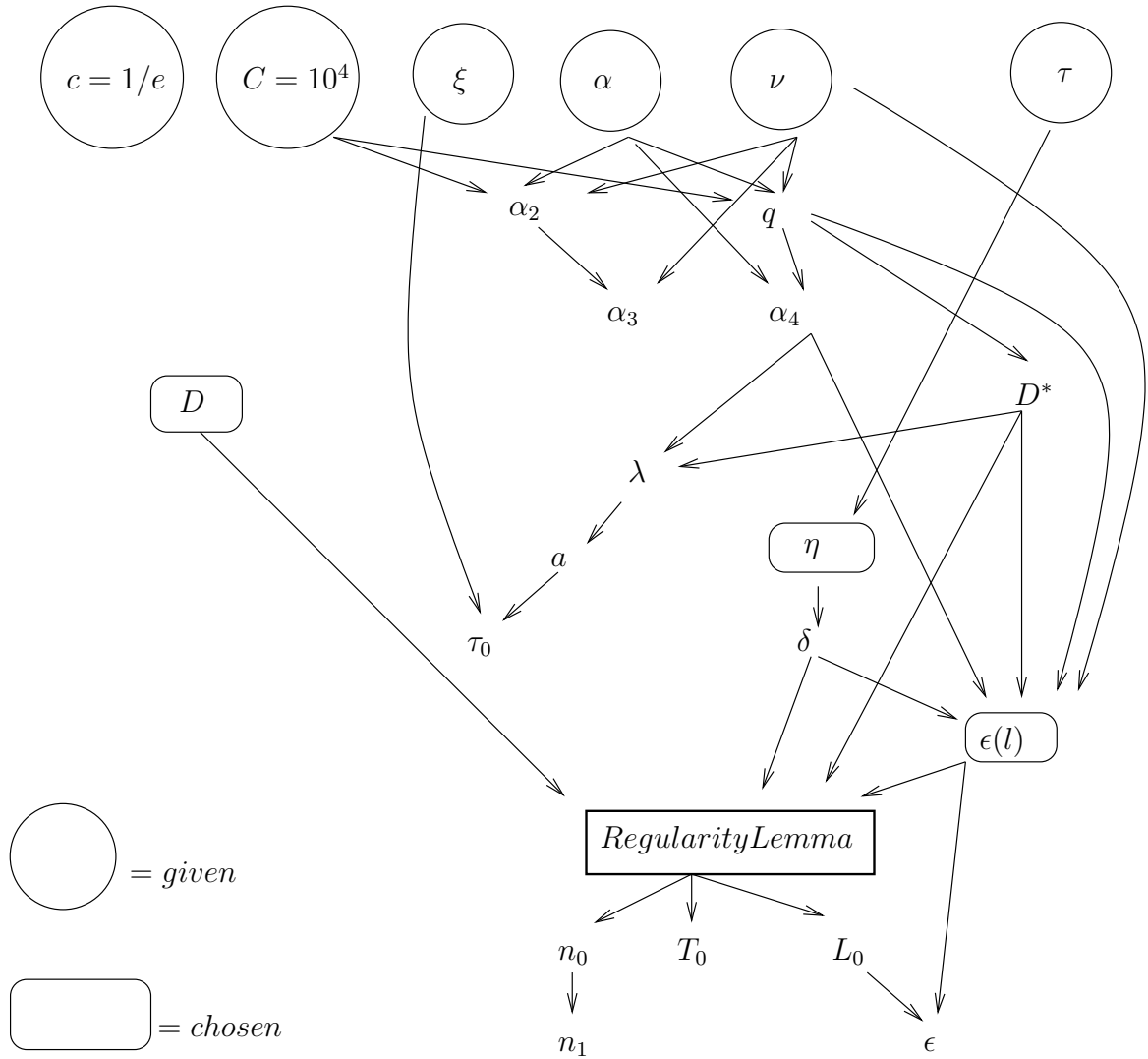}
      \caption{The flowchart of constants}
\end{center}
\end{figure}

\noindent Ehud Friedgut, Institute of Mathematics, Hebrew University,
Jerusalem, Israel.
Email: \texttt{ehudf@math.huji.ac.il.}
\\
\\ Vojtech R\"odl, Emory University, Atlanta, Ga. U.S.A.
Email:  \texttt{rodl@mathcs.emory.edu}
\\
\\ Andrzej Ruci\'{n}ski Adam Mickiewicz University, Pozna\'n, Poland.
Email: \texttt{rucinski@amu.edu.pl}
\\
\\ Prasad Tetali School of Mathematics and College of Computing
Georgia Institute of Technology Atlanta, GA 30332-0160, USA.
Email: \texttt{tetali@math.gatech.edu}

 \end{document}